\documentclass{acta_2022}
\usepackage{amssymb}
\usepackage{amsmath}
\usepackage{etoc}
\usepackage[longnamesfirst]{natbib} \usepackage{har2nat}
\setcitestyle{comma,aysep={}}
\shortcites{samplepreprint} 

\usepackage{amssymb,amsmath}
\usepackage{newtxtext}
\usepackage{newtxmath}
\DeclareSymbolFont{CMlargesymbols}{OMX}{cmex}{m}{n}
\DeclareMathDelimiter{(}{\mathopen} {operators}{"28}{CMlargesymbols}{"00}
\DeclareMathDelimiter{)}{\mathclose}{operators}{"29}{CMlargesymbols}{"01}
\DeclareMathAlphabet\mathcal{OMS}{cmsy}{m}{n}
\SetMathAlphabet\mathcal{bold}{OMS}{cmsy}{b}{n}
\pagerange{\pageref{firstpage}--\pageref{lastpage}}
\allowdisplaybreaks
\doi{XXXXXXXX}  
\usepackage[UKenglish]{babel} 
\newcommand{\ignore}[1]{}
\usepackage[twoside,
	paperheight=247mm,
	paperwidth=174mm,
	marginratio={3:5,2:3},
	left=18mm,
	top=20mm]{geometry}
\usepackage[hang]{subfigure}
\numberwithin{figure}{section}
\numberwithin{table}{section}
\usepackage{subfigure}
\usepackage{graphics,graphicx}
\usepackage[cmyk]{xcolor}
\usepackage{enumitem} 

\usepackage{amscd,MnSymbol}

\newcommand{\pp}[2]{\frac{\partial #1}{\partial #2}} 
\newcommand{\dede}[2]{\frac{\delta #1}{\delta #2}}
\newcommand{\dd}[2]{\frac{\diff#1}{\diff#2}}

\newcommand{\jump}[1]{\left[\!\!\left[ #1 \right]\!\!\right]}
\DeclareMathOperator{\diff}{d}
\DeclareMathOperator{\Ro}{Ro}
\DeclareMathOperator{\ddiv}{div}
\DeclareMathOperator{\DIV}{DIV}
\DeclareMathOperator{\ccurl}{curl}
\DeclareMathOperator{\Diff}{Diff}
\DeclareMathOperator{\ad}{ad}
\DeclareMathOperator{\vol}{vol}
\DeclareMathOperator{\tr}{tr}
\DeclareMathOperator{\Tr}{Tr}

\title[Compatible finite element methods for geophysical fluid dynamics]{Compatible finite element methods for geophysical fluid dynamics}

\author[C. J. Cotter]{Colin J. Cotter\\
Department of Mathematics, Imperial College London,\\ South Kensington Campus, London SW7 2AZ, \\ United Kingdom of Great Britain and Northern Ireland\\
E-mail: {colin.cotter@imperial.ac.uk}}

\begin{document}

\label{firstpage}
\maketitle

\begin{abstract}
  This article surveys research on the application of compatible
  finite element methods to large scale atmosphere and ocean
  simulation. Compatible finite element methods extend Arakawa's
  C-grid finite difference scheme to the finite element world. They
  are constructed from a discrete de Rham complex, which is a sequence
  of finite element spaces which are linked by the operators of
  differential calculus. The use of discrete de Rham complexes to
  solve partial differential equations is well established, but in
  this article we focus on the specifics of dynamical cores for
  simulating weather, oceans and climate.  The most important
  consequence of the discrete de Rham complex is the Hodge-Helmholtz
  decomposition, which has been used to exclude the possibility of
  several types of spurious oscillations from linear equations of
  geophysical flow.  This means that compatible finite element spaces
  provide a useful framework for building dynamical cores. In this
  article we introduce the main concepts of compatible finite element
  spaces, and discuss their wave propagation properties. We survey
  some methods for discretising the transport terms that arise in
  dynamical core equation systems, and provide some example
  discretisations, briefly discussing their iterative solution. Then
  we focus on the recent use of compatible finite element spaces in
  designing structure preserving methods, surveying variational
  discretisations, Poisson bracket discretisations, and consistent
  vorticity transport.
\end{abstract}

\setcounter{tocdepth}{1}
\tableofcontents 
\setcounter{tocdepth}{2}
\renewcommand{\contentsname}{Section contents}

\section{Introduction}
Atmosphere and ocean models used in weather forecasting and climate
simulation are built around dynamical cores, which predict the
fundamental quantities of fluid motion: fluid velocity, pressure,
density, temperature, and, in the case of the ocean,
salinity. Dynamical cores are computer implementations of numerical
discretisations of partial differential equation models of geophysical
fluid dynamics in the absence of viscosity. These models are then
coupled with physics parameterisations that describe additional
physics as well as fluid dynamical processes involving scales that
are too small and fast to represent explicitly in the dynamical core.
In the case of atmosphere models, this can include radiation
processes, cloud models, moisture and precipitation in various phases,
models of unresolved convection, boundary layers, and momentum
transfer due to unresolved internal waves. In the case of ocean
models, this includes vertical mixing due to convection, and
parameterisation of unresolved eddy motions. Atmosphere and ocean
models can be coupled together using parameterisations of air-sea
interaction processes, and with other process models to form climate
models. Other process models include models of land albedo, land ice,
sea ice, atmospheric chemistry etc. Additionally, in operational
forecasting systems, atmosphere and ocean models are blended with
observed data using data assimilation algorithms.

A very important aspect of global atmosphere and ocean models is the
wide range of timescales. In the rotating compressible Euler equations
that model atmospheric flow (although usually with additional
approximations), high frequency motions are possible in the form of
pressure waves and internal waves. However, in the solutions relevant
to atmosphere models these frequencies are observed to have very low
amplitude, with the velocity remaining very close to being horizontal
and divergence-free. Since the numerical solution is never anywhere
near being resolved in an atmosphere model, and due to the coupling
with all of the other modelled processes described above, plus the
modifications made to the solution by the data assimilation
algorithms, it is critically important that unphysical large amplitude
wave motions are avoided. Similar issues arise in the ocean model
context. All of this means that it is critically important to pay
attention to what happens in the dynamical core at the gridscale, to
avoid numerical error triggering these motions. These concerns lie at
the heart of the decisions about which discretisations to use in
building a dynamical core.

In the past decades, atmosphere and ocean dynamical cores have been
mainly built around finite difference, finite volume and
pseudospectral methods. In this article, we describe a more recent
approach to building dynamical cores using compatible finite element
methods. Finite element methods have the benefit that they do not
depend strongly on the structure of the underlying mesh for their
consistency and rate of convergence under mesh refinement. They can
use polygonal cells of different types, with triangulations in the
horizontal being particularly flexible. This allows mesh refinement in
regions of focus, and the constuction of meshes that conform to coastlines and
areas of high topography curvature. It also allows adaptive mesh
refinement although the advantages of this are less clear in large
scale weather and climate simulation.

The relaxation of mesh structure also allows for more uniform meshes
in spherical geometry. Instead of using a latitude longitude grid, which has
very thin cells at the poles due to the convergence of lines of
latitude, a cubed sphere mesh (quadrilateral refinement of a cube
mapped to a sphere) or an icosahedral mesh (triangular refinement of
an icosahedron mapped to a sphere) can be used. These pseudouniform
grids avoid the parallel scalability bottlenecks
\citep{zangl2015icon,adams2019lfric}. 

Finally, finite element methods allow for increasing the degree of
polynomials used in each cell, leading to higher order accuracy.
Higher order finite element methods have a denser, more structured
data layout in each cell, which can be exploited to try to achieve a
higher computational intensity, doing more computational work whilst
fetching data from nearby cells
\citep{dennis2012cam,giraldo2013implicit}.

Compatible finite elements, the subject of this survey, address the
problem of spurious numerical waves that cause problems when coupled
with numerical errors, physics parameterisations and data assimilation
schemes at the gridscale \citep{staniforth2012horizontal}. They can be
seen as an extension of the C grid finite difference method
\citep{arakawa1977computational}, which avoids many spurious wave
issues when used with quadrilateral cells (and slightly less so with
triangular cells). The C grid method achieves this by placing
different fields (velocity components, pressure, \emph{etc,}) on
different entities of the grid (cells, edges, vertices) so that the
discretised vector calculus operators (div, grad, curl) retain kernels
of the appropriate size. This was formalised as the Discrete Exterior
Calculus (DEC) \citep{hirani2003discrete}. Compatible finite element
methods do the same thing, but at the level of the finite element
spaces and the vector calculus mappings between them. These methods
have been unified in the Finite Element Exterior Calculus
\citep{arnold2006finite}\footnote{The author uses the term
``compatible finite elements'' as a way to discuss them with
practitioners when using the language of vector calculus instead of
differential forms, which are more unifying but require more
background material to discuss.}, with a long history going back to
the 1970s of applications to porous media, elasticity and fluid
dynamics. Their application to geophysical fluid dynamics also has a
long history, especially in the ocean
\citep{le1998finite,walters2005coastal}, and their representation of
exact geostrophic balance has been well known from numerical
dispersion analysis for some time
\citep{le2007analysis,rostand2008raviart}. However, it was
\citet{cotter2012mixed} that first noticed that it was the compatible
structure that was behind this property, and that this structure can
be used to understand the numerical dispersion properties of
compatible finite element schemes. Following earlier work by the
C grid and DEC communities, compatible finite element methods have
also been used to build structure preserving discretisations
that embed conservation laws at the discrete level. 

This article will introduce compatible finite element methods for
geophysical fluid dynamics and their properties, discuss how to build
atmosphere and ocean dynamical cores out of them, before focussing on
work on structure preserving discretisations. The rest of the article
is organised as follows. Section \ref{sec:spaces} will review
compatible finite element spaces and their fundamental
properties. Section \ref{sec:wave} will discuss their application to
linear wave models in geophysical fluid dynamics and the use of the
compatible structure to understand discrete wave propagation
properties. Moving towards nonlinear dynamical cores, Section
\ref{sec:transport} will discuss the discretisation of transport
schemes for compatible finite elements, and Section \ref{sec:examples}
will use them to describe discretisations for nonlinear geophysical
fluid dynamics models and survey their use.  One common theme is that
the analysis of the convergence and stability of these schemes is very
scarse and there is a lot of opportunity for future research in this
area. The next four sections concentrate on the use of compatible
finite element methods in structure preserving discretisations of
various types.  Section \ref{sec:variational} surveys methods
constructed from a discrete Hamilton's Principle.  Section
\ref{sec:poisson} covers methods constructed from almost Poisson
brackets that ensure conservation of energy and other quantities.
Section \ref{sec:consistent PV} discusses methods that embed a
diagnosed vorticity or potential vorticity field with its own
consistent dynamics.  Finally, Section \ref{sec:nonaffine} discusses
some issues common to all structure preserving methods related to
nonaffine meshes such as meshes approximating the sphere as well as
terrain following meshes.  We then end the survey with a very brief
summary in Section \ref{sec:summary}.

\section{Compatible finite element spaces}
\label{sec:spaces}
In this section we describe compatible finite element spaces as they
are used in geophysical fluid dynamics applications. Compatible finite
element spaces are sequences of spaces that form a discrete
differential complex (which we shall describe below). In applications
to geophysical fluid dynamics, the focus is on the de Rham complex.
There is a growing body of work for compatible finite elements spaces
built around the Stokes complex \citep[for
  example]{tai2006discrete,falk2013stokes,neilan2020stokes,hu2022family},
the elasticity complex, and the Regge complex for general relativity
\citep[for
  example]{christiansen2011linearization,christiansen2020discrete},
but we shall not discuss those spaces here.

\localtableofcontents

\subsection{Preliminary notation} First we establish some notation. Having defined the usual space
$L^2(\Omega)$ of square integrable scalar functions on $\Omega$ (and
writing $L^2(\Omega)^N$ as the space of vector functions on $\Omega$
in dimension $N$), for a domain $\Omega$ in three dimensions we have
\begin{align}
  H^1(\Omega) &= \left\{ \phi \in L^2(\Omega): \nabla\phi \in
    L^2(\Omega)\right\}, \\
    H(\ccurl;\Omega) &= \left\{ u \in
      L^2(\Omega)^3: \nabla\times u \in L^2(\Omega)^3 \right\},
      \\
      H(\ddiv;\Omega) &= \left\{ u \in L^2(\Omega)^3: \nabla\cdot u
        \in L^2(\Omega)\right\},
\end{align}
where $\nabla$, $\nabla\cdot$ and $\nabla\times$ are the appropriately
defined weak derivative operators (see a textbook on analysis of PDEs
for details, or proceed just assuming that these are the usual
derivatives for smooth functions and the cellwise derivatives for
finite element functions with the appropriate continuity, as we shall
discuss below). These spaces are accompanied by norms
defined as
\begin{align}
  \|\phi\|^2_{L^2(\Omega)} & = \int_\Omega \phi^2 \diff x, \\
  \|u\|^2_{L^2(\Omega)^N} & = \int_\Omega |u|^2 \diff x, \\
  \|\phi\|^2_{H^1(\Omega)} & = \|\phi\|^2_{L^2(\Omega)}
  + \|\nabla \phi\|^2_{L^2(\Omega)^N}, \\
  \|u\|^2_{H(\ccurl;\Omega)} & = \|u\|^2_{L^2(\Omega)^N}
  + \|\nabla\times u\|^2_{L^2(\Omega)^N}, \\
  \|u\|^2_{H(\ddiv;\Omega)} & = \|u\|^2_{L^2(\Omega)^N}
  + \|\nabla\cdot u\|^2_{L^2(\Omega)}.
\end{align}
In two dimensions, we similarly define
\begin{align}
  H(\ccurl;\Omega) &= \left\{ u \in
  L^2(\Omega)^2: \nabla^\perp\cdot u \in L^2(\Omega) \right\},
  \\
  H(\ddiv;\Omega) &= \left\{ u \in
  L^2(\Omega)^2: \nabla\cdot u \in L^2(\Omega) \right\},
\end{align}
having defined the operators $\nabla^\perp \phi = (-\partial_{x_2} \phi,
\partial_{x_1} \phi)$ and $\nabla^\perp\cdot u = -\partial_{x_2}u_1 +
\partial_{x_1}u_2$. Frequently we drop the $\Omega$ from this notation
when the meaning is clear. We will also use the $L^2$ inner product
notation,
\begin{align}
  \left\langle p, q \right\rangle &= \int_\Omega pq \diff x, \, \forall p,q\in L^2(\Omega), \\
  \left\langle u, v \right\rangle & = \int_\Omega u\cdot v \diff x, \, \forall u,v\in
  L^2(\Omega)^N, 
\end{align}
noting that
\begin{equation}
  \|p\|_{L^2(\Omega)}^2 = \left\langle p, p \right\rangle, \quad
  \|u\|_{L^2(\Omega)^N}^2 = \left\langle u, u\right\rangle.
\end{equation}
We also define $L^2(\Omega)$ inner products for
functions on boundaries,
\begin{align}
  \left\llangle p, q \right\rrangle &= \int_{\partial \Omega} pq \diff S, \, \forall p,q\in L^2(\partial\Omega), \\
  \left\llangle u, v \right\rrangle & = \int_{\partial \Omega} u\cdot v \diff S, \, \forall u,v\in
  L^2(\partial\Omega)^N.
\end{align}

\subsection{Discrete de Rham complexes} In three dimensions, discrete de Rham complexes on a domain
$\Omega$ consist of subspaces $\mathbb{W}_0\subset H_1(\Omega)$,
$\mathbb{W}_1\subset H(\ccurl;\Omega)$, $\mathbb{W}_2\subset
H(\ddiv;\Omega)$, $\mathbb{W}_3\subset L^2(\Omega)$, satisfying the
following commutation relations,
\begin{equation}
  \begin{CD}
    \mathbb{W}^0 = {H}^1 @> \diff^1 = \nabla  >>
    \mathbb{W}^1 = {H}(\ccurl) @> \diff^2 = \nabla\times >>
           \mathbb{W}^2 = {H}(\textrm{div}) @>
           \diff^3 = \nabla\cdot
  >> L^2 \\
  @VV{\pi_0}V @VV{\pi_1}V @VV{\pi_2}V @VV{\pi_3}V \\
  {\mathbb{W}}^0_h @> \diff^1 = \nabla  >> {\mathbb{W}}^1_h @>
  \diff^2 = \nabla\times
  >> \mathbb{W}^2_h @> \diff^3= \nabla\cdot >> \mathbb{W}^3_h, \\
\end{CD}
\end{equation}
where $\pi_i$, $i=0,1,2,3$ are surjective projections, satisfying
a bound appropriate to their domain, e.g.
\begin{equation}
  \|\pi_1u\|_{H(\ccurl)} \leq C\|u\|_{H(\ccurl)}, \, \forall
  u \in H(\ccurl).
\end{equation}
The commutation property means that $\diff^{k+1}\pi_k u = \pi_{k+1}\diff^{k+1}
u$, $k=0,1,2$. These projection operators (see
\citet{arnold2006finite} for a guide to their construction) do not
play a role in computations, but they ensure that the finite element
spaces are compatible in the sense that the differential operators
$\nabla$, $\nabla\times$ and $\nabla\cdot$ map surjectively onto the
kernel of the next operator in the sequence, just as is the case for
the infinite dimensional spaces at the top of the diagram.

In two dimensions, the de Rham complex is shorter, and there
are two possible ways to write it,
\begin{align}
  \label{eq:2d complex div}
  \begin{CD}
    \mathbb{V}^0 = {H}^1 @> \diff^1 = \nabla^{\perp}  >>
    \mathbb{V}^1 = {H}(\ddiv) @> \diff^2 = \nabla\cdot >>
           \mathbb{V}^2 =  L^2 \\
  @VV{\pi_0}V @VV{\pi_1}V @VV{\pi_2}V \\
  {\mathbb{V}}^0_h @> \diff^1 = \nabla^{\perp}  >> {\mathbb{V}}^1_h @>
  \diff^2 = \nabla\cdot
  >> \mathbb{V}^2_h, \\
  \end{CD} \\
  \label{eq:2d complex curl}
  \begin{CD}
    \tilde{\mathbb{V}}^0 = {H}^1 @> \diff^1 = \nabla  >>
    \tilde{\mathbb{V}}^1 = {H}(\ccurl) @> \diff^2 = \nabla^\perp\cdot >>
           \tilde{\mathbb{V}}^2 =  L^2 \\
  @VV{\pi_0}V @VV{\pi_1}V @VV{\pi_2}V \\
  \tilde{{\mathbb{V}}}^0_h @> \diff^1 = \nabla  >>
  \tilde{{\mathbb{V}}}^1_h @>
  \diff^2 = \nabla^\perp
  >> \tilde{\mathbb{V}}^2_h.
\end{CD}
\end{align}
The equivalence stems from the fact that any vector field in
$H(\ddiv;\Omega)$ in two dimensions can be transformed into a vector
field in $H(\ccurl;\Omega)$ by rotating it by $\pi/2$. In compatible
finite element methods for geophysical fluid dynamics we tend to use
\eqref{eq:2d complex div}, since it allows for local mass conservation
and exact application of flux boundary conditions.

\subsection{Discrete Hodge-Helmholtz decomposition}
Crucially for the geophysical fluid dynamics setting, the bounded,
commuting, surjective projections ensure a discrete version of the
Hodge-Helmholtz decompositions. At the infinite dimensional level,
these decompositions are
\begin{equation}
  \label{eq:hodge}
  \mathbb{W}^k = B^k \oplus \mathfrak{h}^k \oplus (B^*)^k,
\end{equation}
where
\begin{align}
  B^k & = \left\{ u \in \mathbb{W}^k: \exists \phi\in \mathbb{W}^{k-1},
  \, \mbox{with}\, u=\diff^k\phi\right\}, \\
  \mathfrak{h}^k & = \left\{ u \in \mathbb{W}^k: \diff^{k+1}u=0, \,
  \delta^ku=0\right\}, \\
  (B^*)^k & = \left\{ u \in \mathbb{W}^k: \exists \phi\in \mathbb{W}^{k+1},
  \, \mbox{with}\, u=\delta^k\phi\right\},
\end{align}
defining the dual operator $\delta^k:\mathbb{W}^k\to \mathbb{W}^{k-1}$
such that
\begin{equation}
  \left\langle \phi, \delta^k u \right\rangle = -\left\langle \diff^{k+1}\phi, u \right\rangle.
\end{equation}
When $u\in \mathbb{W}^k$ is appropriately constructed so that boundary
integrals vanish under integration by parts (and is sufficiently smooth
that integration by parts is well defined) then we can make the following
associations,
\begin{equation}
  \delta^0 = -\nabla\cdot, \, \delta^1 = \nabla\times, \,
  \delta^2 = -\nabla.
\end{equation}
For example, considering $\delta^1$, we have
\begin{equation}
  \left\langle v, \delta^1 u \right\rangle
  = \left\langle \nabla\times v, u \right\rangle
  = \left\langle v, \nabla\times u \right\rangle
  - \left\llangle v, n\times u \right\rrangle,
\end{equation}
where $n$ is the outward pointing normal to $\partial\Omega$.
If we choose $u\in \mathring{\mathbb{W}}^2$, where
\begin{equation}
  \mathring{\mathbb{W}}^2 = \left\{ u\in \mathring{\mathbb{W}}^2: u\cdot n = 0 \mbox{ on }
  \partial\Omega\right\},
\end{equation}
then we obtain that $\delta^1 u = \nabla\times u$ for
$u\in \mathring{\mathbb{W}}^2\cap \mathbb{W}^1$.
Similarly, we define
\begin{align}
  \mathring{\mathbb{W}}^0 & =
  \left\{ \phi \in \mathbb{W}^0: \phi=0 \mbox{ on }
    \partial\Omega\right\}, \\
    \mathring{\mathbb{W}}^1 & =
  \left\{ u \in \mathbb{W}^1: u\times n=0 \mbox{ on }
    \partial\Omega\right\}, \\
    \mathring{\mathbb{W}}^3 & =
\mathbb{W}^0.
\end{align}
See Chapter 3 of \citet{arnold2018finite} for an accessible
description of the full functional setting for these aspects.

Returning to \eqref{eq:hodge}, we explain the $\oplus$ notation.  It
means that any $u\in
\mathbb{W}^k$ can be written uniquely as $b + h + c$ with $b\in B^k$,
$h\in \mathfrak{h}^k$, and $c \in (B^*)^k$. In fact, further we have
that the spaces $B^k$, $\mathfrak{h}^k$ and $(B^*)^k$ are
mutually orthogonal under the $L^2$ inner product.
When $k=2$, we recognise
this as the Helmholtz decomposition for vector fields, which says that
a vector field can be written uniquely as $u=\nabla\times v + h +
\nabla \phi$, with $\nabla\cdot h=0$ and $\nabla\times h=0$; h is
referred to as a ``harmonic vector field''. This decomposition
is crucial to the understanding of rapidly rotating fluid dynamics.

As described in \citet{arnold2006finite}, three dimensional compatible
finite element spaces satisfy a discrete Hodge-Helmholtz
decomposition,
\begin{equation}
  \mathbb{W}^k_h = B^k_h \oplus \mathfrak{h}^k_h \oplus (B^*)^k_h,
\end{equation}
where
\begin{align}
    \mathfrak{h}^k_h & = \left\{ u \in \mathbb{W}^k_h: \diff^{k+1}u=0, \,
  \delta^k_hu=0\right\}, \\
  (B^*)^k_h & = \left\{ u \in \mathbb{W}^k_h: \exists \phi\in \mathbb{W}^{k+1}_h,
  \, \mbox{with}\, u=\delta^k_h\phi\right\},
\end{align}
and we have the discrete dual operator $\delta^k_h:\mathbb{W}^{k+1}_h
\to \mathbb{W}^k_h$, such that
\begin{equation}
  \left\langle \phi, \delta^k_h u \right\rangle = \left\langle \diff^{k+1}\phi, u \right\rangle,
  \quad \forall \phi \in \mathbb{W}^k_h, u\in \mathbb{W}^{k+1}_h.
\end{equation}
We note the asymmetry under discretisation: $\diff^k$ has the same
definition at the discrete level, but $\delta^k$ is replaced by the
approximation $\delta^k_h$.

The discrete Hodge-Helmholtz decomposition inherits some important
properties from the infinite dimensional decomposition.  First, we
have $B^k_h\subset B^k$. Second, although $\mathfrak{h}^k_h\neq
\mathfrak{h}^k$, we do have that
$\dim(\mathfrak{h}^k_h)=\dim(\mathfrak{h}^k)$. Further,
$\mathfrak{h}^k_h$ converges to $\mathfrak{h}^k$ as the mesh is
refined. These two properties have important consequences for the
correct representation of inertial oscillations in geophysical fluid
dynamics models. We do not have $(B^*)^k_h\subset (B^*)^k$.

For later discussion we also define
\begin{equation}
  \zeta^k = \left\{u\in \mathbb{W}^k: d^{k+1}u=0 \right\}, 
\end{equation}
noting that $\zeta^k=B^k\oplus \mathfrak{h}^k$. Similarly,
we write $\zeta^k_h=B^k_h\oplus \mathfrak{h}^k_h$.

The two dimensional
discrete de Rham complexes also provide an analogous Helmholtz-Hodge
decomposition for $\mathbb{V}_h^k$, which is most significant for
$\mathbb{V}_h^1$.

\subsection{Two dimensional compatible finite element spaces used in geophysical fluid dynamics}

For compatible finite element methods for geophysical fluid dynamics,
we are mainly focussed on the goal of producing three dimensional
fluid models of the atmosphere and ocean. However, model development
usually starts by consideration of the rotating shallow water
equations, since these equations encompass many of the challenges of
designing numerical schemes for atmosphere and ocean but without the
additional computational challenges of three dimensional models.  On
quadrilateral meshes (such as the ``cubed sphere'' meshing of the
sphere), the most commonly used spaces are $\mathbb{V}_h^0=Q_k$
(tensor product Lagrange elements, e.g. bilinear, biquadratic,
\emph{etc}), $\mathbb{V}_h^1=RT_{k-1}$ (Raviart-Thomas elements on
quadrilaterals, noting here that we use the traditional numbering
convention according to the largest complete polynomial space
contained by the finite element, not the UFL/FIAT numbering according
to the highest degree polynomial in the space), and
$\mathbb{V}_h^2=DQ_1$ (discontinuous tensor product Lagrange
elements).  However, it may be interesting to consider the trimmed
serendipity family of elements, especially at higher order
\citet{gillette2019trimmed}.  On triangular meshes (which are much
more flexible in allowing local mesh refinement), the main
possibilities are $\mathbb{V}_h^0=P_{k+1}$ (Lagrange elements),
$\mathbb{V}_h^1=BDM_k$ (Brezzi-Douglas-Marini elements on triangles),
and $\mathbb{V}_h^2=DG_{k-1}$, or $\mathbb{V}_h^0=P_{k+1}$,
$\mathbb{V}_h^1=RT_k$ (Raviart-Thomas elements on triangles), and
$\mathbb{V}_h^0=DG_k$.Whilst in most applications, the second grouping
is preferred, because it requires less degrees of freedom for the same
accuracy, in geophysical fluid dynamics this grouping has wave
propagation issues related to the Coriolis term, and we tend to prefer
the first grouping based on $BDM$ elements (we shall discuss this
later). In all of these discrete de Rham complexes, we see the
reduction in interelement continuity properties moving across the
discrete de Rham complex: $\mathbb{V}_h^0$ contains only continuous
functions, $\mathbb{V}_h^1$ requires continuity of normal components
its vector valued functions across cell facets (but not tangential
components), and $\mathbb{V}_h^2$ has functions with no interelement
continuity constraints. See \citet{boffi2013mixed} for definitions of
the finite elements introduced in this paragraph.

\subsection{Three dimensional compatible finite element spaces used in
  geophysical fluid dynamics}

When moving to three dimensional models, since the Earth's atmosphere
and ocean are much larger in horizontal extent than the vertical, good
preservation of hydrostatic balance (balance between gravitational
acceleration and vertical pressure gradient) requires the use of
prismatic meshes that are constructed by extruding a two dimensional
base mesh into layers. When the base mesh is constructed from
quadrilaterals, this produces hexahedra, and when the base mesh is
constructed from triangles, this produces triangular prisms. Three
dimensional discrete de Rham complexes are constructed on these spaces
by a tensor product of a two dimensional de Rham complex on the base
mesh with a one dimensional de Rham complex for the vertical
direction, given by
\begin{equation}
    \begin{CD}
    \mathbb{U}^0 = {H}^1 @> \diff^1 = \partial_x  >>
    \mathbb{V}^1 =  L^2 \\
  @VV{\pi_0}V @VV{\pi_1}V  \\
  {\mathbb{U}}^0_h @> \diff^1 = \partial_x  >> {\mathbb{U}}^1_h.
  \end{CD} \\
\end{equation}
Excluding splines \emph{etc.}, the main family of options for these
one dimensional spaces is continuous Lagrange elements of degree $k+1$
for $\mathbb{U}_h^0$ and discontinuous Lagrange elements of degree $k$
for $\mathbb{U}_h^1$. The 3D discrete de Rham complex is then formed
on the reference cell $\hat{K}_3=\hat{K}_2\times \hat{K}_1$ follows,
\begin{align}
  \mathbb{W}_h^0(\hat{K}_3) &= \mathbb{V}_h^0(\hat{K}_2) \otimes \mathbb{U}_h^0(\hat{K}_1), \\
  \mathbb{W}_h^1(\hat{K}_3) &= \underbrace{\hat{k}\mathbb{V}_h^0(\hat{K}_2)\otimes
      \mathbb{U}_h^1(\hat{K}_1)}_{=\mathbb{W}_h^{1,V}}
  \oplus
  \underbrace{
    \iota(\mathbb{V}_h^1(\hat{K}_2))^{\perp} \otimes \mathbb{U}_h^0(\hat{K}_1)}_{=\mathbb{W}_h^{1,H}}, \label{eq:W1 tensor}
  \\
  \mathbb{W}_h^2(\hat{K}_3) &=
  \underbrace{\mathbb{V}_h^2(\hat{K}_2) \otimes \hat{k}\mathbb{U}_h^0(\hat{K}_1) }_{=\mathbb{W}_h^{2,V}}
  \oplus
  \underbrace{\iota(\mathbb{V}_h^1(\hat{K}_2)) \otimes \mathbb{U}_h^1(\hat{K}_1)}_{=\mathbb{W}_h^{2,H}}, \label{eq:W2 tensor} \\
  \mathbb{W}_h^3(\hat{K}_3) &= \mathbb{V}_h^2(\hat{K}_2) \otimes \mathbb{U}_h^1(\hat{K}_1).
\end{align}
Here we write coordinates on $\hat{K}_3$ as $(x,z)$ with $x\in
\hat{K}_2$ and $z\in \hat{K}_1$, defining the tensor product
$V(K_2)\otimes U(K_1)$ as the span of function products $u(x)v(z)$
with $u\in V(K_2)$ and $v\in V(K_1)$. Further, $\hat{k}$ is the unit
upward pointing vector, $\iota$ is the inclusion operator that maps 2D
vectors into equivalent 3D vectors with zero vertical part, and
$\perp$ is the operator that rotates vectors by a quarter of a
rotation in the horizontal direction. These latter technicalities
involving $\iota$, $\hat{k}$ and $\perp$ can all be avoided in the
unified presentation of spaces of discrete differential forms in the
finite element exterior calculus, which also unifies many other aspects
across dimensions, numbering of the spaces, \emph{etc.} See 
\cite{arnold2014finite} for a full presentation of tensor product
discrete differential forms.

Again, in all of these three dimensional discrete de Rham complexes,
we see the reduction in interelement continuity properties moving
across the discrete de Rham complex: $\mathbb{W}_h^0$ contains only
continuous functions, $\mathbb{W}_h^1$ requires continuity of
tangential components its vector valued functions across cell facets
(but not normal components), $\mathbb{W}_h^2$ requires continuity of
normal components its vector valued functions across cell facets (but
not tangential components), and $\mathbb{W}_h^3$ has functions with no
interelement continuity constraints.

As indicated in (\ref{eq:W1 tensor}-\ref{eq:W2 tensor}). the
$\mathbb{W}_h^1$ and $\mathbb{W}_h^2$ spaces can be split into
vertical and horizontal parts, indicated by the $V$ and $H$ superscripts
respectively. As suggested by this notation, the vertical part
contains vector fields that point in the vertical direction, whilst
the horizontal part contains vector fields that point in the
horizontal direction. After Piola transformation to mesh cells (discussed in Section \ref{sec:piola}), this
decomposition is not preserved in general. However, if mesh cells are
arranged into a global tensor product mesh (i.e, flat vertical
layers), then this decomposition is preserved. This will also occur in
a spherical annulus meshed by radially extruding a 2D surface mesh of
the sphere (or an approximation of one). If a terrain following mesh
is used, so that the side walls of the mesh are arranged vertically,
but the horizontal layers move up and down to conform to mountain
ranges on the surface, the decomposition is only partially
preserved: $\mathbb{W}_h^{1,H}$ remains tangential to the ``up''
direction, and $\mathbb{W}_h^{2,V}$ remains normal to the ``up''
direction, but $\mathbb{W}_v^{1,V}$ will contain some horizontal
component and $\mathbb{W}_h^{2,V}$ will contain some vertical
component.

In discretisations for geophysical fluid dynamics it is important that
the gravity term can be balanced by a vertical pressure gradient.  To
avoid degeneracy in this balance, discussed in Section \ref{sec:hydros},
it is necessary to use a finite element space $\mathbb{W}_\theta$ for
temperature variables (entropy, temperature, potential temperature,
\emph{etc}.)  that is adapted to the vertical part of
$\mathbb{W}_h^2$, which is used to represent the velocity in this
framework. We choose $\mathbb{W}^\theta(\hat{K})$
so that $\theta \in \mathbb{W}^\theta(\hat{K}) \implies
\hat{k}\theta\in \mathbb{W}_h^{2,V}(\hat{K})$,
\emph{i.e}.,
\begin{equation}
  \mathbb{W}_\theta(\hat{K}) = \mathbb{V}_h^0(\hat{K}_2) \otimes
  \mathbb{U}_h^1(\hat{K}_1).
\end{equation}

It is also useful to construct vertical slice models by making
analogous tensor product constructions combining the
$(\mathbb{U}_h^0,\mathbb{U}_h^1)$ de Rham complex with itself.  Since
the only possible 1D meshes are intervals, this just leads to the
usual tensor product elements for quadrilaterals that we have already
discussed above.

For further description of the construction and efficient
implementation of these tensor product elements within an automated
system, see \citet{mcrae2016automated}. The use of extruded meshes
also has computational benefits that offset the additional
computational cost of using unstructured meshes in the horizontal. If
a semistructured data layout is used (i.e. a horizontal unstructured
index and a vertical structured one) then for a reasonable number of
vertical layers (20 is already enough in numerical experiments) the
lookups to find data in the unstructured grid data structure are
negligible compared to the computational work done on data loaded into
memory. Thus there is no significant performance penalty to using an
unstructured data structure in the horizontal
\citep{bercea2016structure}. This means that the benefits of the
flexibility and mesh invariance of the unstructured grid data
structure in the horizontal can be exploited in performant 3D
geophysical fluid models.

\subsection{Local to global mappings}
\label{sec:piola}
Whilst there has been work on e.g. H(div) elements on quadrilaterals
where the polynomials are defined directly on the mesh elements
\citep{arbogast2016two}, here we mostly restrict discussion to discrete
de Rham complexes that are constructed on reference elements and
mapped to mesh elements using Piola maps. In three dimensions,
this corresponds to the following set of relations,
\begin{align}
  \psi\in \mathbb{W}_h^0(K)\!:\,
  &\psi\circ g_K = \hat{\psi}, \,
  \hat{\psi}\in \mathbb{W}^0_h(\hat{K}), \\
  u\in \mathbb{W}_h^1(K)\!:\, &u\circ g_K = J^{-T}\hat{u}, \, \hat{u}\in \mathbb{W}^1_h(\hat{K}), \\
  w\in \mathbb{W}_h^2(K)\!:\, &u\circ g_K = J\hat{w}/\det J, \, \hat{w}\in \mathbb{W}^2_h(\hat{K}), \\
  \phi\in \mathbb{W}_h^3(K)\!:\, &\phi\circ g_K = \hat{\phi}/\det J, \, \hat{\phi}\in \mathbb{W}^3_h(\hat{K}),
\end{align}
where $g_K$ is the map between reference cell $\hat{K}$ and mesh cell
$K$, with derivative $J$, $J^{-T}$ means the inverse of the transpose of $J$,and $\circ$ indicates function composition
\emph{i.e.}  $f\circ g$ is another function with $(f\circ
g)(x)=f(g(x))$. This ensures that the discrete de Rham complex
property is preserved under the mapping from $\hat{K}$ to $K$. To
define $\mathbb{W}_\theta(\Omega)$, we just use straightforward
composition with the reference to cell map as is done for
$\mathbb{W}_h^0$. This leads to some differences between
$\mathbb{W}_h^\theta$ and the vertical compoent of
$\mathbb{W}_h^2(\Omega)$ when terrain following coordinates are used.

In two
dimensions when using \eqref{eq:2d complex div}, this becomes
\begin{align}
  \psi\in \mathbb{V}_h^0(K)\!:\,
  &\psi\circ g_K = \hat{\psi}, \,
  \hat{\psi}\in \mathbb{V}^0_h(\hat{K}), \\
  w\in \mathbb{V}_h^1(K)\!:\, &u\circ g_K = J\hat{w}/\det J, \, \hat{w}\in \mathbb{V}^1_h(\hat{K}), \\
  \phi\in \mathbb{V}_h^2(K)\!:\, &\phi\circ g_K = \hat{\phi}/\det J, \,
  \hat{\phi}\in \mathbb{V}^2_h(\hat{K}).
\end{align}

In geophysical fluid dynamics applications, being able to solve
equations on the surface of a sphere is important. In general, two
dimensional complexes can be extended to orientable manifolds embedded
in three dimensions by restricting vector fields to be tangential to
the manifold at each point. Then, for such a vector field $u$, we
define $u^\perp=k\times u$, $\nabla^\perp u = k\times \nabla u$ and
$\nabla^\perp\cdot u= k\cdot \nabla \times u$, where $\nabla$ is now
the projection of the gradient onto the tangent plane. In fact, these
operations can all be given intrinsic definitions on any two
dimensional manifold without reference to an external space
$\mathbb{R}^3$ containing the manifold, best expressed using the
language of differential forms.  However, we do not do this here (see
\citet{arnold2006finite} for intrinsic constructions using
differential forms).

The geometric factors in these formulae introduce complications.
These are related both to the approximation properties of the spaces and to
their computer implementation, the latter due to the resulting
nonpolynomial integrands in weak formulations. When $g_K$ is an affine
transformation, $J$ is constant on each cell, and no alterations to
the approximation properties arise. However, when $g_K$ is nonaffine,
$J$ is nonconstant. This means that the transformed basis functions may
not span the same polynomial spaces as they do on the reference
cell. \citep{arnold2005quadrilateral,falk2011hexahedral} showed that
this occurs for transformed H(div) and H(curl) elements on nonaffine
quadrilaterals and hexahedra. This interferes with the standard approximation
theory error estimates, because they apply the Bramble-Hilbert lemma
considering the largest polynomial space spanned by the basis. The
degradation of approximation theory was demonstrated in practice in
those papers. This presents a concern for the applicability of these
spaces in geophysical fluid dynamics, because we encounter nonaffine
transformations when using quadrilaterals or higher order triangular
cells (i.e. triangles that have been curved to better approximate the
sphere) to approximate the surface of the sphere. These spaces are
also used when extruding the sphere radially to make a spherical
annulus, required for three dimensional atmosphere and ocean models;
this leads to nonaffine prismatic meshes. Nonaffine cells also arise
when terrain following meshes are used. These are meshes that slope
layers up and down to conform to mountain ranges at the Earth surface
in the atmosphere and ocean. Fortunately, these approximation issues
can be avoided through the framework of \citet{holst2012geometric},
who used Strang type estimates to consider the ``variational crime''
of a sequence of meshes that only conform to a manifold in the limit.
Provided that the meshes can be obtained by piecewise smooth mappings
from an affine mesh (i.e. a mesh consisting of cells mapped to the
reference cell by affine transformations), then approximation error
bounds can be obtained that match those of the reference cell.  For
meshes of interest in geophysical fluid dynamics, this covers the
sphere meshes described above, as well as terrain following meshes
obtained by smooth transformation (smoothing of topography is standard
practice in atmosphere and ocean modelling). For the case of meshes (nonaffine or otherwise) of the sphere
extruded into a spherical annulus, \citet{natale2016compatible} showed
that these meshes can be obtained via transformation from an affine
mesh embedded in four dimensions.

Additionally, on unstructured meshes one must take care that the
degrees of freedom are correctly matched up on facets, and with the
correct sign (since they involve normal and tangential components in
general). A systematic approach for this was set out by
\citet{rognes2010efficient}, which is now implemented in a number of
finite element software systems, such as FEnICS
\citep{logg2012automated} and Firedrake
\citep{rathgeber2016firedrake}. Such systems are very useful as an aid
to productivity when dealing with the complexities of compatible
finite elements.

\subsection{Replacing $\nabla\cdot$ with $\DIV$}
One practical modification to this framework on nonaffine meshes is to
replace the transformation for the $L^2$ space ($V = \mathbb{V}^2_h$
in 2D
or $V=\mathbb{W}^3_h$ in 3D) with straightforward composition, i.e.
\begin{align}
  \phi\in {V}(K)\!:\, &\phi\circ g_K = \hat{\phi}, \,
  \hat{\phi}\in {V}(\hat{K}).
\end{align}
Then, the discrete de Rham complex can be restored by replacing
$\nabla\cdot$ with $\DIV:=P_V\nabla\cdot$, where $P_V$ is the $L^2$
projection into $V$. This is always a local operation since $V$ is a
discontinuous space with no interelement coupling.  The diagram
still commutes upon replacing $d^N$ (where $N$ is the dimension) with $\DIV$,
since if $u\in H(\ddiv; \Omega)$, then $\pi_N\circ \nabla\cdot=\DIV$.  This
idea was originally proposed in \citet{bochev2009rehabilitation}, who
defined $\DIV$ equivalently in the case of $RT_0$ spaces on quadrilateral
grids using the discrete exterior calculus framework by evaluating
fluxes through edges, after which the definition as projection of the
divergence can be obtained using the divergence theorem. The extension
to more general H(div) finite element spaces on nonaffine meshes was
probably clear to those authors, but it was also discussed from a
practical perspective in \citet{shipton2018higher}.

\subsection{Primal dual grids}
There have been various attempts to build complementary spaces on dual
grids, in order to avoid the global mass solves that are required to
compute $\delta$. These spaces are linked by ``discrete Hodge star''
operators $(\star_0,\star_1,\star_2)$, with mappings indicated in the
following diagram,
\begin{equation}
  \begin{CD}
  {\mathbb{V}}^0_h @> \diff^1 = \nabla^{\perp}  >> {\mathbb{V}}^1_h @>
  \diff^2 = \nabla\cdot
  >> \mathbb{V}^2_h, \\
  @VV{\star_0}V @VV{\star_1}V @VV{\star_2}V \\
  \tilde{\mathbb{V}}^2_h @< \diff^2 = \nabla^{\perp}\cdot  <<
  \tilde{\mathbb{V}}^1_h @<
  \diff^1 = \nabla
  << \tilde{\mathbb{V}}^0_h.
\end{CD}
\end{equation}
It is important that the discrete Hodge star maps are invertible.  The
approach is clearest for lowest order spaces, where there is at most
one nodal degree of freedom per edge, vertex or cell in the definition
of each space, in analogy with the discrete exterior calculus
\citep{hirani2003discrete}. \citet{thuburn2015primal,melvin2017wave}
constructed complementary spaces on triangulations and their duals by
subdividing cells into triangles and placing $RT_0$ elements in each
subcell. A constraint is applied so that the $d_2$ operators produce
constant functions over the primal and dual cells respectively. This
scheme produces an extension of the primal dual finite
difference C grid approach to the consistent finite element setting
(inconsistencies in the Coriolis term occur on dual icosahedral and
cubed sphere setups in the framework of
\citet{thuburn2012framework,thuburn2014mimetic}).  The framework of
mimetic spectral elements extends this idea to higher order polynomial
spaces \citep{lee2018discrete}.

\section{Wave propagation properties}
\label{sec:wave}
In this section we review the properties of compatible finite element
discretisations applied to linearised geophysical fluid dynamics.
We will see that the favourable properties of the C grid finite difference
approach to finite element methods. It is these properties that
underpin the Met Office's choice to use compatible finite element
methods to build their ``Gung Ho'' atmospheric dynamical core that
lies at the centre of their next generation LFRic modelling system,
for example \citep{adams2019lfric}.

\localtableofcontents

\subsection{Compatible discretisation of the linear rotating shallow
  water equations}
For now, we consider compatible finite element methods applied to the
linearised rotating shallow water equations on the doubly periodic
plane with constant Coriolis parameter,
\begin{equation}
  u_t + \underbrace{fu^{\perp}}_{\mbox{Coriolis}} + \underbrace{g\nabla \eta}_{\mbox{pressure gradient}} = 0, \,
  \eta_t + H\nabla\cdot u = 0,
\end{equation}
where $u$ is the horizontal velocity, $h=\eta+H$ is the layer height
with $H$ constant and $\eta$ integrating to zero, $f$ is the
(constant) Coriolis parameter, and $g$ is the acceleration due to
gravity. The $\perp$ operator rotates vectors one quarter of a
rotation, to obtain $u^{\perp}=(-u_2,u_1)$.  These equations are
solved in two dimensions. 

We use the Helmholtz decomposition $u=\nabla^\perp \psi + \bar{u} +
\nabla \phi$, where $\bar{u}$ is a spatially constant (but possibly
time dependent) vector field, as the harmonic vector fields in the
doubly periodic plane are of this form. Then, we may write
\begin{align}
  \label{eq:phi linear}
  \phi_t - f\psi + g \eta & = 0, \\
  \label{eq:psi linear}
   \psi_t + f\phi & = 0, \\
  \bar{u}_t + f\bar{u}^{\perp} & = 0, \\
  \label{eq:h linear}
  \eta_t + H\nabla^2\phi & = 0.
\end{align}
First, we observe that the harmonic component $\bar{u}$ is decoupled
and rotates at frequency $f$. These are called inertial oscillations.
Second, we can find steady solutions with $\phi=0$, $\bar{u}=0$,
$\psi=g\eta/f$.  These solutions correspond to states of ``geostrophic
balance'', where the velocity is divergence free and the pressure
gradient term cancels out the Coriolis term. The remaining solutions
are inertia-gravity waves with $\phi$ nonzero.  Applying time
derivatives to \eqref{eq:h linear} and \eqref{eq:phi linear}, we
eliminate $\psi$ and $\eta$ (by also using \eqref{eq:psi linear}, to
obtain
\begin{equation}
  \phi_{tt} + f^2 \phi - gH\nabla^2\phi = 0,
\end{equation}
which is the inertia-gravity wave equation. When $f=0$, this becomes
the wave equation with wavespeed $\sqrt{gH}$.  When $f\neq 0$ (a
positive sign would be used for Northern hemisphere dynamics) then the
equation becomes a Klein-Gordon equation, which is dispersive.

Numerically induced oscillations (physical and spurious) in shallow
water models have been extensively examined by Le Roux in a series of
papers using dispersion analysis
\citep{le2007analysis,roux2008analysis,rostand2008raviart,le2012spurious},
including for discretisations in the compatible finite element family.
Much of our vocabulary we use in this section has been taken from
that work.

To discretise this equation in space using compatible
finite elements, we pick a two dimensional discrete
de Rham complex $(\mathbb{V}_h^0,\mathbb{V}_h^1,\mathbb{V}_h^2)$,
and choose $u\in \mathbb{V}_h^1$, $\eta\in \mathbb{V}_h^2$.
Then the usual introduction of inner products with test
functions and integration by parts leads to the discrete
formulation,
\begin{align}
  \label{eq:fem linear u}
  \left\langle w, u_t \right\rangle + \left\langle w, f u^{\perp}\right\rangle
  - \left\langle \nabla\cdot w, g\eta \right\rangle & = 0, \quad \forall w\in
  \mathbb{V}_h^1, \\
  \label{eq:fem linear eta}
  \left\langle \alpha, \eta_t + H\nabla\cdot u \right\rangle & = 0,
  \quad \forall \alpha \in \mathbb{V}_h^2.
\end{align}
Since $\eta_t + H\nabla\cdot u \in \mathbb{V}_h^2$, we may choose
$\alpha=\eta_t + H\nabla\cdot u$ in \eqref{eq:fem linear eta},
concluding that
\begin{equation}
  \label{eq:eta_t L2}
  \eta_t+H\nabla\cdot u=0\mbox{ in }L^2(\Omega).
\end{equation}
This is a useful property that we will use later.
\subsection{Geostrophic balance}
In large scale atmosphere and ocean applications, it is very important
that discretisations can preserve states of geostrophic balance well;
this is usually tested on the linear rotating shallow water equations
as discussed here. \citet{cotter2012mixed} proved that
compatible finite element discretisations will have exactly steady
geostrophic balanced solutions. To be precise, for any divergence-free
velocity field $u$ with zero harmonic component $\bar{u}$, there
exists an $\eta$ such that $(u,\eta)$ form a steady solution in geostrophic
balance. To show this, we just find $\psi \in \mathbb{V}_h^0$ such
that $u=\nabla^\perp\psi$ (possible from the discrete Helmholtz
decomposition), and pick $\eta$ as the $L^2$ projection of $\psi$ into
$\mathbb{V}_2$ before multiplying by $g/f$. Then,
\begin{align}
  \left\langle \nabla\cdot w, g\eta \right\rangle =
  \left\langle \nabla \cdot w, f\psi \right\rangle
  = -\left\langle w, f\nabla \psi \right\rangle
  = \left\langle w, f u^\perp \right\rangle, \quad \forall w\in \mathbb{V}_h^1,
\end{align}
where the first equality holds from the projection, since
$\nabla\cdot w\in \mathbb{V}_h^2$. The second equality holds
by integration by parts, which is exact because $w\in H(\ddiv)$
and $\psi \in H^1$. The final equality follows from $u=\nabla^\perp\psi$.

What is not true is that for every $\eta$ there exists a $\psi$ giving
a steady state solution. This is because the $L^2$ projection from
$\mathbb{V}_h^0$ to $\mathbb{V}_h^2$ is not a bijection. However, this
does not hold in the C grid finite difference case either.

It is also important that discretisations correctly represent inertial
oscillations. In the linear rotating shallow water equations, the only
solutions with $\eta=0$ are the inertial oscillations with spatially
constant $u=\bar{u}$ rotating at frequency $f$. It is important that
discretisations are free of spurious additional inertial modes,
\emph{i.e.} solutions with $\eta=0$ but with spatially varying
$u$. \citet{le2012spurious} examined spurious modes in various finite
element discretisations, showing that when they are present they lead
to degraded error convergence rates. They have also been observed to
lead to problems in practical ocean model simulations, where they can
be excited by nonlinearity in baroclinic jets. This results in the
formation of spurious gridscale oscillatory patterns that do not
change the pressure/layer depth \citep{danilov2019geometric}. In a
closed bounded domain with boundary condition $u\cdot n=0$, we do not
expect inertial oscillations because the space of harmonic vector
fields only contains 0. However, it is possible for spurious inertial
oscillations to satisfy the boundary condition, leading to their
excitation.

\subsection{Inertial oscillations}
\citet{natale2016compatible} showed that compatible finite element
discretisations applied to the rotating shallow water equations in the
periodic plane have the following property: the only time-varying
solutions of (\ref{eq:fem linear u}-\ref{eq:fem linear eta}) with
$\eta=0$ have spatially constant $u_t$, corresponding to inertial
oscillations oscillating with frequency $f$. Any time independent
solutions are in the kernel of the discrete Coriolis operator,
i.e. $u\in \mathbb{V}_h^1$ such that
\begin{equation}
  \left\langle w, u^{\perp} \right\rangle = 0, \quad \forall w\in \mathbb{V}_h^1.
\end{equation}
In other words, these discretisations are free
from inertial oscillations. To see this, first note that if $\eta=0$,
then \eqref{eq:eta_t L2} implies that $\nabla\cdot u=0$, so that
$u=k+\nabla^\perp\psi$ for $k\in\mathfrak{h}_h$ and $\psi\in
\mathbb{V}_h^0$, from the discrete Helmholtz decomposition. Using
$w=\nabla^\perp\gamma$ in \eqref{eq:fem linear u}, we get
\begin{equation}
  \left\langle \nabla^\perp \gamma, u_t \right\rangle =
  -f\left\langle \nabla^\perp \gamma, u^\perp \right\rangle
  = f\left\langle \gamma, \nabla\cdot u \right\rangle = 0, \quad
  \forall \gamma \in \mathbb{V}_h^0.
\end{equation}
Hence, $u_t$ is orthogonal to $B_h^k$, \emph{i.e.}, $u_t\in
\mathfrak{h}_h$. We know from the discrete Hodge-Helmholtz
decomposition that $\dim(\mathfrak{h}_h)=\dim(\mathfrak{h})=2$ on the
periodic plane; $\mathfrak{h}$ are the constant vector fields. In
fact, since the constant vector fields are in $\mathbb{V}_h^1$,
and they are divergence free, and in the kernel of $\delta_h$.
Time independent divergence free solutions satisfy
\begin{equation}
  0 = \left\langle w, u_t \right\rangle = -f \left\langle w, u^\perp\right\rangle, \quad
  \forall w\in \mathbb{V}_h^1,
\end{equation}
\emph{i.e.} they are in the kernel of the discrete Coriolis
operator. Vector fields in this kernel are referred to as ``Coriolis
modes''. They are the main downside of compatible discretisations, but
the dimension of this kernel is always found to be finite and
resolution independent in analyses by Le Roux
\citet{rostand2008raviart}; in fact this number is typically very
small. These Coriolis modes are also always found in C grid finite
difference discretisations.

\subsection{Inertia gravity waves}
\citet{cotter2012mixed} also examined the discrete inertia-gravity
waves that correspond to solutions of the Klein Gordon equation above.
Writing $u=\nabla^\perp\psi+\delta \phi$ for $\psi\in \mathbb{V}_0$
and $\phi\in \mathbb{V}_2$ (having already discarded the harmonic
component since it decouples), and choosing both $w=\nabla^\perp\gamma$
for $\gamma\in \mathbb{V}_h^0$ and $w=\delta \alpha$ for $\alpha\in
\mathbb{V}_h^2$, we obtain
\begin{align}
  \left\langle \delta \alpha, \delta \phi_t \right\rangle - f\left\langle \delta \alpha,
  \nabla \psi \right\rangle - g\left\langle \delta \alpha, \delta \eta \right\rangle
  & = 0, \quad \forall \gamma\in \mathbb{V}_2, \\
  \left\langle \nabla \gamma, \nabla \psi_t \right\rangle + f\left\langle \nabla\gamma,
  \delta\phi\right\rangle  & = 0, \quad \forall \gamma \in \mathbb{V}_h^0, \\
  \left\langle \phi, \eta_t + H \nabla\cdot \delta \phi \right\rangle & = 0, \quad
  \forall \phi \in \mathbb{V}_h^2.
\end{align}
Using the definition of $\delta$, and the fact that the height
equation holds in $L^2$, we get
\begin{align}
  -\left\langle \alpha, \nabla\cdot\delta \phi_t \right\rangle + f\left\langle \alpha,
  \nabla^2 \psi \right\rangle + g\left\langle \alpha, \nabla\cdot\delta \eta \right\rangle
  & = 0, \quad \forall \gamma\in \mathbb{V}_2, \\
  \left\langle \nabla \gamma, \nabla \psi_t \right\rangle - f\left\langle \gamma,
  \nabla\cdot\delta\phi\right\rangle
  & = 0, \quad \forall \gamma \in \mathbb{V}_h^0, \\
  \eta_t + H\nabla\cdot \delta \phi = 0.
\end{align}
We can recognise the operator $\nabla\cdot\delta$ as the mixed approximation
$\tilde{\nabla}^2$ of the Laplacian defined by
\begin{align}
  \left\langle \alpha, \tilde{\nabla}^2\phi \right\rangle
  - \left\langle \alpha, \nabla\cdot \sigma \right\rangle & = 0, \quad \forall \alpha \in
  \mathbb{V}_h^2.
\end{align}
After restricting to $\overline{\mathbb{V}}_h^2$, defined by
\begin{equation}
  \overline{\mathbb{V}}_h^2 = \left\{ \phi \in \mathbb{V}_h^2: \int_\Omega
  \phi \diff x = 0 \right \},
\end{equation}
this discretisation is well known to be invertible, stable and
convergent. Hence, when $f=0$, we can deduce
\begin{equation}
  \phi_{tt} - gH\tilde{\nabla}^2\phi = 0,
\end{equation}
\emph{i.e.} the mixed approximation of the wave equation for $\phi\in
\overline{\mathbb{V}}_h^2$.

When $f\neq 0$, we have to introduce projection operators $P_0:\overline{\mathbb{V}}_h^2\to
\overline{\mathbb{V}}_h^0$, $P_2:\overline{\mathbb{V}}_h^0\to
\overline{\mathbb{V}}_h^2$, 
defined by 
\begin{align}
  \left\langle \nabla\gamma, \nabla P_0\phi \right\rangle &= \left\langle \nabla \gamma,
  \delta \phi \right\rangle, \quad \forall \gamma \in \overline{\mathbb{V}}_h^0, \\
  \left\langle \alpha, -\tilde{\nabla}^2P_2\psi \right\rangle :=
  \left\langle \delta\alpha, \delta P_2\psi \right\rangle &= \left\langle \delta \alpha,
  \nabla \phi \right\rangle, \quad \forall \gamma \in \overline{\mathbb{V}}_h^2.
\end{align}
$P_0$ is well posed since it just involves solving the usual Galerkin
discretisation of the Laplacian on $\overline{\mathbb{V}}_h^0$, whilst
$P_2$ is well posed since it involves solving the mixed discretisation
as we have already discussed.

Using $P_0$ and $P_2$, we get
\begin{align}
  \phi_t - f P_2 \psi + gh & = 0, \\
  \psi_t + f P_0 \phi & = 0, \\
  \eta_t + H\tilde{\nabla}^2 \phi &= 0,
\end{align}
from which we can deduce the discrete Klein-Gordon equation
\begin{equation}
  \phi_{tt} + f^2 P_2P_0\phi - gH \tilde{\nabla}^2 \phi = 0.
\end{equation}
The behaviour of the numerical dispersion relation depends on the
kernel of the composition $P_2P_0$, for which a lower bound is
obtained by considering the relative sizes of $\bar{\mathbb{V}}_h^0$
and $\bar{\mathbb{V}}_h^2$. For $RT_k$ on quadrilaterals, these two
spaces have the same dimension in the periodic domain, so this
projection is not too harmful. For $RT_k$ spaces on triangles,
$\dim(\bar{\mathbb{V}}_2) > \dim(\bar{\mathbb{V}}_0)$, so $P_2P_0$ is
not surjective. In the case of the C grid finite difference method on
triangles, a similar issue arises, causing high and low frequency
branches of the numerical dispersion relation to intertangle, leading
to numerical noise when $f$ is sufficiently large
\citep{danilov2010utility}. Further analysis is required to really pin
down these issues in the compatible finite element case.  For $BDM_k$
spaces on triangles, we have the opposite situation,
$\dim(\bar{\mathbb{V}}_2) < \dim(\bar{\mathbb{V}}_0)$, so there is at
least the chance for $P_2P_0$ to be surjective (although both
projection operators will have checkerboard modes in their kernel on
structured meshes). This suggests that $BDM$ spaces are more appropriate
for geophysical fluid dynamics using triangular meshes, but further
analysis of these issues is needed to make these statements more
precise.

\subsection{Spectral gaps and zero group velocity}
\label{ssec:spectral gap}
\citet{staniforth2013analysis,melvin2014two} examined the numerical
dispersion relation for the cases of $RT_0$ and $RT_1$ on quadrilaterals,
motivated by building discretisations on the sphere using a cubed
sphere grid. This type of dispersion analysis allows to focus on group
($\partial \omega/\partial k$, where $\omega$ is the frequency and $k$
is the wavenumber) and phase velocity ($\omega k/|k|^2$) for numerical
discretisations. For $RT_0$, the numerical dispersion relation is very
similar to the C grid finite difference numerical dispersion relation
on quadrilaterals, with no turning points for the group velocity
except at maximum wavenumbers. For $RT_1$, the numerical dispersion
relation has two roots for each wave number, corresponding to the
resolution of higher wavenumbers in the gridcell. When these are
properly interpreted, the group velocity is mostly well behaved except
for a jump in the dispersion relation at $\Delta x$ wavelengths where
the group velocity goes to zero before and after the jump. Remarkably,
the jump is repairable by modifying the coefficients in the mass
matrix in such a way that the convergence rate is not eroded. It is
not really the jump itself which is the problem, but the repair also
makes the group velocity become nonzero through a L'H\^opital's rule
type cancellation. This fix is independent of the value of $f$ and
$H$, so is useable in practice. In numerical experiments,
\citet{melvin2014two} showed that this modification leads to
propagation of a wave packet with $\Delta x$ wavelength which
otherwise stays in the same location, spuriously.

\citet{eldred2019quasi} introduced an alternative approach to avoiding
these dispersion relation spectral gaps, in which the nodal variables
are the same as $RT_0$ and $DG_0$ on quadrilaterals, but a higher order
polynomial expansion is constructed by using nodal variables from
surrounding cells. This can be seen as a form of spline, but does not
increase the degree of continuity, just the polynomial degree.  The
effect on the numerical dispersion relation is that there is only one
branch, so there cannot be jumps. The downside of this approach is
that the stencil of the operators is extended to more cells, and the
standard approach of finite element assembly becomes more complicated.
\citet{eldred2018dispersion,eldred2019dispersion} examined the
spectral gaps in $RT_k$ elements for larger k, and showed that this spline
approach also fixes the problem for higher k.

\subsection{Rossby waves}
Following the approach of \citet{thuburn2008numerical},
\citet{cotter2012mixed} also examined the Rossby wave propagation
properties of compatible finite element spaces. Rossby waves occur in
the situation where the Coriolis parameter $f$ is spatially varying
(as is the case on the sphere under the ``traditional approximation'',
where $f=2\Omega\cdot n$, with $\Omega$ being the rotational velocity
of the sphere, and $n$ being the normal to the sphere surface). To
perform the analysis, we consider solutions on an infinite plane, with
$f=f_0+ \beta y$, with $f_0$ and $\beta$ being constants. Then, if the
Rossby number $\Ro=U/fL$ is small (where $U$ is a typical velocity scale
and $L$ is a typical spatial scale), and also $\beta L/f_0=\mathcal{O}(\Ro)$,
then we neglect $u_t$ and $\beta y u^\perp$ in the velocity equation.
This gives the geostrophic balanced states $(u,h)=(u_g,\eta_g)$ satisfying
\begin{equation}
  f_0 u^{\perp}_g = -g\nabla \eta_g,
\end{equation}
so that $u_t=0$, and $\eta_t=-H\nabla\cdot u=H\nabla\cdot
(f_0\nabla^\perp \eta)=0$, as we saw previously. This is the main
reason why it is important for the discretisation to represent these
geostrophic balanced states. To obtain dynamics, we consider
$\mathcal{O}(\Ro)$ corrections to $(u,\eta)$, which we write as
$(u_{ag},\eta_{ag})$ (ageostrophic velocity and height). The equations
at the next order in $\mathcal{O}(\Ro)$ then give
\begin{align}
  (u_g)_t + f_0 u_{ag}^{\perp} + \beta y u_g^\perp
  + g\nabla \eta_{ag} = 0,
  (\eta_g)_t + H\nabla\cdot u_{ag}  = 0.
\end{align}
Applying $-\nabla^\perp\cdot$ to the first equation and using the second gives
\begin{equation}
  -\nabla^2\psi_t + \frac{f_0}{H}(\eta_g)_t - \beta u_g\cdot\hat{y} = 0,
\end{equation}
where we used that $\nabla^\perp u_g=0$, so $u_g=\nabla^\perp\psi$
(having already eliminating inertial oscillations which are fast),
and $\hat{y}$ is the unit vector in the $y$-direction. Finally this becomes
\begin{equation}
  (\frac{f_0^2}{gH}-\nabla^2)
  \psi_t - \beta \pp{\psi}{x} = 0,
\end{equation}
which is the Rossby wave equation, which exhibits waves propagating
Westwards when $\beta>0$ (\emph{i.e.}, in the Northern Hemisphere).

Now we examine what happens with compatible finite element
discretisations in the low $\Ro$ limit. At leading order in Rossby
number, we obtain the geostrophic balance equation,
\begin{align}
  \left\langle w, fu^\perp_g \right\rangle - g\left\langle \nabla\cdot w, \eta_g \right\rangle
  & = 0, \quad \forall w\in \mathbb{V}_h^1, \\
  \left\langle \alpha, H\nabla\cdot u_g\right\rangle & = 0, \quad \forall \mathbb{V}_h^2.
\end{align}
We have already seen that this has solutions $u_g=\nabla^\perp\psi$,
$h = P_2(f\psi)/h$ for $\psi\in \mathbb{V}_h^0$. At the next order we have
\begin{align}
  \label{eq:discrete Rossby u_g}
  \left\langle w, (u_g)_t \right\rangle + \left\langle w, f_0 u_{ag}^\perp + \beta y u_g \right\rangle - \left\langle \nabla\cdot w,
  g\eta_{ag} \right\rangle & = 0, \quad \forall w \in \mathbb{V}_h^2, \\
  \left\langle \alpha, (\eta_g)_t + H\nabla\cdot (u_g) \right\rangle & = 0, \quad
  \forall \alpha \in \mathbb{V}_h^1.
  \label{eq:discrete Rossby eta_g}
\end{align}
Choosing $w=\nabla^\perp \gamma$ with $\gamma\in \mathbb{V}_h^0$, and
integrating by parts in \eqref{eq:discrete Rossby u_g} (permissible
since $\gamma\in H^1$ and $u_{ag}\in H(\ddiv)$) gives
\begin{equation}
  \label{eq:curl discrete Rossby u_g}
  \left\langle \nabla \gamma, \nabla \psi \right\rangle - \left\langle \gamma, f\nabla\cdot u_{ag}\right\rangle
  - \left\langle \gamma, \beta \underbrace{u_g\cdot \hat{y}}_{=\pp{\psi}{x}} \right\rangle 
   = 0, \quad \forall \gamma \in \mathbb{V}_h^0.
\end{equation}
\eqref{eq:discrete Rossby eta_g} implies that
\begin{equation}
  \label{eq: discrete Rossby eta_g L^2}
  f_0P_2(\psi)/g + H\nabla\cdot u_g = 0, \mbox{ in } L^2.
\end{equation}
Combining this with \eqref{eq:curl discrete Rossby u_g} then gives
the discrete Rossby wave equation,
\begin{equation}
  \label{eq:discrete Rossby psi}
  \left\langle \nabla \gamma, \nabla \psi_t \right\rangle + \left\langle \gamma, \frac{f_0^2}{gH}
  P_2\psi\right\rangle
  - \left\langle \gamma, \beta \pp{\psi}{x} \right\rangle 
   = 0, \quad \forall \gamma \in \mathbb{V}_h^0.  
\end{equation}
Without the $P_2$ projection operator, this would just be a regular $H^1$ finite element
approximation of the Rossby wave equation. With it, there is the possibility of some
projection errors altering the numerical dispersion relation, especially in the case
of $BDM$ elements where $\dim(\mathbb{V}_h^0) > \dim(\mathbb{V}_h^2)$. However,
this will only occur for high wavenumber waves, where the dynamics is already dominated
by Laplacian term, and the Rossby waves will have very slow phase and group velocities
in either case. Hence, we do not believe that this causes a problem for discrete
Rossby wave propagation. A similar argument was made in \citet{thuburn2008numerical}
when considering C grid discretisations on hexagons.

\citet{rostand2008raviart} examined the wave propagation properties of
the $RT_0$-$DG_0$ and $BDM_1$-$DG_0$ compatible finite element discretisations
(and $RT_0$-$CG_1$ and $BDM_1$-$CG_1$ discretisations that are not compatible)
using discrete dispersion relations computed through Fourier analysis.
These dispersion relations revealed the steady geostrophic modes for
constant $f$, and two branches of the dispersion relation for
$BDM_1$-$DG_0$, a primary one attached to the origin (zero frequency for
zero wavenumber) and a secondary one, which they described as
spurious. It is possible that this second branch can be interpreted as
corresponding to higher wavenumbers resolved in the cell, just as for
$RT_1$ on quadrilaterals as explored by \citet{staniforth2013analysis},
but doing these calculations is difficult on triangles and more work
is needed to clarify this. In their dispersion analysis,
\citet{rostand2008raviart} identified ``CD modes'' in the $BDM_1$-$DG_0$
discretisation, which are modes in the intersection of the Coriolis
operator and the divergence operator. These precisely correspond to
the modes in the kernel of the $P_2$ operator appearing in
\eqref{eq:discrete Rossby psi}. In experiments with the linear
rotating shallow water equations with balanced initial data and
$f=f_0+\beta y$, designed to examine Rossby wave propagation, they
observed accurate solutions on structured grids but very noisy
solutions on unstructured grids. The noise was attributed to
interactions with the CD modes. However, when we have repeated these
experiments using modern automated finite element systems such as
FEniCS and Firedrake, we have not observed this noisy behaviour with
$BDM_1$ on unstructured grids. It seems likely that
\citet{rostand2008raviart} had bugs in their implementation related to
the identification of the two nodal variables between two cells across
each edge (since no such problem arose with $RT_0$, which only has one
nodal variable per edge). This is understandable, because the problem
of how to systematically assemble $BDM_1$ and higher order H(div) spaces
on triangles was not solved until \citet{rognes2010efficient}; these
things are very difficult to implement by hand. These results may have
discouraged the adoption of the $BDM$ family for atmosphere and ocean
modelling, but it seems like a good option for lowest order spaces to
avoid the issues with spurious inertia-gravity wave propagation with
$RT_0$.

\subsection{Consistent linear tidal response}
\citet{cotter2016mixed} considered solutions of the linearised barotropic
tide equations, which take the form,
\begin{equation}
  u_t + fu^{\perp} + g\nabla (D+b) = -cu + F, \qquad
  D_t + H\nabla\cdot u = 0.
\end{equation}
These are the linear rotating shallow water equations with additional
topography $b$, spatiotemporal lunar forcing $F$ and linear friction
coefficient $c$. They used the Helmholtz equation to obtain an
exponentially damping lower and upper bound on the energy in the
absence of forcing. In the presence of time dependent forcing
(quasiperiodic forcing is appropriate for tidal models), they proved
that the solution converges at exponential rate to a time dependent
solution as $t\to \infty$, independent of the initial condition. This
is the solution that is of interest when predicting tides. When
compatible finite element methods are used to discretise the tide
equations, they proved that the continuous time discrete space
solution also converges exponentially to a time dependent numerical
solution, independent of the initial condition. Finally they showed
that this discrete attracting solution converges to the unapproximated
attracting solution as the mesh is refined. \citet{cotter2018mixed}
extended this analysis to a nonlinear model with $cu$ replaced by
$c|u|u$, which is the more realistic damping model that is actually
used by oceanographers. This nonlinear case is surprisingly subtle but
they were able to prove long time stability of the system and obtain
rates of damping in the unforced case. These were used to prove error
estimates for the discrete solution obtained using compatible finite
element methods. \citet{kirby2021preconditioning} used the compatible
finite element framework to design a preconditioner for the implicit
solver for the tidal equations, proving that the convergence rates
are independent of mesh resolution. \citet{cotter2022weighted} extended
this approach to the multiple layer version of this model.

\subsection{Hydrostatic balance}
\label{sec:hydros}
Finally in this section, we discuss the discrete hydrostatic balance
properties of compatible finite element methods.
Later, we shall introduce three dimensional geophysical fluid dynamics
models with gravity and pressure gradient terms, so that the velocity
equation takes the form
\begin{equation}
  \pp{u}{t} + \ldots + \underbrace{\nabla p}_{\mbox{pressure gradient}}
  = \underbrace{b\hat{k}}_{\mbox{gravity}},
\end{equation}
in the case of the Boussinesq equations (typically used in ocean modelling),
where $p$ is the pressure, $b$ is the buoyancy, and $\hat{k}$ is
the unit normal vector in the ``up'' direction. In the case of the
compressible Euler equations (typically using in atmosphere modelling),
we have
\begin{equation}
  \pp{u}{t} + \ldots + c_p\underbrace{\theta\nabla \Pi}_{\mbox{pressure gradient}} =
  \underbrace{-g\hat{k}}_{\mbox{gravity term}},
\end{equation}
where $\theta$ is the potential temperature, $\Pi$ is the Exner
pressure, and $g$ is the acceleration due to gravity. In both cases,
we are concerned with states of hydrostatic balance, which is when the
vertical component of the pressure gradient term balances the gravity
term. Either we are considering hydrostatic models, where this balance
is enforced exactly in the model, or we are considering nonhydrostatic
models where it is important that these hydrostatic states can be
accurately represented. In the compatible finite element case, to
study the vertical part of the velocity equation, we restrict the test
function in the velocity equation to $\mathbb{W}_h^{2,V}$ (assuming a
tensor product discrete de Rham complex), the vertical part of the
space $\mathbb{W}_h^2$ containing the discretised velocity $u$.

In the case of the Boussinesq equations, the discrete hydrostatic
balance is written (after integration by parts) as
\begin{equation}
  \label{eq:hydros}
  -\left\langle \nabla\cdot w, p \right\rangle = \left\langle w\cdot \hat{k}, b\right\rangle,
  \quad \forall w \in \mathring{\mathbb{W}}_h^{2,V},
\end{equation}
for pressure $p \in \mathbb{W}_h^3$ and buoyancy $b\in
\mathbb{W}_h^\theta$, where $\mathring{\mathbb{W}}_h^{2,V}$ is the
subspace of the vertical space $\mathbb{W}_h^{2,V}$ requiring the
boundary condition $u\cdot n=0$ at the top and the bottom.  Despite
appearances, this is actually only defining the vertical part of the
pressure gradient term, since $w\in \mathring{\mathbb{W}}_h^{2,V}$
always points in the vertical direction.

Given $p$, there is a unique $b$ that satisfies this hydrostatic
balance. To see this, we note that if the layers of the mesh are flat,
then $w\cdot k \in \mathbb{W}_h^\theta$ for all $w\in
\mathbb{W}_h^{2,V}$.  If the layers are not flat, \emph{i.e.} for
terrain following coordinates, then there exists $0 \leq \kappa \leq
\infty$ such that $\kappa w\cdot k\in \mathbb{W}_h^\theta$ for
all $w\in \mathbb{W}_h^{2,V}$. After replacing $w\cdot k=\kappa^{-1}\gamma$
for $\gamma \in \mathbb{W}_h^\theta$, we recognise the right hand side
of \eqref{eq:hydros} as a nondegenerate weighted $L^2$ inner product,
hence $b$ is unique.

To discuss the nature of the uniqueness of $p$, we consider an
alternative boundary condition with $u\cdot n = 0$ on the bottom, and
$p=p_0$ on the top (for some chosen $p_0$ which may depend on the
horizontal coordinate).  The equation after integration by parts and
use of the top boundary condition gives
\begin{equation}
  \label{eq:hydros bc}
  -\left\langle \nabla\cdot w, p \right\rangle
  + \left\llangle w\cdot n, p_0 \right\rrangle
  = \left\langle w\cdot \hat{k}, b\right\rangle,
  \quad \forall w \in \mathring{\mathbb{W}}_h^{2,V},
\end{equation}
where $\mathring{\mathbb{W}}_h^{2,V}$ is now the subspace with
vanishing normal component on the bottom only. To analyse this
problem, \citet{natale2016compatible} introduced the following
formulation, defining $(v,p)\in \mathbb{W}_h^{2,V}\times \mathbb{W}_h^3$
such that
\begin{align}
  \left\langle w, v\right\rangle -\left\langle \nabla\cdot w, p \right\rangle
  &= \left\langle w\cdot \hat{k}, b\right\rangle- \left\llangle w\cdot n, p_0 \right\rrangle,
  \quad \forall w \in \mathring{\mathbb{W}}_h^{2,V}, \\
  \left\langle \phi, \nabla\cdot v \right\rangle & = 0, \qquad \forall \phi
  \in \mathbb{W}_h^3,
\end{align}
which we recognise as a mixed problem defined on
$\mathbb{W}_h^{2,V}\times \mathbb{W}_h^3$. Using the same arguments as
we have previously, at the solution we have $\nabla\cdot v=0$ in
$L^2$. Since $v\cdot n=0$ on the bottom surface, and $v$ points in the
vertical direction, we conclude that $v=0$, and therefore $p$ solves
\eqref{eq:hydros bc}. \citet{natale2016compatible} showed that this
type of vertical mixed problem has a unique solution $(v,p)$. Hence,
there is a one to one correspondence between $p$ and $b$, as required.
If there are boundary conditions $u\cdot n=0$ on both top and bottom
surfaces, $p$ is only determined up to the value $p_0$ restricted
to the upper surface, also as required.

In the case of the compressible Euler equations, taking boundary
conditions $u\cdot n=0$ on the bottom surface and $\Pi=\Pi_0$
on the top surface, the discrete hydrostatic
balance is written (after integrating by parts) as
\begin{equation}
  -\left\langle \nabla\cdot (\theta w), \Pi \right\rangle
  = \left\langle w\cdot \hat{k}, g\right\rangle
  - \left\llangle w\cdot n, \Pi_0 \right\rrangle,
  \quad \forall w \in \mathring{\mathbb{W}}_h^{2,V},
\end{equation}
for $\Pi\in\mathbb{W}_h^3$ and $\theta\in \mathbb{W}_h^\theta$.  Using
an extension of the techniques described for the Boussinesq equation,
\citet{natale2016compatible} proved similar results. This motivates
the use of $\mathbb{W}_h^\theta$ for temperature variables like $b$ or
$\theta$. \citet{melvin2018choice} showed through linear
dispersion analysis applied to the compressible Boussinesq equations
that this choice does indeed lead to an absence of spurious
hydrostatic modes that would appear if $\theta \in \mathbb{W}_h^3$ or
$\mathbb{W}_h^0$.

\section{Transport and stabilisation}
\label{sec:transport}
Hopefully it is clear from Section \ref{sec:wave} that it might be
interesting to consider designing a numerical atmosphere or ocean
model using compatible finite element methods. One very important
aspect of these models is the choice of transport schemes, \emph{i.e.}
the discretisation of the advection operators. Since different fields
are restricted to different spaces from the discrete de Rham complex
(or $\mathbb{W}_h^\theta$) with different continuity constraints, we
need to consider a diverse range of transport schemes, some
of which we briefly survey in this section.

\localtableofcontents

\subsection{Transport of $H^1$ fields}
For scalar fields in $\mathbb{W}_0$ or $\mathbb{V}_0$ (we will call it
$V$ here), we consider the discretisation of the advection equation
\begin{equation}
  \pp{q}{t} + u\cdot \nabla q = 0,
\end{equation}
for some specified $u\in \mathbb{W}_h^2$ or $\mathbb{V}_h^1$, where we
assume that $u\cdot n=0$ on exterior boundaries. Since $V$ is a
continuous finite element space, we simply take the $L^2$ inner
product with a test function and integrate by parts to obtain the
standard continuous finite element approximation, seeking $q\in V$
such that
\begin{equation}
  \left\langle \phi, \pp{q}{t} \right\rangle - \left\langle \nabla \cdot (u\phi), q \right\rangle
  = 0, \quad \forall \phi \in V.
\end{equation}
As is well known, this discretisation tends to produce spurious
oscillations at regions of rapid changes in $q$. One way to suppress
these oscillations is to use the Streamline Upwind Petrov-Galerkin
(SUPG) method
\citep{brooks1982streamline,tezduyar1988petrov,tezduyar1989finite}.
In this approach, the test function $\phi$ is replaced by $\phi + \tau
u\cdot \nabla \phi$, which biases it in the upwind direction.  $\tau$
is some chosen stabilisation parameter which depends on the mesh, $u$
and other parameters. This leads to
\begin{equation}
  \left\langle \phi + \tau u\cdot \nabla \phi, \pp{q}{t} \right\rangle - \left\langle
  \nabla \cdot (u\phi), q \right\rangle
  + {\left\langle \tau u\cdot \nabla \phi, u\cdot \nabla q \right\rangle}
  = 0, \quad \forall \phi \in V.
\end{equation}
The final term performs diffusion along streamlines of $u$, which
tends to reduce spurious oscillations. By applying this modification
in the $\partial q/\partial t$ term as well as the $u\cdot \nabla q$
term, we obtain a consistent approximation (\emph{i.e.}, substituting
a smooth exact solution of the unapproximated equation produces zero).

Another possibility is the edge stabilisation approach proposed for
advection equations and analysed (when combined with diffusion) in
\citet{burman2004edge}, resulting in the formulation
\begin{equation}
  \left\langle \phi, \pp{q}{t} \right\rangle - \left\langle \nabla \cdot (u\phi), q \right\rangle
  + \left\llangle \gamma h^2 \jump{\nabla \phi},\jump{\nabla q}\right\rrangle_{\Gamma}
  = 0, \quad \forall \phi \in V,
\end{equation}
where
\begin{equation}
  \left\llangle u,v \right\rrangle_\Gamma = \int_\Gamma u\cdot v\diff S,
\end{equation}
$\Gamma$ is the union of all interior facets $f$ in the mesh
(\emph{i.e.}, facets joining two cells), $\jump{v}=v^+n^++v^-n^-$ for
vector fields $v$, each facet $f$ has been arbitrarily assigned $+$
and $-$ labels to its two sides, $a^{\pm}$ indicates the restriction
of the discontinuous function $a$ to the $\pm$ side of the facet,
respectively, $h$ is a mesh edge length parameter, and $\gamma$ is a
(possibly $u$ and or $q$ dependent) stabilisation parameter.  This
term has a diffusive effect across interior facets, penalising jumps
in $q$, without sacrificing consistency as $h\to 0$.

\subsection{Transport of $L^2$ fields}
The spaces $\mathbb{W}_h^3$ and $\mathbb{V}_h^2$ have no continuity
constraints, which allows for upwind stabilisation \emph{via} a
discontinuous Galerkin formulation. Here we consider the continuity
equation,
\begin{equation}
  D_t + \nabla\cdot(uD) = 0,
\end{equation}
where $u$ is as above. We introduce the discretisation by
multiplying by a test function and integrating over a single mesh cell
$e$,
\begin{equation}
  \int_e \phi D_t -\nabla \phi \cdot u D \diff x
  + \int_{\partial e} \tilde{D}u\cdot n \phi \diff S = 0,
  \quad \forall \phi \in V(e),
\end{equation}
where $\partial e$ is the boundary of $e$ with outward pointing normal
$n$, $V$ is the chosen discontinuous space, and $\tilde{D}$ is the
upwind value of $D$, which must be defined in terms of the values
of $D$ on the inside and the outside of $e$.  When $u\cdot n>0$,
then $\tilde{D}$ is equal to $D$ from inside $e$, otherwise the
outside value is used. If we sum over all the cells $e$ in the mesh,
we obtain
\begin{equation}
  \label{eq:DG}
  \left\langle \phi, D_t \right\rangle - \left\langle \nabla_h \phi, uD \right\rangle
  + \left\llangle \jump{\phi u}, \tilde{D} \right\rrangle_\Gamma
= 0, \quad \forall \phi \in V(\Omega),
\end{equation}
where $\nabla_h$ indicates the cellwise ``broken'' gradient,
\begin{equation}
  \nabla_h|_eD = \nabla|_e D,
\end{equation}
for each cell $e$ in the mesh. An analysis of how this choice of
upwind $\tilde{D}$ introduces stabilisation is provided in
\citet{brezzi2004discontinuous}. This method is locally conservative.

\subsection{Transport of $H(\ddiv)$ fields}
In consideration of the $(u\cdot \nabla)u$ term in the velocity
equation of geophysical models, discretisations for vector advection
equations of the form
\begin{equation}
  \pp{v}{t} + (u\cdot \nabla)v = 0,
\end{equation}
for a vector field $v$ in $H(\ddiv)$ spaces $\mathbb{W}_h^{2,V}$ or
$\mathbb{V}_h^1$, and $u$ is again as above. When the equation is
solved on surface of the sphere, $v$ is constrained to be tangential
to the sphere (or the mesh approximating the sphere in the discrete
case). Then, the equation becomes
\begin{equation}
  \pp{v}{t} + \mathbb{P}_S\left((u\cdot \nabla)v\right) = 0,
\end{equation}
where $\mathbb{P}_S$ is the Euclidean projection into the
tangent plane to the sphere.

Functions in $H(\ddiv)$ spaces are only partially continuous
(in the normal component across facets) so we need to start
by considering an upwind formulation on a single cell again,
\begin{equation}
  \int_e w\cdot v_t - \nabla\cdot (u\otimes w)\cdot v \diff x
  + \int_{\partial e} n\cdot u w\cdot\tilde{v} \diff S
  = 0, \quad \forall w \in V(e),
\end{equation}
where $(a\otimes b)_{ij}=a_ib_j$ for vectors $a$ and $b$, $\tilde{v}$
is the upwind value of $v$, $A:B=\sum_{ij}A_{ij}B_{ij}$ for two
matrices $A$ and $B$, and $V$ is the chosen $H(\ddiv)$ space.  Summing
up over all of the cells in the mesh gives
\begin{equation}
  \label{eq:H div advection}
  \int_\Omega w\cdot v_t - \nabla_h\cdot (u\otimes w)\cdot v \diff x
  + \int_{\partial \Gamma} \jump{u\otimes w}\cdot\tilde{v} \diff S
  = 0, \quad \forall w \in V(\Omega),
\end{equation}
where $\jump{u\otimes w}=(n^+\cdot u^+)w^+ + (n^-\cdot u^-)w^-$.
Since $v\in V$ has continuous normal components, $\tilde{v}$ only
differs from $v$ in the tangential component. Hence, the upwind
stabilisation may be insufficient to adequately suppress oscillations,
depending on the shape of the mesh cells and the direction of the
velocity. On meshes approximating the sphere (and other manifolds),
these formulae require modification when $n^+\neq -n^-$ on an edge.
The modification rotates $\tilde{u}$ into the tangent plane of the
cell $e$, as described in \citet{bernard2009high}. The projection of
the advection equation into the tangent to the sphere is naturally
dealt with in \eqref{eq:H div advection}, because $w\in
\mathbb{V}_h^2$ is always tangential to the surface mesh.

The vorticity form is an alternative form of the vector advection
equation, given by
\begin{equation}
  \pp{v}{t} + (\nabla\times v)\times u + \frac{1}{2}\nabla(u\cdot v)
  + \frac{1}{2}\left(
  (\nabla v)^Tu - (\nabla u)^Tv
  \right) = 0,
\end{equation}
in three dimensions, where $(\nabla v)^T_{ij} = \partial u_j/\partial
x_i$. In two dimensions this is written
\begin{equation}
  \label{eq:vector eqn}
  \pp{v}{t} + (\nabla^\perp \cdot v)\cdot u^\perp  + \frac{1}{2}\nabla(u\cdot v)
  + \frac{1}{2}\left(
  (\nabla v)^Tu - (\nabla u)^Tv
  \right) = 0,
\end{equation}
where $w^\perp=(-w_2, w_1)$, $\omega=\nabla^\perp\cdot
w:=-\pp{w_2}{x_1} + \pp{w_1}{x_2}$ for a vector field $w$ in planary
geometry. On the sphere, with outward pointing normal $\hat{k}=x/|x|$
(and $x$ is the three dimensional coordinate with origin at the centre
of the sphere), we have $w^\perp=\hat{k}\times w$ and
$\nabla^\perp\cdot w=\hat{k}\cdot \nabla \times w$, where $\nabla$ is
now the projection of the gradient into the tangent to the mesh
surface (see \citet{rognes2013automating} for implementation details).

When $v=u$, we have the ``vector invariant form'',
\begin{equation}
  \pp{u}{t} + \underbrace{(\nabla\times u)\times u}_{\mbox{or }
    (\nabla^\perp\cdot u)u^\perp}+ \frac{1}{2}\nabla |u|^2=0,
\end{equation}
which is particularly useful on the surface of the sphere because it
avoids the need to rotate upwinded vectors. To see this, we multiply
by a test function $w\in V$ and integrate over one cell, integrating
by parts to get
\begin{equation}
  \int_e w\cdot \pp{u}{t} - \nabla^\perp(w\cdot u^{\perp})\cdot u
  - \nabla\cdot w \frac{1}{2}|u|^2 \diff x
  + \int_{\partial e} w\cdot u^{\perp}n^\perp\cdot \tilde{u} \diff S
  = 0, \quad \forall w \in V.
\end{equation}
Here, no rotation is required because the tangent to the edge between
cells agrees on both sides, it is just the facet normal (the normal
to the cell edge that is tangential to the cell surface)
that can change on
manifold meshes. Summing over all of the cells in the mesh gives
\begin{equation}
  \label{eq:Hdiv vi}
  \int_\Omega w\cdot \pp{u}{t} - \nabla^\perp(w\cdot u^{\perp})\cdot u
  - \nabla\cdot w \frac{1}{2}|u|^2 \diff x
  - \int_{\Gamma} \jump{w\cdot u^{\perp}}\cdot \tilde{u}^\perp \diff S
  = 0,
\end{equation}
where for scalars $\phi$, $\jump{\phi}=\phi^+n^+ + \phi^-n^-$.
This upwinded vector invariant form for $H(\ddiv)$
spaces first appeared in \citet{natale2018variational} for the
incompressible Euler equations and was used for the rotating shallow
water equations on the sphere in \citet{gibson2019compatible}.

The polynomial spaces for lowest order RT elements do not span all
linear vector fields. This means that if we use the above scheme then
it will only be first order accurate. \citet{bendall2023improving}
looked at using an auxiliary field $q\in \mathbb{V}^0_h$ with
$q=-\delta_0 v$, \emph{i.e.} $q$ approximates $\nabla^\perp\cdot v$ if
the solution domain $\Omega$ has no boundary. Then, we can use
this in an approximation of \eqref{eq:vector eqn},
\begin{equation}
  \label{eq:unstabilised vorticity form}
  \left\langle w, \pp{v}{t} \right\rangle + \left\langle w, qu^\perp \right\rangle
  - \frac{1}{2}\left\langle \nabla\cdot w, u\cdot v \right\rangle
  + G'(v; w), \quad \forall w \in \mathbb{V}_h^1.
\end{equation}
This was inspired by the energy enstrophy conserving schemes that we
discuss in Section \ref{sec:poisson}. Here, $G'$ represents the
discretisation of the last two terms on the left hand side of
\eqref{eq:vector eqn}, which we do not go into here (standard upwind
discontinuus Galerkin approaches were used, similar to those above).
To obtain the dynamics for $q$, we can select $w=-\nabla^\perp \gamma$
for $\gamma\in \mathbb{V}_h^0$, and substitute into
\eqref{eq:unstabilised vorticity form} to obtain
\begin{equation}
  \left\langle \gamma , \pp{q}{t} \right\rangle - \left\langle \nabla \gamma, qu \right\rangle
  - G'(v; \nabla^\perp\gamma), \quad \forall \gamma \in \mathbb{V}_h^0.
\end{equation}
We recognise the first two terms as the standard continuous finite
element approximation of $\partial q/\partial t + u\cdot \nabla q$.
As we discussed above, some method of stabilisation is usually needed
to suppress oscillations in this approximation. Modifying the test
function according to the SUPG approach is ungainly here, because
of the surface terms in $G'$. Instead, \citet{bendall2023improving}
used a residual based approach, writing
\begin{equation}
  \left\langle \gamma , \pp{q}{t} \right\rangle - \left\langle \nabla \gamma, q^*u \right\rangle
  - G'(v; \nabla^\perp\gamma), \quad \forall \gamma \in \mathbb{V}_h^0,
\end{equation}
where
\begin{equation}
  q^* = q - \tau\left( \pp{q}{t} + \nabla\cdot(uq) +
  \frac{1}{2}\nabla_h^\perp \left( (\nabla_hv)^Tu - (\nabla_h u)^Tv
  \right)\right),
\end{equation}
where $\tau$ is a stabilisation parameter. We note that this change
preserves the consistency of the discretisation since the quantity in
the brackets is just the curl ($\nabla^\perp\cdot$) of
\eqref{eq:vector eqn}. Then, the equation for $v$ becomes
\begin{equation}
  \left\langle w, \pp{v}{t} \right\rangle + \left\langle w, q^*u^\perp \right\rangle
  - \frac{1}{2}\left\langle \nabla\cdot w, u\cdot v \right\rangle
  + G'(v; w), \quad \forall w \in \mathbb{V}_h^1.
\end{equation}
\citet{bendall2023improving} showed in numerical experiments
(including a nonlinear rotating shallow water equation test case on
the sphere) that this discretisation produces second order accurate
solutions using lowest order RT quadrilateral elements on a cubed
sphere grid, whilst \eqref{eq:Hdiv vi} only produces first order
accurate solutions. A predecessor of this scheme was considered
in \citet{kent2022mixed}, but without the consistent definition
of $q$ and $u$.

Transport schemes for $\mathbb{W}_h^1$ are discussed in
\citet{wimmer2022structure} in the context of magnetohydrodynamics, using
similar ideas to those discussed in this section for $\mathbb{W}_h^2$.

\subsection{Temperature space transport schemes}
The temperature space $\mathbb{W}^\theta$ produces similar challenges,
since it is continuous in the vertical and discontinuous in the
horizontal. \citet{yamazaki2017vertical} combined an upwind
discontinuous Galerkin discretisation in the horizontal with an SUPG
discretisation in the vertical, with the modification of test
functions $\gamma\mapsto \gamma + \tau \hat{k}\cdot u\pp{\tau}{z}$,
producing a second order scheme when the $RT_1$ discrete de Rham complex
is used, so that temperature is continuous quadratic in the vertical and
discontinuous linear in the horizontal.

In staggered grid weather models, it is standard practice to collocate
thermodynamic tracers such as moisture, \emph{etc.}, with temperature.
This makes it easier to localise thermodynamic processes that alter,
and depend on, the temperature (see \citet[for
  example]{bush2020first}). Hence, we need to use
$\mathbb{W}_h^\theta$ transport schemes for those tracers as well. For
many of these tracers (moisture, chemical species, \emph{etc.}) it is
important to avoid numerical over- and undershoots leading to negative
humidity, for example. Hence, it is important to be able to
incorporate limiters into $\mathbb{W}_h^\theta$ transport schemes.
\citet{cotter2016embedded}  proposed such a scheme for the vertically
quadratic, horizontally linear $\mathbb{W}_h^\theta$ also considered
by \citet{yamazaki2017vertical}. In that scheme, at the beginning of
the timestep, the vertical continuity conditions are relaxed, and an
upwind discontinuous Galerkin transport scheme is applied over one
timestep in $\hat{\mathbb{W}}_h^\theta$, the corresponding
discontinuous space. A slope limiter, such as the one in
\citet{kuzmin2013slope}, can then be used to avoid under- and
overshoots in this step. Then, an element based flux corrected
remapping is used to transform $\theta$ back to the vertically
continuous space $\mathbb{W}_h^\theta$.  The flux correction switches
between a high order and low order mapping into order to maximise the
use of the high order solution unless under- or overshoots would
otherwise occur.

\subsection{Recovered space schemes} The lowest order RT de Rham complex
on hexahedra is attractive because it allows storage of field values
at the same grid locations as for the C grid finite difference
approximation. This is why this de Rham complex is being used for the
Met Office ``Gung Ho'' dynamical core \citep{melvin2019mixed}.
However, as we have already mentioned above, standard upwind finite
element schemes on these spaces are only first order accurate, because
only $\mathbb{W}_h^0$ has element shape functions that span the
complete linear polynomial space; $\mathbb{W}_h^i$, $i=1,2,3$, and
$\mathbb{W}_h^\theta$ do not. \citet{bendall2019recovered} chose to
address this by using recovery operators to construct higher order
approximations of the solution based on averaging cell values around
vertices. If the original low order solution is obtained by
interpolating a smooth function, this recovery step produces a higher
order continuous finite element solution
\citep{georgoulis2018recovered}. Following \citet{cotter2016embedded},
they then relaxed the continuity of the higher order recovered
solution and applied a discontinuous Galerkin transport scheme step
before remapping back to the original low order finite element
spaces. This was demonstrated in numerical experiments to produce
second order convergence of solutions. \citet{bendall2023improving}
introduced modifications to extend this approach to the sphere.

This recovery process also allows the introduction of limiters to
prevent over- and undershoots.  This was done in
\citet{bendall2020compatible}, applied to compressible Euler solutions
with moisture, where limiters are critical for stability. This
produced the first atmospheric simulations using compatible finite
elements with moist physics. \citet{bendall2022solution} then showed
how to achieve this framework whilst conserving mass and total
moisture.

\section{Example discretisations and iterative solution strategies}
\label{sec:examples}
In this section, we survey some compatible finite element
discretisations of geophysical fluid dynamics models, concentrating on
approaches that can be considered as evolutions of existing approaches
using more ``traditional'' discretisations. More advanced structure
preserving discretisations are discussed in Sections
\ref{sec:variational}-\ref{sec:nonaffine}.

\localtableofcontents

\subsection{Rotating shallow water equations on the sphere}
We start with the rotating shallow water equations on the sphere, written
\begin{align}
  \pp{u}{t} + \mathbb{P}_S\left((u\cdot \nabla)u\right) +
  fu^\perp + g\nabla (D + b)
  & = 0, \\
  \pp{D}{t} + \nabla\cdot (Du) & = 0,
\end{align}
where $u$ is the horizontal velocity tangential to the sphere, $D$ is the
depth of the fluid layer, and $b$ is the height of the bottom surface.

\citet{gibson2019compatible} introduced a spatial discretisation built
around the vector invariant formulation \eqref{eq:Hdiv vi} for
velocity advection $u$, and the discontinuous Galerkin formulation
\eqref{eq:DG} for depth $D$. This leads a spatial discretisation seeking
$(u,D)\in \mathbb{V}_h^1\times \mathbb{V}_h^2$ such that
\begin{align} \nonumber
  \left\langle w, \pp{u}{t} \right\rangle
  - \left\langle \nabla^\perp_h(w\cdot u^{\perp}), u \right\rangle
  + \left\llangle \jump{w\cdot u^{\perp}}, \tilde{u}^\perp \right\rrangle_{\Gamma} & \\
   \qquad - \nabla\cdot w \left(\frac{1}{2}|u|^2 + g(D+b)\right) \diff x
   &= 0, \quad \forall w \in \mathbb{V}_h^1, \\
    \left\langle \phi, D_t \right\rangle - \left\langle \nabla_h \phi, uD \right\rangle
  + \left\llangle \jump{\phi u}, \tilde{D} \right\rrangle_\Gamma
  &= 0, \quad \forall \phi \in \mathbb{V}_h^2.
\end{align}
They used a semi implicit timestepping scheme, which is best described as
some form of iteration towards the fully implicit timestepping scheme
\begin{align} \nonumber
  \left\langle w, u^{n+1} - u^n \right\rangle - \Delta t\left\langle \nabla^\perp_h(w\cdot
  \bar{u}^{\perp}), (u^{n+1/2})^{\perp} \right\rangle &
  \\ \nonumber
  \qquad + \Delta t\left\llangle \jump{w\cdot
    \bar{u}^{\perp}}, (\tilde{u}^{n+1/2})^\perp \right\rrangle_{\Gamma} &
  \\ \qquad - \Delta t\nabla\cdot w \left(\frac{1}{2}|\bar{u}|^2
  + g(\bar{D}+b)\right)
  \diff x &= 0, \quad \forall w \in \mathbb{V}_h^1,
  \label{eq:swe midpoint u}
  \\ \left\langle \phi,
  D_t \right\rangle - \Delta t\left\langle \nabla_h \phi, \bar{u}D \right\rangle +
  \Delta t\left\llangle
  \jump{\phi \bar{u}}, \tilde{D} \right\rrangle_\Gamma &= 0, \quad \forall \phi
  \in \mathbb{V}_h^2,
  \label{eq:swe midpoint D}
\end{align}
where $u^{n+1/2}=\bar{u}=(u^{n+1}+u^n)/2$. This mixture of two symbols
for the same thing is introduced to describe an iteration based around
this implicit discretisation. We write $v^0=u^n,v^1,v^2,\ldots$ and
$D^0=D^n,D^1,D^2,\ldots$ as a sequence of iterative approximations
to $u^{n+1}$ and $u^n$ respectively. For each iteration $k$, we write
$\bar{u}=(u^n+v^k)/2$, $\bar{D}=(D^n+D^k)/2$, solving
(\ref{eq:swe midpoint u}-\ref{eq:swe midpoint D})
for $u^{n+1}$ and $D^{n+1}$ using those values of $\bar{u}$ and $\bar{D}$.
Then, we use the linearisation about the state of rest to compute iterative
corrections $(\Delta u, \Delta D)$, according to
\begin{align}
  \nonumber
  \left\langle w, \Delta u \right\rangle + \frac{\Delta t}{2}\left\langle w, f\Delta u^\perp \right\rangle & \\
  \qquad - \frac{\Delta t}{2}\left\langle \nabla\cdot w, g\Delta D \right\rangle
  & = -R_u[w] := -\left\langle w, u^{n+1} - u^n \right\rangle, \quad \forall w \in \mathbb{V}_h^1,    \label{eq:swe u linear update} \\
  \label{eq:swe D linear update}
  \left\langle \phi, \Delta D + \frac{H\Delta t}{2}\nabla\cdot \Delta u \right\rangle
  & = -R_D[\phi] := -\left\langle \phi, D^{n+1} - D^n\right\rangle, \quad \forall \phi\in \mathbb{V}_h^2.
\end{align}
We will discuss the solution of this linear system later. The time
integration scheme applies a fixed number of iterations of this type
(typically $2\leq k_{\max} \leq 4$). We then take
$(u^{n+1},D^{n+1}) =(v^{k_{\max}}, D^{k_{\max}})$
before moving to the next timestep.  It is not intended to converge to
the solution of the implicit midpoint rule but just to produce a
second order semi implicit scheme that is stable conditional on the
advective Courant number $|u|\Delta t/\Delta x$, but unconditionally
in the wave Courant number $\sqrt{gH}\Delta t/\Delta x$.  A probably
more stable approach is to use the form \eqref{eq:vector eqn} for the
velocity equation, substituting $\bar{u}$ in for $u$, and $u$ in for
$v$. This was done using quadrilateral $RT_0$ elements in
\citet{bendall2023improving}.  Another approach is to replace the
implicit midpoint rule update for $u^{n+1}$ and $D^{n+1}$ given
$\bar{u}$, instead using an explicit transport step (or several
substeps). This can facilitate more sophisticated transport schemes
with limiters that are hard to implement in implicit schemes. This was
also done using quadrilateral $RT_0$ elements in \citet{bendall2023improving}.

\subsection{Rotating incompressible Boussinesq equations}
For ocean models, the most common setting is the incompressible
Boussinesq equations,
\begin{align}
  \pp{u}{t} + (u\cdot \nabla)u +
  2\Omega\times u + \nabla p - b\hat{k} & = 0, \\
  \pp{b}{t} + (u\cdot \nabla)b & = 0, \\
  \nabla \cdot u & = 0,
\end{align}
where $\Omega$ is the Earth's rotation rate, $p$ is the pressure, and
$b$ is the buoyancy. Here we limit discussion to the rigid lid
approximation with boundary conditions $u\cdot n=0$ on all boundaries.
In full ocean models there are generally two thermodynamic tracers,
the potential temperature $\theta$ and the salinity $S$, which are
both transported by advection equations as $b$ is above, and $b$ is
then a specified function of $\theta$ and $S$ \emph{via} an equation
of state. However, here we keep things to the simple formulation
above. Further, many models make the hydrostatic approximation, but we
do not discuss that here. Finally, these equations are extended in
operational models to include mixing parameterisations and
representations of other physical processes; we do not discuss those
either.

\citet{yamazaki2017vertical} introduced a hybrid approach for velocity
advection in the vertical slice setting (the velocity is three
dimensional but all fields are independent of $y$, so the equations
can be solved on a two dimensional mesh in the $x-z$ plane), with
\eqref{eq:Hdiv vi} as the transport scheme for the $x-z$ components of
velocity, and a standard upwind discontinuous Galerkin scheme for the
$y$ component. For bouyancy $b$ (represented in $\mathbb{W}^\theta$)
the hybrid scheme with upwind discontinuous Galerkin in the horizontal
and SUPG in the vertical was used. This spatial discretisation was
combined with a similar timestepping scheme to the one above,
resulting in a linear system of the form,
\begin{align}
  \nonumber
  \left\langle w, \Delta u \right\rangle + \frac{\Delta t}{2} \left\langle w, 2\Omega \times \Delta u \right\rangle & \\
\qquad  - \frac{\Delta t}{2}
  \left\langle \nabla\cdot w, \Delta p \right\rangle - \left\langle w, b\hat{k}\right\rangle
  \label{eq:bouss linear u}  & = -R_u[w], \quad \forall w\in \mathbb{W}_h^2, \\
  \label{eq:bouss linear b}
  \left\langle \gamma, \Delta b\right\rangle
  + \frac{\Delta t}{2}\left\langle \gamma, \Delta u\cdot \hat{k}B_z\right\rangle
  & = - R_b[\gamma], \quad \forall \gamma \in \mathbb{W}_h^\theta, \\
  \label{eq:bouss linear p}
  \left\langle \phi, \nabla\cdot \Delta u \right\rangle & = 0, \quad \forall
  \phi \in \mathbb{W}_h^3,
\end{align}
where $B_z$ is the vertical derivative of a reference buoyancy
profile, to compute the iterative linear corrections to $u^{n+1}$,
$p^{n+1}$ (actually an approximation to pressure at time level
$t^{n+1/2}$) and $b^{n+1}$, analogously to (\ref{eq:swe u linear
  update}- \ref{eq:swe D linear update}).

Equations (\ref{eq:bouss linear u}-\ref{eq:bouss linear p}) were
solved by eliminating $\Delta b$. This is possible without introducing
errors on an extruded mesh with flat layers when $B_z$ is constant,
since $\Delta b \hat{k} \in \mathbb{W}_h^{2,V}$. When that is not the
case, the errors from the approximate elimination can be removed by
using the elimination as a preconditioner for a Krylov method on
(\ref{eq:bouss linear u}-\ref{eq:bouss linear p}).  At the time, the
equations were modified by setting $\Omega$ to zero on the left hand
side (which prevented the iterative solver from being robust to large
values of $\Omega$), and using GMRES applied to the coupled system,
\begin{align}
  \label{eq:u bouss reduced} \nonumber
  \left\langle w, \Delta u + \frac{\Delta t^2}{4}\hat{k}\hat{k}\Delta uB_z
  \right\rangle
  + \frac{\Delta t}{2} \left\langle w, 2\Omega \times \Delta u \right\rangle & \\
  \qquad
  - \frac{\Delta t}{2}
  \left\langle \nabla\cdot w, \Delta p \right\rangle - \left\langle w, b\hat{k}\right\rangle
  & = -\tilde{R}_u[w], \quad \forall w\in \mathbb{W}_h^2, \\
  \nabla\cdot \Delta u & = 0.
  \label{eq:p bouss reduced}
\end{align}
This was preconditioned by an $H(\ddiv)$ block preconditioner, which
we do not describe here. However, it is much better to precondition
(\ref{eq:u bouss reduced}-\ref{eq:p bouss reduced}) using a hybridised
solver, which we describe later in this section.

\citet{yamazaki2017vertical} showed that this suite of discretisation
and solver choices produces a scheme that can resolve fronts in the
Eady vertical slice frontogenesis problem to a similar degree as
C grid finite difference methods used previously.  After this paper was
published, the authors experimented with replacing the hybrid
advection scheme for velocity with a full vector invariant form in all
3 components of velocity. In this case, oscillations in the velocity
field quickly emerge when the front sharpens. These oscillations do
not appear when \eqref{eq:H div advection} is used instead. It is
possible that this is related to the Hollingsworth instability
associated to the vector invariant form used with C grid finite
difference methods \citep{hollingsworth1983internal}, but this
requires further investigation.

\subsection{Rotating compressible Euler equations}
For global atmosphere models, a standard approach is to solve
the rotating compressible Euler equations, given by
\begin{align}
  u_t + (u\cdot \nabla)u + 2\Omega \times u
  + c_p\theta \nabla \Pi + g\hat{k} & = 0, \\
  \label{eq:theta transport}
  \theta_t + u\cdot \nabla \theta & = 0, \\
  D_t + \nabla\cdot(uD) & = 0, \\
  \Pi &= E(\theta, D),
\end{align}
where $\theta$ is the potential temperature, $\Pi$ is the Exner
pressure, $D$ is now the density for 3D models, $c_p$ is the specific heat at
constant pressure (a constant parameter in the ideal gas law), and $E$ is a
prescribed function describing the thermal equation of state relating
$\Pi$, $D$ and $\theta$. This form of the equations is known as the
``theta-Pi'' formulation. Other thermodynamic formulations make use of
pressure and temperature directly, or other combinations of variables,
but we do not discuss them here. For simplicity we
consider boundary conditions $u\cdot n=0$ on the bottom and top
boundaries of the domain, although representations of the top of the
atmosphere can be rather more complicated (since the real atmosphere
has a density that decreases with height until it is so low that
assumptions underlying the fluid dynamics model do not hold).

A compatible finite element formulation of the compressible Euler
equations uses $u\in \mathbb{W}_h^2$, $D\in \mathbb{W}_h^3$,
and $\theta\in\mathbb{W}_h^\theta$. Either $\Pi$ is solved as an
independent variable in $\mathbb{W}_h^3$, with the equation
of state being projected into $\mathbb{W}_h^3$, or $\Pi$ can
just be replaced by $E(\theta, D)$ in the velocity equation,
leading to a system for $u$, $\theta$ and $D$. 

\citet{natale2016compatible} proposed a formulation using
\eqref{eq:Hdiv vi} as the transport scheme for velocity, upwind
discontinuous Galerkin for density transport and the hybrid scheme for
potential temperature used in \citet{yamazaki2017vertical}. The main
additional challenge is the discretisation of the pressure gradient
term $-\theta \nabla \Pi$ (which is scaled by $c_p$).  In
\citet{yamazaki2017vertical}, the pressure gradient $\nabla p$ was
integrated by parts, but this is more complicated in the compressible
Euler case because of the presence of $\theta\in \mathbb{W}_h^\theta$,
which can have discontinuities in the horizontal direction across
vertical facets. We need to integrate by parts because $\Pi$ is
discontinuous, whether it is an independent variable in
$\mathbb{W}_h^3$ or the evaluation of $E(D, \theta)$ (since
$D\in \mathbb{W}_h^3$ is discontinuous).
\citet{natale2016compatible} proposed to apply integration by parts
separately in each cell $e$ for this term, obtaining
\begin{equation}
  \int_e \nabla\cdot (w\theta)\Pi \diff x
  - \int_{\partial e}\theta w\cdot n \{\Pi\} \diff S, \quad \forall w\in \mathbb{W}_h^2(e),
\end{equation}
where $\{\Pi\}$ is the average value of $\Pi$
between the inside and the outside of $e$ (since $\Pi$ takes two
values on $\partial e$). Summing this over all cells gives
\begin{equation}
  \left\langle \nabla_h\cdot (w\theta),\Pi\right\rangle
  - \left\llangle \jump{w\theta}, \{\Pi\} \right\rrangle_{\Gamma},
  \quad \forall w\in \mathbb{W}_h^2(\Omega),
\end{equation}
where $\{\Pi\}$ is now defined as $(\Pi^++\Pi^-)/2$. Again, analogously to (\ref{eq:swe u
  linear update}- \ref{eq:swe D linear update}), a semi implicit
timestepping scheme can be used, this time built around a
linearisation about a state of rest with $u=0$, $\theta=\bar{\theta}$
and $\Pi=\bar{\Pi}$. These reference profiles vary in the vertical
only (or following the ENDGame approach \citep{wood2014inherently},
the values of $\bar{\theta}$ and $\bar{\Pi}$ can be used from $\theta$
and $\Pi$ at the previous timestep, but whilst neglecting their
horizontal derivatives in the linearisation to facilitate efficient
solution). Finally, we neglect horizontal derivatives of $\Delta \theta$
appearing in the $\Delta \theta\nabla \bar{\Pi}$ and term,
for the same reason.
This results in the following linear iteration, presented
in \citet{gibson2019} (and used for the numerical results in
\citet{natale2016compatible}),
\begin{align}
  \nonumber
  \left\langle w, \Delta u \right\rangle + \frac{\Delta t}{2}\left\langle
  w, 2\Omega \times \Delta u\right\rangle & \\
  - \frac{c_p\Delta t}{2} \left\langle \nabla_h\cdot (\bar{\theta}w),
  \Delta\Pi \right\rangle 
  \nonumber
  + \frac{c_p\Delta t}{2}\left\llangle \jump{\bar{\theta}w}, \{\Delta \Pi\}
  \right\rrangle_\Gamma & \\
\qquad   -  \frac{c_p\Delta t}{2}\left\langle \nabla\cdot (\hat{k}\Delta \theta\theta w\cdot\hat{k}), \bar{\Pi} \right\rangle
& = -R_u[w], \quad \forall w\in\mathbb{W}_h^2, \\
\left\langle \gamma, \Delta \theta \right\rangle
+ \frac{\Delta t}{2} \left\langle \gamma,\pp{ \bar{\theta}}{z}\Delta u\cdot \hat{k}\right\rangle & = - R_\theta[\gamma], \quad \forall \gamma \in \mathbb{W}_h^\theta,
\\
\left\langle \phi, \Delta D\right\rangle - \frac{\Delta t}{2}
\left\langle \nabla_h\phi, \bar{D}\Delta u\right\rangle
+ \frac{\Delta t}{2} \left\llangle \jump{\phi \Delta u},\{\bar{D}\}\right\rrangle_\Gamma
& = - R_D[\phi], \quad \forall \phi \in \mathbb{W}_h^3, \\
\Delta \Pi & = \pp{E}{\theta}\Delta\theta+\pp{E}{D} \Delta D,
\label{eq:Pi pointwise}
\end{align}
for the iterative updates $(\Delta u, \Delta \theta, \Delta D) \in
\mathbb{W}_h^2\times \mathbb{W}_h^\theta \times \mathbb{W}_h^3$ to
$(u^{n+1}, \theta^{n+1}, D^{n+1})$. Alternatively, as was done
in \citet{melvin2019mixed}, we can add $\Delta \Pi\in \mathbb{W}^3_h$
to the list of independent variables and replace \eqref{eq:Pi pointwise}
with
\begin{equation}
  \left\langle \alpha, \Delta \Pi \right\rangle  = \left\langle \alpha, \pp{E}{\theta}\Delta\theta+\pp{E}{D} \Delta D \right\rangle, \quad \forall \alpha \in \mathbb{W}_h^3.
\end{equation}
\citet{melvin2019mixed} adopted a hybrid approach, using finite volume
methods to approximate the transport terms in a discretisation otherwise
built using compatible finite element methods.

\subsection{Iterative solver strategies}
Now we focus on the iterative solver strategies for these linear
implicit systems. In all of the strategies we discuss here,
the temperature $\Delta \theta$ is first (approximately) eliminated
following our description of the approach to incompressible Boussinesq
equations discussed above. This leads to the following system
\begin{align}
  \nonumber
  \left\langle w, \Delta u \right\rangle + \frac{\Delta t}{2}\left\langle
  w, 2\Omega \times \Delta u\right\rangle & \\
  - \frac{c_p\Delta t}{2} \left\langle \nabla\cdot (\bar{\theta}w),
  \Delta\Pi \right\rangle
  \nonumber
   + \frac{c_p\Delta t}{2}\left\llangle \jump{\bar{\theta}w}, \{\Delta \Pi\}
  \right\rrangle_\Gamma & \\
\qquad   -  \frac{c_p\Delta t}{2}\left\langle \nabla\cdot (\hat{k}\Delta \theta\theta w\cdot\hat{k}), \bar{\Pi} \right\rangle
& = -R_u[w], \quad \forall w\in\mathbb{W}_h^2,
\label{eq:euler u reduced}
\\
\left\langle \phi, \Delta D\right\rangle - \frac{\Delta t}{2}
\left\langle \nabla_h\phi, \bar{D}\Delta u\right\rangle
+ \frac{\Delta t}{2} \left\llangle \jump{\phi \Delta u},\{\bar{D}\}\right\rrangle_\Gamma
& = - R_D[\phi], \quad \forall \phi \in \mathbb{W}_h^3, \\
\Delta \theta &= -\frac{\Delta t}{2} \bar{\theta}\Delta u\cdot \hat{k} +
r_\theta, \\
\Delta \Pi & = \pp{E}{\theta}\Delta\theta+\pp{E}{D} \Delta D,
\label{eq:euler Pi reduced}
\end{align}
for $(\Delta u,\Delta D) \in \mathbb{W}_h^2\times \mathbb{W}_h^3$,
where $r_\theta\in \mathbb{W}_h^\theta$ such that
\begin{equation}
  \left\langle r_\theta, \gamma \right\rangle = R_\theta[\gamma], \quad
  \forall \gamma \in \mathbb{W}_h^\theta.
\end{equation}
In other words, $r_\theta$ is the $L^2$ Riesz representer of $R_\theta$.
\citet{mitchell2016high} proposed to solve this reduced system using
GMRES with a Schur complement preconditioner, using an approximate
Schur complement formed from the lumped velocity mass matrix (and
setting $\Omega=0$). This was incorporated into a horizontal multigrid
scheme (coarsening the mesh in the horizontal but not the vertical)
using line smoothers (direct solves neglecting horizontal coupling
between columns) for the approximate Schur complement on the levels.
This combination of horizontal multigrid and vertical line smoothers
is necessary because of the small aspect ratio of the atmosphere, and
scalable parallel performance was observed over a large range of
resolutions. This solve approach was successfully implemented in the
Met Office system using the discretisation approach of
\citet{melvin2019mixed} by \citet{maynard2020multigrid}.

There is an alternative solution approach that has been applied to
compatible discretisations of elliptic problems since the mid 20th
Century: hybridisation. In hybridisation, the continuity constraints
of the $H(\ddiv)$ space are relaxed, and are enforced through Lagrange
multipliers as part of the solution formulation. The Lagrange
multiplier space $\Tr$, known as the trace space since it is supported
only on cell facets, is chosen to match the $H(\ddiv)$ space when
restricted to a facet and dotted with the normal component. This means
that it is discontinuous between facets that meet at a vertex in two
dimensions, or at an edge or vertex in three dimensions. For
example, the hybridisable formulation of (\ref{eq:swe u linear
  update}-\ref{eq:swe D linear update}) seeks
$(\Delta u, \Delta D, \lambda)\in \hat{\mathbb{V}}_h^1\times \mathbb{V}_h^2
\times \Tr(\mathbb{V}_h^1)$ such that
\begin{align}\nonumber
  \left\langle w, \Delta u \right\rangle + \frac{\Delta t}{2}\left\langle w, f\Delta u^\perp \right\rangle & \\
  \quad - \frac{\Delta t}{2}\left\langle \nabla\cdot w, g\Delta u \right\rangle
  + \left\llangle \jump{w}, \lambda \right\rrangle_\Gamma
  \label{eq:hyb swe u}  & = -\tilde{R}_u[w], \quad \forall w \in \hat{\mathbb{V}}_h^1, \\
  \left\langle \phi, \Delta D + \frac{H\Delta t}{2}\nabla\cdot \Delta u \right\rangle
  & = -R_D[\phi], \quad \forall \phi\in \mathbb{V}_h^2, \\
  \left\llangle \gamma, \jump{\Delta u}\right\rrangle_\Gamma & = 0, \quad
  \forall \gamma \in \Tr(\mathbb{V}_h^1),
  \label{eq:hyb swe lambda}
\end{align}
where
\begin{equation}
  \tilde{R}_u[w] = R_u[w], \quad \forall w \in \mathbb{V}_h^1.
\end{equation}
To see that this is an equivalent formulation to (\ref{eq:swe u linear
  update}-\ref{eq:swe D linear update}), note that
$\hat{\mathbb{V}}_h^1 \subset \mathbb{V}_h^1$, so we may choose $w\in
\mathbb{V}_h^1$ in \eqref{eq:hyb swe u}. In that case, the $\lambda$
term vanishes because $\jump{w}=0$, and we recover \ref{eq:swe u
  linear update}. Further, \eqref{eq:hyb swe lambda} ensures that
$u\in \mathbb{V}_h^1$ at the solution, since taking
$\gamma=\jump{\Delta u}$ implies that $\jump{\Delta u}=0$ in
$L^2(\Gamma)$; we note that
\begin{equation}
  \mathbb{V}_h^1(\Omega) = \{u \in \hat{\mathbb{V}}_h^1(\Omega):
  \|\jump{u}\|_{L^2(\Gamma)}=0\}.
\end{equation}
The advantage of this formulation is that $\Delta u$ and $\Delta D$
can now both be eliminated elementwise, leading to a sparse system for
$\lambda$; the reduced system is referred to as the hybridised system.
This is possible because we can take $w$ and $\phi$ supported in only
one cell, and then given $\lambda$, we can solve for $\Delta u$ and
$\Delta D$ independently in each cell (this is referred to as the
local solver). Material discussing the well posedness of the
hybridised system is surveyed in \citet{boffi2013mixed}, along with
postprocessing techniques for obtaining improved approximations using
$\lambda$. A nonrigorous intuitive explanation for this is that
$\lambda$ gives an approximation of $\Delta D$ evaluated on facets,
and hence the hybridised equation has properties of an approximation
to the Helmholtz equation satisfied by $\Delta D$ after eliminating
$\Delta u$ from the linear PDE. This idea is built upon in
\citet{cockburn2004characterization}, which provides an explicit weak
form characterisation of the hybridised method (including the
nonsymmetric term here containing $\Omega\times u$ is a
straightforward extension of that work). This idea was used to
demonstrate smoothing properties of standard iterative methods in
\citet{gopalakrishnan2003schwarz} and a convergent multigrid scheme in
\citet{gopalakrishnan2009convergent}. \citet{gibson2020slate} applied
the hybridisation technique to the linear compressible Boussinesq
equations, which are a minor modification of the linearisation of the
incompressible Boussinesq equations described above, incorporating
linear acoustic waves.

This strategy cannot be applied to (\ref{eq:euler u reduced}
-\ref{eq:euler Pi reduced}), because the averaging $\{\Delta \Pi\}$ of
$\Pi$ on facets couples the values of $\Delta \Pi$ between cells (and
so there is no local solver). \citet{gibson2019} proposed a modification
to address this, in which we seek
$(\Delta u, \Delta D, \lambda)\in \hat{\mathbb{W}}_h^2\times\mathbb{W}_h^3
\times \Tr(\mathbb{W}_h^2)$ such that
\begin{align}
  \nonumber
  \left\langle w, \Delta u \right\rangle + \frac{\Delta t}{2}\left\langle
  w, 2\Omega \times \Delta u\right\rangle & \\
  - \frac{c_p\Delta t}{2} \left\langle \nabla_h\cdot (\bar{\theta}w),
  \Delta\Pi \right\rangle
  \nonumber
  \qquad + \frac{c_p\Delta t}{2}\left\llangle \jump{\bar{\theta}w}, \lambda
  \right\rrangle_\Gamma & \\ \nonumber
  \qquad   -  \frac{c_p\Delta t}{2}\left\langle \nabla\cdot (\hat{k}\Delta \theta\theta w\cdot\hat{k}), \bar{\Pi} \right\rangle & \\
  \qquad +\underbrace{\frac{c_p\Delta t}{2}\left\llangle \jump{\hat{k}\Delta \theta\theta w\cdot\hat{k}}, \{\Pi\} \right\rrangle_{\Gamma_H}}_{\star}
& = -R_u[w], \quad \forall w\in\hat{\mathbb{W}}_h^2,
\\
\left\langle \phi, \Delta D\right\rangle - \frac{\Delta t}{2}
\left\langle \nabla\phi, \bar{D}\Delta u\right\rangle
+ \frac{\Delta t}{2} \left\llangle \jump{\phi \Delta u},\{\bar{D}\}\right\rrangle_\Gamma
& = - R_D[\phi], \quad \forall \phi \in \mathbb{W}_h^3, \\
\left\llangle \gamma, \jump{u}\right\rrangle_\Gamma & = 0, \quad
\forall \gamma \in \Tr(\mathbb{W}_h^2), \\
\Delta \theta &= -\frac{\Delta t}{2} \bar{\theta}\Delta u\cdot \hat{k}, \\
\Delta \Pi & = \pp{E}{\theta}\Delta\theta+\pp{E}{D} \Delta D,
\end{align}
where $\Gamma_H$ is the set of horizontal faces between cells in the
same vertical column.  Here, the idea is that $\lambda$ is an
approximation of $c_p\Delta t\Delta\Pi/2$ on mesh facets. This system
is not equivalent to (\ref{eq:euler u reduced} -\ref{eq:euler Pi
  reduced}), though. Note the addition of the term indicated with
$\star$. This term vanishes when $w\in\mathbb{W}_h^2$, so it does not
change the solution, but it was found that without it, iterative
solvers do not perform well; it appears to be required for the
coercivity of the solution. This system has not been analysed yet, but
it was demonstrated to produce comparible results to standard test
cases when applied to the full nonlinear compressible Euler equations.
Scalable multigrid behaviour for the hybridised system was also
demonstrated. \citet{bendall2020compatible} extended this solver
approach to the compressible Euler equations with moisture, where the
mass lumping approach was found not to work well when $RT_1$ spaces
were used. \citet{betteridge2022hybridised} showed that this
hybridisation approach produces scalable results for the Met Office
formulation of \citet{melvin2019mixed}.

\subsection{Computing hydrostatic balanced states}
The hybridisation approach also provides a useful way to solve
for hydrostatic balanced states, satisfying
\begin{equation}
  -\left\langle \nabla\cdot (\theta w), \Pi \right\rangle
  = \left\langle w\cdot \hat{k}, g\right\rangle
  - \left\llangle w\cdot n, \Pi_0 \right\rrangle
  \quad \forall w \in \mathring{\mathbb{W}}_h^{2,V}.
\end{equation}
Following the technique of \citet{natale2016compatible} previously
discussed, for given $\theta\in \mathbb{W}_h^\theta$, we seek
$(v,D)\in\mathbb{W}_h^2\times \mathbb{W}_h^3$ such that
\begin{align}
  \left\langle w, v\right\rangle -\left\langle \nabla\cdot (\theta w), \Pi \right\rangle
  &= \left\langle w\cdot \hat{k}, g\right\rangle
  - \left\llangle w\cdot n, \Pi_0 \right\rrangle
  \quad \forall w \in \mathring{\mathbb{W}}_h^{2,V}, \\
  \left\langle \phi, \nabla\cdot v \right\rangle & = 0, \quad \forall \phi
  \in \mathbb{W}_h^3, \\
  \Pi &= E(\theta, D),
\end{align}
which is independent between columns. An equivalent hybridisable
formulation seeks 
$(v,D)\in\mathbb{W}_h^2\times \mathbb{W}_h^3$ such that
\begin{align}\nonumber
  \left\langle w, v\right\rangle -\left\langle \nabla\cdot (\theta w), \Pi \right\rangle & \\
\qquad  + \left\llangle \jump{w}, \lambda \right\rrangle_\Gamma
  &= \left\langle w\cdot \hat{k}, g\right\rangle
  - \left\llangle w\cdot n, \Pi_0 \right\rrangle
  \quad \forall w \in \hat{\mathbb{W}}_h^{2,V}, \\
  \left\langle \phi, \nabla\cdot v \right\rangle & = 0, \quad \forall \phi
  \in \mathbb{W}_h^3, \\
  \left\llangle \gamma, \jump{v}\right\rrangle_\Gamma  & = 0, \quad
  \forall \gamma \in \Tr(\mathbb{W}_h^{2,V}), \\
  \Pi &= E(\theta, D),
\end{align}
where now $\Gamma$ includes the bottom boundary (but not the top),
and $\Tr(\mathbb{W}_h^{2,V})$ is only supported on horizontal facets
between neighbouring cells in the same column. This can be solved
using Newton's method, with solution of the Jacobian system \emph{via}
the hybridised system for $\delta \lambda$.

\subsection{Monolithic solvers}
One more recent solver approach has been investigated in numerical experiments
in \citet{cotter2022compatible}, who used a similar suite of
discretisations to that of \citet{natale2016compatible}, with the
exception of using edge stabilisation for potential temperature
instead of vertical SUPG.  In this work, the fully nonlinear implicit
midpoint rule is solved using Newton's method, and GMRES is applied to
the coupled system for $\Delta u,\Delta \theta, \Delta D$ without
elimination. The system is preconditioned by an additive Schwarz
method, computing exact solutions of the Jacobian system restricted to
overlapping columnar patches; each patch consists of the cells
surrounding one vertical edge (excluding degrees of freedom attached
to the vertical facets on the side boundaries of the patch).  This
scheme also requires further analysis, but was shown to produce mesh
independent iteration counts in numerical experiments.

\section{Variational discretisations}
\label{sec:variational}

Variational discretisations are discretisations that are derived from
a discrete Hamilton's principle. They were originally introduced in
the setting of ordinary differential equations (ODEs), taking the name
``variational integrators'', surveyed in \citet{marsden2001discrete}.
The idea behind variational discretisations is that rather than
discretising the equations directly, we instead discretise the
action functional from which the equations are derived. In the
case of ODEs, this means replacing the time integral with a discrete
quadrature rule involving the solution at discrete points in time.
The discretisation of the equations is then obtained by finding
stationary points of the discretised Lagrangian.

The advantage of variational discretisations is that if the
discretised action has symmetries, then these symmetries give rise to
corresponding conserved quantities via (the discrete) Noether's
theorem. In the context of variational integrators for mechanical
systems, this yields discrete conservation of momentum and angular
momentum, for example. Further, after Legendre transformation the
discrete timestepping map is symplectic, leading to the conservation
(up to exponentially small terms in $\Delta t$) of a modified
energy/Hamiltonian obtained through backward error analysis
\citep{hairergeometric,hairer2003geometric,leimkuhler2004simulating,sanz1992symplectic}.

In principle, the variational discretisation approach can be extended
directly to partial differential equations by simply discretising the
action functional in space as well as in time (or one may consider
spatial semidiscretisation by discretising in space only, as we shall
mostly do here). For fluid dynamics, the situation is more
challenging, because the underlying variational principle is defined
in terms of the Lagrangian flow map, rather than the Eulerian
quantities. This is discussed in the following subsection.

\localtableofcontents

\subsection{Hamilton's Principle for fluid dynamics: continuous theory}
The Lagrangian flow map is treated formally (we avoid
discussions of smoothness \emph{etc.} here) as a diffeomorphism
$\chi(\cdot, t):\Omega_0\to \Omega$, mapping labels in a configuration
space $\Omega_0$ to fluid particle locations at time $t$ in the
physical domain $\Omega$. We have a time-$t$ family of
maps\footnote{Here we use the notation $\chi(\cdot, t)$ to indicate
the function $x\mapsto\chi(x,t)$ for given $t$.}.  The variational
formulation then follows by writing an action as an integral over
$\Omega_0$, and variations in $\chi$ are considered subject to the
usual endpoint conditions in time plus the requirement that $\chi$ be
a diffeomorphism. Working with these Lagrangian flow maps is difficult
both for theory and numerical computation; this was addressed by
Arnold's geometric formulation in Eulerian variables for the
incompressible Euler equations \citep{arnold1966geometrie}. For each $(x,t)\in \Omega\times
[0,T]$, where $T$ is the time interval over which the equations are
being solved, $u(x,t)$ is a vector tangent to $x$ at $\Omega$\footnote{We have
$u(x,t)\in\mathbb{R}^N$ for the case where $\Omega$ is an $N$
dimensional subset of $\mathbb{R}^N$, but in geophysical fluid
dynamics we are also interested in solving problems where $\Omega$ is
the surface of a sphere.} For each $t$, we say that $u(\cdot, t) \in
\mathfrak{X}(\Omega)$, the space of vector fields on $\Omega$.

The formulation stems from the observation that the Eulerian velocity
$u$ satisfies
\begin{equation}
    \label{eq:phi_t}
  \pp{}{t}\chi = u\circ \chi,
\end{equation}
where for functions of $(x,t)$ we write $(u\circ
\chi)(x,t)=u(\chi(x,t),t)$.  To construct Hamilton's principle, we
need to consider perturbations to $\chi$ that are still
diffeomorphisms for each $t$.  If we consider such a continuous
2-parameter family of perturbed maps $\tilde{\chi}(x,t,\epsilon)$, with
$\tilde{\chi}(x,t,\epsilon=0)=\chi(x,t)$, then there exists a continuous
family of vector fields $\tilde{w}(\cdot, t,\epsilon) \in \mathfrak{X}(\Omega)$ such
that
\begin{equation}
  \pp{}{\epsilon}\tilde{\chi} = \tilde{w}\circ\tilde{\chi},
\end{equation}
where we now extend the $\circ$ notation to the case of two
parameters, $\epsilon$ and $t$, writing $(\tilde{w}\circ
\tilde{\chi})(x,t,\epsilon)=\tilde{w}(\tilde{\chi}(x,t,\epsilon),t,\epsilon)$. Differentiating
with respect to $\epsilon$ and evaluating at $\epsilon=0$, we get the
infinitesimally perturbed $\chi$,
\begin{equation}
  \label{eq:delta phi}
  \delta \chi :=
  \pp{}{\epsilon}|_{\epsilon=0}\tilde{\chi} =
   w \circ \chi,
\end{equation}
where 
$w(x,t):=\tilde{w}(x,t,\epsilon=0)$. Here $\delta \chi$ depends on the
direction of perturbation in the parameter $\epsilon$, so we always
consider $\delta \chi$ being defined with respect to a particular
choice of $w$. By differentiating \eqref{eq:phi_t} with respect to
$\epsilon$, we obtain
\begin{align}
  \delta \pp{}{t}\chi & = (\nabla u)\circ \chi \cdot \delta{\chi} + \delta{u}\circ
  \chi,  \\
  & = (\nabla u)\circ \chi \cdot w \circ \chi + \delta{u}\circ
  \chi,  \label{eq:delta dot phi}
\end{align}
where
\begin{equation}
  \delta u:= \pp{}{\epsilon}|_{\epsilon=0} \pp{}{t}{\chi}\circ \chi^{-1}.
\end{equation}
Similarly, by differentiating \eqref{eq:delta phi} with respect
to $t$, we obtain 
\begin{align}
  \pp{}{t} \delta \chi & = (\nabla w)\circ \chi\cdot \pp{}{t}{\chi} + \pp{}{t}{w}\circ
  \chi, \\
  & = (\nabla w)\circ \chi\cdot u \circ \chi + \pp{}{t}{w}\circ
  \chi. \label{eq:dot delta phi}
\end{align}
By subtracting Equations \eqref{eq:delta dot phi} and \eqref{eq:dot delta phi},
noting the symmetry of second derivatives and composing with $\chi^{-1}$,
we obtain
\begin{equation}
  \delta u = \pp{}{t}{w} + \underbrace{u\cdot \nabla w -
    w \cdot \nabla u}_{=[u,w]}.
\end{equation}
Hence, we have successfully related infinitesimal perturbations
(``variations'') in $u$ to variations in $\chi$ entirely in terms
of the Eulerian quantities $u$ and $w$. Arnold used this
calculation to derive the incompressible Euler equations by seeking
$u \in \mathfrak{X}_{\vol}(\Omega)$ (the subspace of $\mathfrak{X}(\Omega)$
containing divergence-free vector fields), that extremise the
reduced action functional 
\begin{equation}
  S = \int_{0}^T \int_{\Omega} \frac{1}{2}|u|^2 \diff x \diff t,
\end{equation}
subject to endpoint conditions $u|_{t=0}=u_0$, $u|_{t=T}=u_T$ (for
some chosen $u_0$ and $u_T$), and the boundary conditions $u\cdot n=0$
on the boundary $\partial\Omega$ of the domain $\Omega$, where $n$ is
the outward pointing normal to $\partial\Omega$. If we perturb $u$
subject to these conditions, this implies that $w_{t=0}=w_{t=T}=0$.

Taking variations gives
\begin{align}
  0 = \delta S & =  \int_{0}^T \int_{\Omega} \delta u \cdot u \diff x \diff t,\\
  \label{eq:EP}
  = & \int_{0}^T \int_{\Omega} (\pp{}{t}{w} + [u,w]) \cdot u \diff x \diff t,\\
  = &\int_{0}^T \int_{\Omega} w\cdot (-\pp{}{t}{u} - \nabla\cdot(u \otimes u)
  \nonumber
 -
  \underbrace{(\nabla u)^Tu}_{=\frac{1}{2}|u|^2}) \diff x \diff t
  \\
  & \quad  + \int_0^T \int_{\partial \Omega} \underbrace{u\cdot n}_{=0}w\cdot u \diff S\diff t  \nonumber \\
  & \qquad + \left[ \int_\Omega u\cdot \underbrace{w}_{=0} \diff x\right]_{t=0}^{t=T},  \quad \forall w \in \mathfrak{X}_{\vol}(\Omega),
\end{align}
where $(u\otimes u)_{ij}=u_iu_j$, and $(\nabla u)_{ij}=\partial u_i/\partial x_j$.
This formally implies that
\begin{equation}
  \mathbb{P}\left(\pp{}{t}{u} + u\cdot \nabla u \right)=0,
\end{equation}
where $\mathbb{P}$ is the $L^2$ projection onto divergence-free vector
fields. In other words, there exists a pressure $p$ such that
\begin{equation}
  \pp{}{t}{u} + u\cdot \nabla u  + \nabla p=0,
\end{equation}
which is the familiar form of the incompressible Euler equation
(together with the divergence free constraint $\nabla\cdot u=0$).

Besides the usual translational and rotational symmetries (assuming
appropriate boundary conditions), the variational formulation has an
additional particle relabelling symmetry, meaning that we can replace
$\chi$ by $\chi\circ \psi$ for any diffeomorphism
$\psi:\Omega_0\to\Omega_0$ and the action remains invariant (in fact,
the value of $u$ does not change). Physically this corresponds to fact
that the fluid physics is independent of the choice of label for a
fluid particle (the label for a fluid particle at time $t$ and
position $x_0$ being $\chi^{-1}(x)$, \emph{i.e.} $\chi(x_0,t)=x$). As
discussed in many places \citep[for
  example]{morrison1982poisson,salmon1998lectures,
  shepherd1990symmetries}, this symmetry leads to conservation of
circulation
\begin{equation}
\dd{}{t}  \int_{C(t)} u \cdot dx = 0,
\end{equation}
for closed loops $C(t)$ that are being transported by the fluid
velocity $u$.

The reduction from the flow map $\chi$ to the Eulerian velocity by
Arnold was characterised in \citet{holm1998euler} as Euler-Poicar\'e
reduction by symmetry. In that work, this framework was extended to
quantities that are advected by the flow, including densities solving
the continuity equation, tracers solving the advection equation, etc.,
leading to the derivation of the full family of geophysical fluid
dynamics models and beyond: quasigeostrophic approximations \citep{holm1998hamilton},
complex fluids \citep{holm2002euler,gay2009geometric},
vertical slice models \citep{cotter2013variational}, pseudocompressible and anelastic approximations \citep{cotter2014variational},
the Boussinesq-alpha model \citep{badin2018geometric},
models with moisture and irreversible processes \citep{gay2019variational}.
The derivation of conservation laws associated with particle
relabelling symmetries through Noether's theorem applied to this
framework was presented in \citet{cotter2013noether,cotter2019particle}.

\subsection{Koopman representation: continuous theory}
The difficulty with adapting Hamilton's Principle as described above to discretisations is that
there does not exist a finite dimensional subspace of flow maps that
closes appropriately under composition. A solution to this, proposed
by \citep{pavlov2011structure}, is to use the Koopman representation
of flow maps. We shall only briefly describe their approach here,
but will provide more detail about the extension to compatible
finite elements shortly.

In the Koopman framework, flow
maps $\chi$ in the group $\Diff(\Omega)$ of diffeomorphisms on
$\Omega$ are represented by elements of $GL(L^2(\Omega))$, the invertible linear maps from $L^2(\Omega)$ to $L^2(\Omega)$.
In particular, flow maps $\chi$ represented by
linear maps $\rho_\chi \in GL(L^2(\Omega))$, defined by
$\rho_\chi\cdot a \equiv \rho_{\chi}(a) = a\circ \chi^{-1}$. Of course, there are many maps in
$GL(L^2(\Omega))$ that can not be written this way, and so
there is not an isomorphism. In fact, the map defines a subgroup of
$GL(L^2(\Omega))$, which we call $G(L^2(\Omega))$, which is
isomorphic to the group $\Diff(\Omega)$. 

To specialise to incompressible flows, we use $GL_0(L^2(\Omega))$, defined as
\begin{align} \nonumber
  GL_0(L^2(\Omega))=
  \Bigg\{
  \rho \in GL(L^2(\Omega)):
  (\rho\cdot a,\rho \cdot b)_{\Omega} = (a,b)_{\Omega}, \\
  \quad \forall a,b,\in L^2(\Omega)\mbox{ and }\rho\cdot c = c, \forall c \in \mathbb{R}
  \Bigg\}.
\end{align}
The $\rho_\chi$ representation defines an isomorphism of a subgroup of
$GL_0(L^2(\Omega))$, which we call $G_0(L^2(\Omega))$, to the group of
volume preserving diffeomorphisms $\Diff_{\vol}(\Omega)$. The group
$G_0(L^2(\Omega))$ was approximated in numerical discretisations of the
incompressible Euler equations in \citet{pavlov2011structure}; the
extension to the full diffeomorphism group was used in the extension
to compressible models of geophysical fluid dynamics in
\citet{desbrun2014variational,bauer2017variational,brecht2019variational,bauer2019towards}.

The Koopman representation also provides an isomorphism between the
Lie algebra $\mathfrak{X}(\Omega)$, corresponding to vector fields on
$\Omega$, and a subspace of the Lie algebra
$\mathfrak{gl}(L^2(\Omega))$, which we call
$\mathfrak{g}(L^2(\Omega))$. For a one parameter family $\chi_s\in
\Diff(L^2(\Omega))$ of maps with $\chi_0=Id$ and
$\dd{}{s}|_{s=0}\chi_s=v$,
\begin{equation}
  \dd{}{s}|_{s=0}\rho_{\chi_s}a = \dd{}{s}|_{s=0}(a\circ \chi_s^{-1})
  := -L_va = v\cdot \nabla a,
\end{equation}
where $L_v$ is called the \emph{Lie derivative}. Corresponding
definitions follow for the restrictions
$\Diff(\Omega)\to\Diff_{\vol}(\Omega)$ and
$GL(L^2(\Omega))$ to $GL_0(L^2(\Omega))$.

\subsection{Previous work on variational discretisations using
  Koopman representation}
\citet{pavlov2011structure} used the Koopman representation to derive a
discrete variational principle that considered finite subspaces
$GL(V)$, where $V$ is the space of cellwise constant functions
defined on a mesh. The subset $G_h \subset GL(V)$ approximating the
flow maps was identified as being generated by maps that only allow
instantaneous fluxes between neighbouring cells (thus approximating
the diffeomorphism property). These fluxes can be described within the
framework of Discrete Exterior Calculus (DEC)
\citep{hirani2003discrete}. The difficulty with trying to discretise
this structure is that it is not possible to find a subset $G_h$
that can be generated by a finite subspace of $\mathfrak{gl}(V)$.
This was addressed by introducing a \emph{nonholonomic constraint} on
the time derivative of the flow map (and corresponding Koopman
representative of $GL(V)$). Remarkably, Hamilton's principle under
these constraints is still reducible (i.e., the $GL(V)$ elements can
be eliminated in favour of their generating vector fields), leading to
a discretisation that can be solved entirely in terms of Eulerian
velocities, yielding a spatial discretisation that corresponds to a
known marker-and-cell scheme when a structured grid of square cells is
used. This spatial discretisation is then combined with quadrature
approximation in time to produce a fully discrete variational
integrator for fluid dynamics.

The principle goal of the variational discretisation is to obtain
schemes with discrete conservation laws.  Before time discretisation,
the discrete action is invariant under time translations, leading to
conservation of energy. Time discretisation breaks this symmetry, but
there is the potential to apply backward error analysis to obtain
exponentially accurate conservation of a modified energy. Whilst
numerical results from these variational schemes do exhibit long time
approximate energy conservation as would be expected from this
backward error analysis, it has not yet been adapted to the type of
nonholonomic constraints occurring in this framework.  The variational
scheme also contains an echo of the circulation theorem.  Instead of
considering loops, we consider currents. These are objects dual to
velocities, defined by the duality pairing
\begin{equation}
  c[v] = \int_\Gamma v \cdot \diff x, \quad\forall v \in \mathfrak{X}(\Omega),
\end{equation}
for some curve $\Gamma$.
Currents are transported by a flow map $\chi$ via
$c\mapsto ((\nabla \chi)c)\circ \chi^{-1}$, allowing Kelvin's
circulation theorem to be reformulated (in the case
of incompressible Euler equations) as
\begin{equation}
  \label{eq:disc kelv}
  \dd{}{t} \left(((\nabla \chi) c)\circ \chi^{-1}\right)[u] = 0.
\end{equation}
There is a discrete analogue of this formula when currents are
approximated by objects dual to velocity fields on the discrete
grid. These discrete currents can be considered to be loops that have
been smoothed out over a finite area. This can be used to define a
discrete circulation that is conserved along the solution, including
the discrete time solution. However, since the approximation of a
discrete current to a continuous current corresponding to an advected
loop gets noisier as time progresses, it is not clear how or whether
these discrete conservation laws constrain the discrete fluid dynamics
in the same way that circulation conservation does in the continuous
case (which occurs through the link to Casimirs on the Hamiltonian
side). This is another important open question about this framework.

\subsection{Compatible finite element variational discretisation of
  incompressible flow}
Inspired by the observation of Dmitry Pavlov that one could extend
this framework to other discretisation methods by simply selecting a
discrete space for velocity fields and a discrete approximation of
their Lie group action on scalar fields, and noting the links between
DEC and FEEC, \citet{natale2018variational} developed such an
extension to compatible finite element methods for the incompressible
Euler case. In this case, we use the space of velocity fields
$\mathring{W}_h^r$ defined by
\begin{equation}
  \mathring{W}_h^r = \left\{
  u \in W_h^r: \nabla\cdot u = 0, \, u\cdot n|_{\partial\Omega}=0
  \right\},
\end{equation}
where $W_h^r$ is a degree $r$ $BDM$ or $RT$ space (they both have the same
divergence free subspace so the distinction is not important here).
For the discrete Koopman representation of the flows generated by
these velocity fields, we select $V^s_h$, the space of discontinuous
piecewise polynomials of degree $s$. Following the discrete Lie
derivative framework (the Eulerian version in particular) of
\citet{heumann2010eulerian}, given $u\in \mathring{W}_h^r$,
we define the discrete advection operator $X_u\in \mathfrak{gl}(V_h^s)$
defined by
\begin{equation}
  \label{eq:X_u}
  \left\langle X_u a, b \right\rangle = -\left\langle a, \nabla\cdot (ub)\right\rangle
  + \left\llangle \jump{au}, \{b\}\right\rrangle_\Gamma,\quad
  \forall a, b \in V_h^s,
\end{equation}
where $\{b\}=(b^++b^-)/2$. This is a centred approximation of the
advection operator on $V_h^s$, meaning that $X_u$ is linear
in $u$ (an upwinded approximation would break this).

To develop the nonholonomic constraint that enforces a dynamics that
converges to fluid motion by diffeomorphism, we seek a Koopman
representation of a time dependent flow map
$\hat{\chi}(\cdot, t)\in G(V_h^s)$ that transports
advected tracers according to $a=\hat{\chi}\cdot a_0$. If we require
that all such advected tracers satisfy the equation,
\begin{equation}
  \dd{}{t}a + X_{u}a_0 = 0,
\end{equation}
for some time dependent velocity field $u\in \mathring{W}^r_h$, then
we obtain the \emph{nonholonomic constraint} on $\hat{\chi}$ that
\begin{equation}
  \dd{}{t} {\hat{\chi}}_h\cdot a +
  X_{u}\hat{\chi}\cdot a = 0, \quad \forall a \in V_h,
\end{equation}
for some $X_{u}\in \mathfrak{gl}(V_h^s)$, i.e.,
\begin{equation}
    \label{eq:nonholo}
  \dd{}{t}\hat{\chi} + X_u\hat{\chi} = 0.
\end{equation}
This constraint describes the approximation $\hat{\chi} \in G_h$ to the
Koopman representation $\hat{\chi}\in G(L^2(\Omega))$ of the flow map
$\chi$. It is nonholonomic because it constrains the time derivative
of $\hat{\chi}$, not $\hat{\chi}$ itself, and this constraint cannot be
integrated to obtain such a constraint $F(\hat{\chi})=0$. This is because
the subspace $S_h^r(V_h^s)\subset \mathfrak{gl}(V_h^s)$ defined by the
image of the map $u\in \mathring{W}_h^r\mapsto X_u$ is not closed
under Lie brackets. In other words $[X_u,X_v]$ is not guaranteed to be
in $S_h^r(V_h^s)$ for all $X_u,X_v\in S_h^r(V_h^s)$.

\citet{natale2018variational} proved that if $r\geq s$, then the map
between $u\in W^r_h$ and $X_u\in S_h^r(V_h^s)\subset \mathfrak{gl}(V_h^s)$
is an isomorphism.  \citet{gawlik2021variational} took this further,
by considering the extension of $X_u$ to the whole of $H(\ddiv)\cap
L^p(\Omega)^n$ (with some $p>2$; this technicality ensures the
existence of traces on individual facets), using the same formula
\eqref{eq:X_u}. They consider the space $\hat{S}_h(V_h^s) \subset
\mathfrak{gl}(V_h^s)$, defined by
\begin{equation}
\hat{S}_h(V_h^s) = \{X_u:u \in H(\ddiv)\cap L^p(\Omega)^n\},
\end{equation}
and proved that $\hat{S}_h(V_h^s)$ is isomorphic to $RT^{2r}$ (via
the isomorphism $u\mapsto X_u$). This shows that subspaces of $RT^{2r}$
are a necessary choice for $W^r_h$. When discontinuous finite element
spaces are chosen for $V$, the compatible finite element framework is
thus a necessity rather than a choice.

To continue the derivation of the discrete incompressible Euler
equations, we need to form the Lagrangian, which is just
the kinetic energy. Defined in terms of the original flow map,
the undiscretised Lagrangian is
\begin{equation}
  \int_\Omega \frac{1}{2}|u|^2\diff x.
\end{equation}
In the finite element framework, we need to write this Lagrangian as a
functional of $X_u \in {S}_h(V_h^s)$, but we need to be able to
recover $u$ from $X_u$ to substitute into the kinetic energy.
As proposed by
Dmitry Pavlov, \citet{natale2018variational} obtained an approximation
to this Lagrangian by applying $X_u$ to each of the Cartesian
coordinates, i.e.
\begin{equation}
  \label{eq:hatX_u}
  u = \overline{X_u} = \sum_{i=1}^N X_u(x_i)e_i,
\end{equation}
where $e_i$ is the unit vector in the direction of increasing
coordinate $x_i$.\footnote{In fact, this was used to prove the
isomorphism between $X_u$ and $u$.} This presentation assumes a
Cartesian metric, and suitable changes need to be made for solution of
the equations on manifolds such as the surface of a sphere.

To properly define Hamilton's
principle, we need to define the Lagrangian on the whole tangent
bundle $TG_h(V)$ of $G_h(V)$ (consisting of pairs $(\pp{}{t}{\hat{\chi}},\hat{\chi})$
with $\hat{\chi}\in G_h(V)$). Thus we extend the definition of the ``overbar map''
\eqref{eq:hatX_u} to the whole of $\mathfrak{gl}(V_h)$,
\begin{equation}
  \overline{\xi} = \sum_{i=1}^N e_i\xi\cdot x_i, \quad \forall \xi \in
  \mathfrak{gl}(V_h),
\end{equation}
and write the Lagrangian $L_h:TG_h(V)\to \mathbb{R}$,
as
\begin{equation}
  \label{eq:euler lagrangian nonholo}
  L_h = \int_\Omega \frac{1}{2}
  \left\|\overline{\left(\pp{}{t}{\hat{\chi}}\circ \hat{\chi}^{-1}\right)}\right\|^2 \diff x.
\end{equation}
We can write this as a reduced Lagrangian on $\mathfrak{gl}(V_h^s)$,
\begin{equation}
  \label{eq:reduced lagrangian nonholo}
  \ell_h[X] = \int_\Omega \frac{1}{2}\left\|\overline{X}\right\|^2\diff x.
\end{equation}

The Lagrange-D'Alembert principle (Hamilton's principle with
nonholonomic constraints) seeks $\hat{\chi}$ with $\pp{}{t}{\hat{\chi}}$
satisfying \eqref{eq:nonholo} such that
\begin{equation}
  \delta \int_0^T L_h(\pp{}{t}{\hat{\chi}},\hat{\chi})\diff t = 0,
\end{equation}
for all variations $\delta \hat{\chi}$ satisfying
\begin{equation}
\delta \hat{\chi} + X_{w}\hat{\chi} = 0,
\end{equation}
for some time dependent $w\in \mathring{W}_h^r$. If we have a
Lagrangian (such as \eqref{eq:euler lagrangian nonholo}) that is
reducible by right action, i.e.
\begin{equation}
  L_h(\pp{}{t}{\hat{\chi}},\hat{\chi}^{-1}) =
  \ell_h(\pp{}{t}\hat{\chi}\circ \hat{\chi}^{-1}),
\end{equation}
then we can perform Euler-Poincar\`e reduction, taking care with
the nonholonomic constraint. This proceeds much as in the unapproximated
case described above, as follows. If we have $g$ such that
\begin{equation}
  \pp{}{t}\hat{\chi} + X\hat{\chi} = 0,
  \, \delta \hat{\chi} + Y \hat{\chi} = 0,
\end{equation}
for $X,Y\in \mathfrak{gl}(V_h^s)$,
then calculations identical to the above lead to
\begin{equation}
  (\delta X)\hat{\chi} + XY\hat{\chi}
  - \pp{}{t}Y \hat{\chi} - YX\hat{\chi} = 0,
\end{equation}
i.e.
\begin{equation}
  \label{eq:variations gl}
  \delta X = \pp{}{t}Y + [X,Y],
\end{equation}
where $[X,Y]=XY-YX$ is the usual commutator for linear operators.
Thus the Lagrange-D'Alembert principle can be reduced to
the corresponding reduced D'Alembert principle, as follows.
Find $X\in \mathfrak{gl}(V_h^s)$ subject to the constraint
\begin{equation}
  X = X_{u},
\end{equation}
for some $u\in \mathring{W}_h^r$, such that
\begin{equation}
  \delta \int_0^T \ell_h(X)\diff t = 0,
\end{equation}
for variations of the form \eqref{eq:variations gl} for all time
dependent $Y\in \mathfrak{gl}(V_h^s)$, subject to the constraint
$Y=X_{w}$ for some time dependent $w\in \mathring{W}_h^r$.  The
nonclosure of $S^r_h(V_h^s)$ under Lie brackets together with the
appearance of the Lie bracket in \eqref{eq:variations gl} is the
reason why we defined a Lagrangian $TG_h(V)$ and not just for
$\hat{\chi}$ satisfying the constraint.\footnote{In fact, it is only
necessary to define the reduced Lagrangian $\ell_h$ for $X=X_{u}$ and
$X=[X_u,X_v]$ for $u,v\in \mathring{W}_h^r$. This is how the problem
was approached for the discrete exterior calculus formulation
\citet{pavlov2011structure}. However, the finite element framework and
the bar map makes it easy enough to extend to the whole of
$\mathfrak{gl}(V_h^s)$.}
From this reduced principle, we can derive the
Euler-Poincar\`e-D'Alembert equation,
\begin{equation}
  \label{eq:EPDA}
  \left\langle \dd{}{t}\dede{l_h}{X}, Y\right\rangle
  + \left\langle \dede{l_h}{X}, [X, Y] \right\rangle = 0,
\end{equation}
for all $Y$ satisfying \eqref{eq:variations gl},
where
\begin{equation}
  \left\langle \dede{l_h}{X}, Y \right\rangle
  := \delta l_h[X; \delta X] = \lim_{\epsilon\to 0}\frac{L[X+\epsilon\delta X]
    - L[X]}{\epsilon}.
\end{equation}

For our reduced Lagrangian \eqref{eq:reduced lagrangian nonholo},
we have
\begin{equation}
  \left\langle \dede{\ell_h}{X}, \delta X\right\rangle =
  \delta \ell_h[X;\delta X] = \left\langle \delta \bar{X}, \bar{X} \right\rangle
  = \left\langle \bar{\delta X}, \bar{X} \right\rangle,
  \quad \forall \delta X \in \mathfrak{gl}(V_h^s).
\end{equation}
This means that solving \eqref{eq:EPDA} is equivalent to finding
$A_h\in S^r_h(V_h)$ such that
\begin{equation}
  \left\langle \dd{}{t}\bar{A}_h, \bar{B}_h\right\rangle_\Omega
  + \left\langle \bar{A}_h, \overline{[A_h,B_h]}\right\rangle_\Omega = 0,
\end{equation}
for all $B_h\in S^r_h(V_h)$, where $[\cdot,\cdot]$ is the
commutator bracket for linear operators. \citet{natale2018variational}
then showed that \eqref{eq:EPDA} is equivalent to finding
$u\in \mathring{W}_s$ such that
\begin{equation}
  \left\langle \pp{}{t}{u}, v \right\rangle_{\Omega} + \left\langle
  X_uu, v \right\rangle = 0, \quad \forall v \in \mathring{W}_h^r,
\end{equation}
where $X_u:\mathring{W}_h^r \to \mathring{W}_h^r$ is defined as
\begin{equation}
  \left\langle X_ua, b \right\rangle =
  \left\langle a,\nabla_h\times(u\times b) - u\nabla\cdot b \right\rangle
  + \left\llangle \{u\}, \jump{n\times (u\times b)} \right\rrangle_\Gamma.
\end{equation}
Surprisingly, given all of the complexity in the formulation, this
takes the form of a conventional finite element approximation without ever
needing to calculate $X_u$.  Some further
manipulation shows that this discretisation is in fact identical to
the centred flux discretisation described in \citet{guzman2017h},
which emerged around the same time (but without the variational
derivation).

\subsection{Compatible finite element discretisation for compressible
  fluids}
The framework was extended to compressible fluid equations in
\citet{gawlik2021variational}, which applied the programme of
\citet{desbrun2014variational,bauer2019towards} to the compatible
finite element case. This involves the introduction of advected
quantities such as temperature (which satisfies a scalar advection
equation) and density (which satisfies a continuity equation).
In the unapproximated equations, this enables us to treat Lagrangians
of the form
\begin{equation}
  L(\pp{\chi}{t},\chi) = \ell(u, a_1, a_2, \ldots, a_n),
\end{equation}
where $a_i$ are advected quantities satisfying
\begin{equation}
  \pp{}{t}a_i + \mathcal{L}_ua_i = 0,
\end{equation}
where $\mathcal{L}_u$ is a Lie derivative of an appropriate type
e.g. for advected scalars $a$,
\begin{equation}
  \mathcal{L}_ua = u\cdot\nabla a,
\end{equation}
as before, whilst for advected densities $D$,
\begin{equation}
  \mathcal{L}_uD = \nabla\cdot(uD).
\end{equation}
In that case, Hamilton's principle leads to the Euler-Poincar\'e equation
with advected quantities \citep{holm1998euler},
\begin{equation}
  \label{eq:epadv}
  \pp{\dede{l}{u}}{t} + u\cdot\nabla\dede{l}{u}
  + (\nabla u)^T\dede{l}{u} = \sum_i a_i \diamond \dede{l}{a_i},
\end{equation}
where the diamond operator $\diamond$ is defined by
\begin{equation}
  \left\langle a_i \diamond \dede{l}{a_i}, w \right\rangle
  = -\left\langle \mathcal{L}_w a_i, \dede{l}{a_i} \right\rangle,
\end{equation}
for all vector fields $w$. This allows for the relaxation to arbitrary
diffeomorphisms instead of volume preserving ones, and enables the
variational derivation of the full range of compressible fluid models
arising in geophysical fluid dynamics and beyond. In particular, the
incompressible Euler equation can be recovered by introducing a
Lagrange multiplier (the pressure) to enforce constant density $D$.

In the Koopman operator framework, discrete advected densities $D\in
V_h^s$ are treated by defining their transport equation as being dual
to that of scalar functions $f\in V_h^s$ i.e.
if $f = \hat{\chi}f_0 \in V_h^s$. 
\begin{equation}
  \left\langle D, f \right\rangle = \left\langle D_0, f_0 \right\rangle.
\end{equation}
Hence,
\begin{align}
  0 & = \pp{}{t}\left\langle D, f \right\rangle =
  \left\langle \pp{}{t}D, f \right\rangle + \left\langle D, \pp{}{t}f \right\rangle, \\
  & = \left\langle \pp{}{t}D, f \right\rangle - \left\langle D, Xf \right\rangle.
\end{align}
Therefore, we conclude that
\begin{equation}
  \left\langle \pp{}{t}D, \phi \right\rangle - \left\langle X^*D, \phi \right\rangle = 0,
  \quad \forall \phi \in V_h^s,
\end{equation}
where $X = (\pp{}{t}g)\circ g^{-1} \in \mathfrak{gl}(V_h^s)$.
When $X$ satisfies the nonholonomic constraint $X=X_{u}$ for
some $u\in {W}_h^r$, this becomes
\begin{align}
  \left\langle \pp{}{t}D, \phi \right\rangle
  &= \left\langle X_{u}^* D, \phi \right\rangle = \left\langle D, X_{u}\phi \right\rangle, \\
  &= -\left\langle \nabla_h\cdot(uD) \phi \right\rangle + \left\llangle \jump{Du}, \{\phi\}\right\rrangle_\Gamma, \\
  & = \left\langle u\cdot \nabla_h \phi, D \right\rangle
  + \int_\Gamma u\cdot \Bigg(\frac{1}{2}(D^+ n^+ + D^-n^-)(\phi^++\phi^-)
  \nonumber \\
  & \qquad\qquad\qquad\qquad\qquad\qquad
  -n^+D^+\phi^+ - n^-D^-\phi^-\Bigg)\diff S, \\ 
  & = \left\langle u\cdot \nabla_h \phi, D \right\rangle
  + \int_\Gamma u\cdot \Bigg(\frac{1}{2}(-D^+n^+\phi^+
  + D^+n^+\phi^- + \nonumber \\
  & \qquad\qquad\qquad\qquad\qquad\qquad
  D^-n^-\phi^+ -D^-n^-\phi^-\Bigg)\diff S, \\
  & = \left\langle u\cdot \nabla_h \phi, D \right\rangle
  + \int_\Gamma u\cdot \Bigg(\frac{1}{2}(-D^+n^+\phi^+
  - D^+n^-\phi^- - \\
  & \qquad\qquad\qquad\qquad\qquad\qquad
  D^-n^+\phi^+ -D^-n^-\phi^-\Bigg)\diff S, \\
  & = \left\langle u\cdot \nabla_h \phi, D \right\rangle
  - \left\llangle \jump{u\phi}, \{D\} \right\rrangle_\Gamma,
\end{align}
which is the standard Discontinuous Galerkin centred flux scheme for
the continuity equation for advected densities.

Similarly, by considering $\delta \left\langle D, f\right\rangle=0$, we obtain
\begin{equation}
  \delta D - Y^*D = 0,
\end{equation}
where $Y$ is such that $\delta \hat{\chi} + Y\hat{\chi}=0$.

We can again use the overbar map to build Lagrangians with advected
quantities defined on the whole of $GL(V_h^s,V_h^s)$. For example,
the shallow water equations (with flat topography) has the (reduced) Lagrangian
\begin{equation}
  \ell = \int_\Omega D\frac{|u|^2}{2} - \frac{gD^2}{2}\diff x.
\end{equation}
The discrete Lagrangian can then be written as
\begin{align}
  L(\pp{}{t}\hat{\chi},\hat{\chi}) = \int_\Omega (\hat{\chi}^*D_0)\frac{|\overline{\pp{}{t}\hat{\chi}
      \circ \hat{\chi}^{-1}}|^2}{2} - \frac{g(\hat{\chi}^*D_0)^2}{2}\diff x.
\end{align}
This has discrete reduced Lagrangian
\begin{equation}
  \ell_h[D, X] = \int_\Omega D\frac{|\overline{X}|^2}{2} - \frac{gD^2}{2}
  \diff x.
\end{equation}

The reduced Hamilton's principle with nonholonomic constraints becomes
\begin{equation}
  \delta \int_0^T \ell_h(X, D) \diff t = 0, 
\end{equation}
with variations
\begin{equation}
  \delta X = \pp{}{t} Y + [X, Y], \quad
  \delta D = Y^*,
\end{equation}
and constraints $X=X_u$ for some $u\in {W}_h^r$, and $Y=Y_w$
for all $w\in {W}_h^r$.

This gives the equation
\begin{equation}
  \left\langle \dd{}{t}\dede{l}{X_u}, X_w \right\rangle
  - \left\langle \dede{l}{X_u}, [X, X_w]\right\rangle + \left\langle \dede{l}{X)_u}, X_w^*D \right\rangle,
  \quad \forall w \in {W}_h^r,
\end{equation}
which can be written as
\begin{equation}
  \dd{}{t}\dede{l}{X_u} + \ad^*_{X_u}\dede{l}{X_u} - \dede{l}{D}\diamond D
  \in (S^r_h)^*,
\end{equation}
where 
\begin{equation}
  \left\langle \ad^*_XY, Z \right\rangle = -\left\langle Y, [X,Z]\right\rangle,\,
  -\left\langle X\diamond Y, Z \right\rangle = \left\langle X, Y^*Z\right\rangle,
\end{equation}
  and $(S^r_h)^*$ is the dual space to $S^r_h$ in $\mathfrak{gl}(V_h)$.
This is the discrete analogue of \eqref{eq:epadv}.

\subsection{Compatible finite element variational discretisations: discrete conservation laws}
The goal of variational discretisations is to derive numerical methods
that provide discrete analogues of conservation laws of the
unapproximated equations. Noether's theorem can still be applied to
nonholonomic variational principles without change, provided that the
constraints are invariant under the relevant symmetry as well as the
Lagrangian. The discrete varational principles discussed here
are invariant under time translation, leading to the conservation of
energy as usual. This is straightforward to check directly 
in the incompressible case of \citet{natale2018variational} since
the energy is
\begin{equation}
  E = \left\langle \bar{A}_h, \bar{A}_h \right\rangle,
\end{equation}
so the energy equation is obtained in \eqref{eq:EPDA} by
taking $\bar{B}_h=\bar{A}_h$ and using antisymmetry of the bracket.

Regarding Kelvin's circulation theorem, \citet{natale2018variational}
demonstrated the same ``echo'' in \eqref{eq:disc kelv} for the
compatible finite element case, which also emerges from 
\citet{gawlik2021variational} in the case where the Lagrangian
depends only on velocity and density, as expected. In the incompressible
case, the conservation law takes the form
\begin{equation}
  \dd{}{t}\left\langle u, \overline{\hat{\chi}X_c\hat{\chi}^{-1}} \right\rangle = 0,
\end{equation}
for all time independent $c\in \mathring{W}_h^r$.

\citet{gawlik2021variational} further extended the framework by
introducing advected tracers (which can represent potential
temperature, and salinity in the ocean), which are discretised as
elements $\theta$ of $V_h$ that can be acted on by velocity via
$X_{u}\theta$ as above. This unlocks the variational discretisation
of all of the main variational models of geophysical fluid dynamics.
The framework was further extended in \citet{gawlik2021structure} to
accommodate advected transported fluxes represented in
${W}^r_h$, leading to variational discretisations of
magnetohydrodynamics.

\subsection{Variational time integrators}
So far in the this section we have only discussed the variational
discretisation in space, leading to a system that is still continuous
in time. As set out in the original vision of
\citet{pavlov2011structure}, the idea is to also discretise Hamilton's
principle in time. Variational integrators arising from time
discretisation in Hamilton's principle have a long history, with the
programme being formally set out in \citet{marsden2001discrete}.
Variational time integration using a finite difference discretisation
in time was investigated in \citet{gawlik2021variational}. Since the
kinetic energy is a nonlinear function of both the coordinate in $V_h$
and its rate of change, any variational integrator will result in a
system that requires the solution of an implicit nonlinear system
(using Newton's method, for example). This contrasts with the case of
classical mechanics, where the Lagrangian $L(z,\pp{}{t}{z})$ splits into a
kinetic energy depending only on $\pp{}{t}{z}$ and a potential energy
depending only on $z$. This split makes explicit variational
integrators possible in that case. Following
\citet{pavlov2011structure}, \citet{gawlik2021variational} made the
choice of replacing $A(t)=\pp{}{t}{\hat{\chi}}\hat{\chi}^{-1}$ with $A_k =
\tau^{-1}(\hat{\chi}_{k+1}\hat{\chi}^{-1}_k)/\Delta t$, where $\tau:\mathfrak{gl}(V_h)
\to G_(V_h)$ is the Cayley transform
\begin{equation}
  \tau(A) = \left(I - \frac{A}{2}\right)^{-1}
  \left(I + \frac{A}{2}\right),
\end{equation}
with other possibilities for $\tau$ being discussed in
\citet{bou2009hamilton}. After this replacement, the variational
integrator is derived by finding the stationary point of the resulting
discrete action principle depending on $\hat{\chi}_{k+1}$ and $\hat{\chi}_k$. It was
found through numerical experiments that the resulting scheme is
only conditional stable, requiring condition $\delta t <
Ch$. This is disappointing given that intensive computation is
required to advance the solution by a small step. The nature of this
stepsize requirement is an open problem in the area, as is the
question of whether a variational integrator can be found that allows
larger timesteps.  \citet{natale2018variational} used the implicit
midpoint rule to discretise the semidiscrete variational
discretisation. The implicit midpoint rule is not known to be a
variational integrator for such systems, but does have the property
that it preserves any quadratic invariants of the time continuous
system, which includes the energy in the incompressible case.

\section{Almost-Poisson brackets}
\label{sec:poisson}
An alternative, but related, route to structure preserving
discretisations is found via Poisson bracket formulations.

Poisson brackets are bilinear, skew symmetric maps that take pairs of
functionals on some space where the solutions of the PDE reside, which also
satisfy the Jacobi identity,
\begin{equation}
  \{A, \{B, C\}\} + \{B, \{C, A\}\} + \{C, \{A, B\}\} = 0,
\end{equation}
for all functionals $A,B,C$. In the process of discretisation of
Poisson brackets for fluid dynamics, the Jacobi identity is lost, for
similar reasons that the nonholonomic constraints are required for the
variational discretisations of Section \ref{sec:variational}. Poisson
brackets without the Jacobi identity are called ``almost Poisson
brackets'' but we shall just use the term Poisson bracket in the rest
of this article for brevity.  Poisson brackets can have special
functionals called Casimirs, which make the Poisson bracket vanish,
i.e. $C$ is a Casimir if $\{F,C\}=0$ for all functionals $F$. Then,
Casimirs are conserved by the dynamics since
\begin{equation}
  \pp{}{t}{C} = -\{C,H\}=0.
\end{equation}
The goal of building discretisations using Poisson brackets is that
they automatically conserve the Hamiltonian, and if any Casimirs
survive the discretisation process then they will be conserved as
well. We shall discuss specific examples later.

The aim of building structure preserving discretisations for fluid
PDEs using Poisson brackets was introduced in
\citet{morrison1982poisson,salmon1983practical}, although in fact the
energy and enstrophy preserving discretisations in
\citet{arakawa1966computational} are the first instance of a Poisson
bracket discretisation (but not presented that way), which was
extended to the rotating shallow water equations and beyond in
subsequent work by Arakawa, Sadourny and others
\citep{sadourny1968integration,sadourny1972conservative,sadourny1975dynamics,arakawa1977computational,arakawa1981potential,arakawa1990energy}).

In the 2000s, when attention was focussed much more on triangular and
polygonal grids to provide a more uniform coverage of the sphere, the
idea of using Poisson brackets was revived to produce energy
conserving schemes (or to at least guide the design of practical
schemes that are as energy consistent as possible). This took place in
a number of groups
\citep{ringler2010unified,skamarock2012multiscale,eldred2017total,gassmann2008towards,gassmann2013global,tort2015energy,dubos2015dynamico}. There
was also work on extending Poisson brackets to Nambu brackets in
pursuit of additional conserved quantities in the method
\citep{sommer2009conservative}.

There have also been some interesting studies about the relevance of
conservation for geophysical models. \citet{dubinkina2007statistical}
demonstrated the benefits of using energy-enstrophy conserving schemes
to obtain correct statistical equilibria, and
\citet{thuburn2014cascades} demonstrated that energy conservation is
important to obtain realistic backscatter in underresolved simulations
of two dimensional turbulence (provided that enstrophy is dissipated
at the small scale). \citet{dubinkina2018relevance} demonstrated that
conserving both energy and enstrophy is important in the context of
data assimilation. Even when structure preserving discretisations lead to
systems of equations that are challenging to solve efficiently, it is
useful to consider how they relate to more standard discretisations to
see where conservation errors are being committed, and to see when
they are likely to be large or small.

In this section, we describe how Poisson brackets can be used to
construct energy (and enstrophy) conserving schemes. This work has
been heavily informed by previous works using finite difference
methods, such as those cited above. We shall start by briefly
discussing the 2D incompressible Euler equations, the the rotating
shallow water equations and then vertical slice and three dimensional
models. We initially assume that we are solving the equations on a
closed manifold (the surface of the sphere, or periodic boundary
conditions, for example), and will return to the treatment of boundary
conditions later.

\localtableofcontents

\subsection{Incompressible Euler equations: continuous theory}
Returning to Arnold's variational formulation for incompressible flow
in \eqref{eq:EP}, we can reformulate as
\begin{equation}
  \label{eq:arnold to poisson}
  \left\langle \pp{m}{t}, v \right\rangle
  + \left\langle \left[
    \dede{H}{m}, v
    \right], m \right\rangle = 0,
  \, \forall v \in \mathfrak{X}_{\vol}(\Omega),
\end{equation}
where
\begin{equation}
  H(m) = \left\langle m, u \right\rangle - \ell(u),
\end{equation}
and we take $m=\dede{l}{u}\in \mathfrak{X}_{\vol}(\Omega)$, now
inverting the relationship so that $u$ is considered as an operator
applied to $m$. Here we use the variational derivative $\dede{F}{u}\in
\mathfrak{X}_{\vol}(\Omega)$ defined by
\begin{equation}
  \left\langle \dede{F}{u}, v\right\rangle = \lim_{\epsilon\to 0}
  \frac{1}{\epsilon}\left(F[u+\epsilon v] - F[u]\right).
\end{equation}

Equation \eqref{eq:arnold to poisson} is equivalent to the Poisson
formulation
\begin{equation}
  \pp{}{t}{F}[m] + \left\{F, H\right\} = 0,
\end{equation}
for arbitrary functionals $F:\mathfrak{X}_{\vol}(\Omega)\to \mathbb{R}$,
with Poisson bracket
\begin{equation}
  \label{eq:m bracket}
  \left\{F,G\right\} =
  \int_\Omega \left[\dede{G}{m},\dede{F}{m}\right]\cdot m
    \diff x.
\end{equation}
Since the Poisson bracket is antisymmetric, this leads immediately to
conservation of the Hamiltonian,
\begin{equation}
  \dd{}{t}{H} = -\{H,H\}+\pp{H}{t}=0,
\end{equation}
provided that $H$ has no explicit dependence on time $t$.

For the 2D incompressible Euler equations, a clear path towards
deriving enstrophy conservation requires modification of 
the bracket \eqref{eq:m bracket} by changing variables $m\to u$.
This produces the equivalent Poisson formulation
\begin{equation}
  \label{eq:u LP}
  \pp{}{t}{F}[u] + \left\{F, H\right\}[u] = 0,
\end{equation}
for all functionals $F:\mathfrak{X}_{\vol}(\Omega)\to \mathbb{R}$,
where
\begin{equation}
  \label{eq:omega bracket}
  \{F,G\} = \int_\Omega \omega \left(\dede{F}{u}\right)
  \cdot \dede{G}{u}^{\perp} \diff x,
\end{equation}
$\omega=\nabla^\perp\cdot u$,
and
\begin{equation}
  \label{eq:incompressible H}
  H = \frac{1}{2}\int_\Omega |u|^2\diff x.
\end{equation}
The derivation of this bracket formulation from \eqref{eq:m bracket}
is discussed in \citet{morrison1982poisson,marsden1983coadjoint}.
Here, we just directly demonstrate that it leads to the incompressible
Euler equations by computation. First, we compute $\dede{H}{u}=u$.
For a linear functional $F[u]=\left\langle
u, w \right\rangle$ with $w\in \mathfrak{X}_{\vol}(\Omega)$,
\eqref{eq:u LP} becomes
\begin{equation}
  \left\langle u_t, w \right\rangle + \left\langle \omega, w\cdot u^\perp \right\rangle = 0,
\quad \forall w \in \mathfrak{X}_{\vol}(\Omega),
\end{equation}
which is a weak formulation of the equation
\begin{equation}
  \label{eq:weak 2d euler}
  u_t + \omega u^\perp + \nabla P = 0,
\end{equation}
where $P$ is some potential chosen so that $u_t\in
\mathfrak{X}_{\vol}(\Omega)$ (because we only test against functions
from $\mathfrak{X}_{\vol}(\Omega)$ in \eqref{eq:weak 2d euler}, so the
equation is projected into $\mathfrak{X}_{\vol}(\Omega)$). Writing
$P=\frac{1}{2}|u|^2+p$, we obtain
\begin{equation}
  u_t + u^\perp \nabla^\perp\cdot u + \nabla \frac{1}{2}|u|^2 + \nabla p = 0,
\end{equation}
which becomes recognisable as the incompressible Euler equations
after recalling the identity
\begin{equation}
  (u\cdot \nabla)u = \omega u^\perp + \frac{1}{2}\nabla |u|^2.
\end{equation}
Returning to the Poisson bracket \eqref{eq:omega bracket}, we find
that it has an infinite number of Casimirs of the form
\begin{equation}
  C_n[u] = \int_\Omega \omega^n \diff x, \quad n=1,2,\ldots
\end{equation}
To verify that $C_n$ is a Casimir, we compute
\begin{align}
  \left\langle \dede{C_n}{u}, v\right\rangle &= n\int_\Omega \omega^{n-1}
  \nabla^\perp v \diff x, \quad \forall v \in \mathfrak{X}_{\vol}(\Omega), \\
  & = \int_\Omega\nabla^\perp(-n\omega^{n-1})\cdot v \diff x,
  \label{eq:dCn}
\end{align}
having integrated by parts (we assume for now that there are now
boundaries and so may ignore the surface term). Hence, we conclude
formally that
\begin{equation}
\dede{C_n}{u} = -n\nabla^\perp\omega^{n-1}.
\end{equation}
Inserting into the Poisson bracket then gives
\begin{align}
  \label{eq:euler casimir calc}
  \left\{F,C_n\right\}
  &= -n\int_\Omega\omega\dede{F}{u}\cdot \left(\nabla^\perp\omega^{n-1}\right)^{\perp}\diff x, \\
  &= n\int_\Omega\omega\dede{F}{u}\cdot \nabla\omega^{n-1}\diff x, \\
  & = (n-1)\int_\Omega \dede{F}{u}\cdot \nabla\omega^n \diff x, \\
  & = -(n-1)\int_\Omega \omega^n\underbrace{\nabla\cdot\dede{F}{u}}_{=0}
  \diff x = 0,
  \label{eq:euler casimir calc end}
\end{align}
having integrated by parts again, for any functional $F$ on
$\mathfrak{X}_{\vol}(\Omega)$. We have $\nabla\cdot\dede{F}{u}=0$
since $\dede{F}{u}\in \mathfrak{X}_{\vol}(\Omega)$. Hence, $C_n$ is a
conserved quantity for the Poisson dynamics from any Hamiltonian.  In
particular, $C_1$ is the \emph{total vorticity}, whilst $C_2$ is the
\emph{enstrophy}, both of which provide strong constraints on two
dimensional incompressible turbulence.

\subsection{Incompressible Euler equations: compatible finite element discretisation}
To make our compatible finite element discretisation, we restrict $u$
and $w$ to the divergence-free subspace $\zeta_h$ of some chosen H(div)
finite element space $\mathbb{V}_h^1$ (such as Raviart-Thomas or Brezzi-Douglas-Marini
on triangles). We have to make a further approximation since
$\nabla^\perp\cdot u$ is not defined for H(div) spaces, and so we
define $\omega_h \in \mathbb{V}_h^0$, such that
\begin{equation}
  \label{eq:omega_h}
  \left\langle \gamma, \omega_h \right\rangle = -\left\langle \nabla^\perp\gamma, u \right\rangle,
  \quad \forall \gamma \in \mathbb{V}_h^0,
\end{equation}
i.e. $\omega_h$ is defined from the approximated weak curl of $u$, the
dual of the $\nabla^\perp$ operator restricted to $\mathbb{V}_h^0$ and $\mathbb{V}_h^1$.
Note that this is where we have used the absence of boundary,
otherwise there would be a boundary term causing complications that we
shall discuss later. Having defined $\omega_h$, we write the discrete
Poisson bracket as
\begin{equation}
  \label{eq:omega incompressible bracket}
  \left\{
  F,G
  \right\} = \int_\Omega \omega_h \dede{F}{u}\cdot \dede{G}{u}^\perp
  \diff x,
\end{equation}
where $F$ and $G$ are now functionals on $\mathbb{V}_h^1$. With the same Hamiltonian
\eqref{eq:incompressible H} now restricted to $\mathbb{V}_h^1$, the Poisson
formulation implies the following dynamics for $u\in \zeta_h$,
\begin{equation}
  \label{eq:u zeta_h}
  \left\langle w, u_t \right\rangle + \left\langle \omega_h w\cdot u^\perp \right\rangle = 0,
  \quad \forall w\in \zeta_h.
\end{equation}
The discrete Hamiltonian is conserved as usual through the antisymmetry
of the bracket and the time-independence of $H$. Concerning Casimirs of
the bracket, we can repeat the earlier calculation for functionals $C_n$
computed with $\omega_h$ substituted for $\omega$,
i.e.
\begin{equation}
C_{n,h}[u] = \int_\Omega \omega_h^n \diff x.
\end{equation}
Following \eqref{eq:dCn}, we obtain
\begin{equation}
  \left\langle \dede{C_{n,h}}{v}\right\rangle
  = \int_\Omega \nabla^\perp(-n\omega^{n-1})\cdot v \diff x,
  \quad \forall v \in \zeta_h,
\end{equation}
but now can only conclude that
$\dede{C_{n,h}}{v}=P_1\left(-n\nabla^\perp\omega^{n-1}\right)$, where
$P_1$ is the $L^2$ projection into $\zeta_h$, which prevents us from
showing that $C_{n,h}$ is a Casimir for $n>2$. However, when $n=1$, we
obtain that $\dede{C_{1,h}}{v}=0$ (so $C_{1,h}$ is trivially conserved,
just as $C_1$ is for the unapproximated case). When $n=2$, we have
$\nabla^\perp\omega\in \mathbb{V}_h^1$ by the embedding property of the discrete
de Rham complex, and then a calculation identical to (\ref{eq:euler
  casimir calc}-\ref{eq:euler casimir calc end}) shows that
$\{C_{2,h},G\}=0$ for any functional $G$ on $\zeta_h$, and hence the
numerical enstrophy $C_{2,h}$ is a Casimir and is conserved for
dynamics generated from any Hamiltonian. Hence, this scheme conserves
energy, total vorticity, and enstrophy.

To make a practical implementation of the scheme, one can follow
two approaches. The main hurdle is that the scheme is defined on
$\zeta_h$, and not the whole of $\mathbb{V}_h^1$. Since $\nabla\cdot$
maps from $\mathbb{V}_h^1$ onto $\mathbb{V}_h^2$, the divergence free subspace $\zeta_h$
is equivalently represented as
\begin{equation}
  \zeta_h = \left\{u\in \mathbb{V}_h^1:\int_\Omega \phi\nabla \cdot u \diff x
  = 0,\quad\forall \phi\in\mathbb{V}_h^2\right\}.
\end{equation}
Hence, we can equivalently write the following system
for $(\omega,u,P)\in \mathbb{V}_h^0\times \mathbb{V}_h^1\times \mathbb{V}_h^2$ such that
\begin{align}
  \left\langle \gamma, \omega \right\rangle - \left\langle \nabla^\perp\gamma,
  u \right\rangle & = 0, \quad \forall \gamma \in \mathbb{V}_h^0, \\
  \left\langle w, u_t \right\rangle + \left\langle w, \omega u^\perp\right\rangle
  - \left\langle \nabla \cdot w, P \right\rangle & = 0, \quad \forall w \in \mathbb{V}_h^1, \\
  \left\langle \phi, \nabla\cdot u \right\rangle & = 0, \quad \forall \phi \in \mathbb{V}_h^2.
\end{align}
It can easily be checked that the solution satisfies $\nabla\cdot u =
0$ in $L^2$. Selecting $w\in \zeta_h\subset \mathbb{V}_h^1$ makes the $P$ term
disappear and we recover \eqref{eq:u zeta_h}. This formulation builds
a bridge to the shallow water and compressible systems that we shall
look at later.

On the other hand, for $u\in \zeta_h$ we can directly parameterise
$u=\nabla^\perp\psi$ for $\psi\in \mathbb{V}_h^0$ and choose $w=\nabla^\perp\beta$
in \eqref{eq:u zeta_h} for all $\beta\in \mathbb{V}_h^0$
and we obtain
\begin{align}
  \left\langle \gamma, \omega \right\rangle + \left\langle \nabla\gamma,\nabla \psi\right\rangle
  & = 0, \quad \forall \gamma \in \mathbb{V}_h^0, \\
  \left\langle \nabla \beta, \nabla \psi_t\right\rangle
  + \left\langle \nabla \beta, \omega \nabla^\perp \psi \right\rangle & = 0,
  \quad \forall \beta \in \mathbb{V}_h^0,
\end{align}
which is a discretisation of the incompressible Euler equation in
vorticity streamfunction form,
\begin{equation}
  \omega_t + \nabla\cdot(\omega\nabla^\perp\psi) = 0,
  \quad -\nabla^2\psi = \omega.
\end{equation}
This formulation was presented in \citet{liu2001simple}, where it was
analysed in the viscous case, but the energy conservation of the
inviscid equations did form an important part of the proof. A related
scheme was presented in \citet{liu2000high}, but with the vorticity in
a discontinuous space, with appropriate jump terms defining the fluxes
between cells. With an average flux, the energy and enstrophy are both
conserved.  If a Lax Friedrichs flux is used, then energy is still
conserved but enstrophy is dissipated (as is appropriate for long time
simulations of cascading 2D vortex dynamics). It is possible to modify
the Poisson bracket so that enstrophy is dissipated in this way, with
the antisymmetric formulation still conserving energy, as we shall see
later. \citet{liu2000high} proved convergence for both types of fluxes.
This formulation and analysis was extended by \citet{bernsen2006dis}
to the quasigeostrophic model of large scale rotating geophysical fluid
dynamics, including the case with islands in the flow (considering
the fluid as an ocean).

To extend the conservation properties to a full discrete method, one
can use the implicit midpoint rule, which conserves all quadratic
invariants of the continuous time system. This includes energy and
enstrophy (or just energy where enstrophy is dissipated) in the
incompressible Euler case.

\subsection{Rotating shallow water equations: continuous theory}
Moving on to the rotating shallow water equations, it is tempting to
continue working with the vorticity-streamfunction formulation above.
However, this places a limitation on the possibilities of extension
to three dimensional models, complicates the boundary conditions,
and is not preferred by practitioners since the prognostic variables
are not quantities that are directly measurable. Hence, we must
address the challenge of finding a compatible finite element
discretisation of the rotating shallow water equations using
velocity $u$ and layer depth $D$.

We start from the Lie-Poisson formulation of the rotating shallow
water equations given by
\begin{align}
  H &= \int_\Omega \frac{1}{2D}|m|^2 + gD\left(\frac{D}{2}+b\right)
  \diff x, \\
  \{F,G\} & = \left\langle [\dede{G}{m}, \dede{F}{m}], m\right\rangle
  + \left\langle \dede{F}{D},\nabla\cdot(D\dede{G}{m}) \right\rangle
  - \left\langle \dede{G}{D},\nabla\cdot(D\dede{F}{m}) \right\rangle,
\end{align}
where $m$ is defined by
\begin{equation}
  \int_\Omega m\cdot v \diff x = \int_\Omega D\left({u} +
  R\right)\cdot v\diff x, \quad \forall v\in \mathfrak{X}(\Omega),
\end{equation}
where $\nabla^\perp\cdot R=f$, the Coriolis parameter. These equations
emerge from the reduced Hamilton's principle with advected density $D$
after applying the Legendre transform as described in
\citet{holm1998euler}. Following the incompressible case, if we want
enstrophy conservation to emerge then we need to change variables
to $(u,D)$, this leads to the equivalent Poisson bracket formulation
\begin{align}
  \label{eq:q H}
  H &= \frac{1}{2}\int_\Omega D\|u\|^2 + gD\left(\frac{D}{2}+b\right) \diff x, \\
  \{F, G\} & = \left\langle q,\dede{F}{u}\cdot\dede{G}{u}^\perp\right\rangle
  - \left\langle \nabla\cdot\dede{F}{u}, \dede{G}{D}\right\rangle
  + \left\langle \nabla\cdot\dede{G}{u}, \dede{F}{D}\right\rangle,
  \label{eq:q bracket}
\end{align}
where $q = \frac{\nabla^\perp\cdot u + f}{D}$ is the potential
vorticity. After computing the variational derivatives
\begin{align}
  \dede{H}{u} & = Du, \\
  \dede{H}{D} & = \frac{1}{2}\|u\|^2 + g(D+b),
\end{align}
and substituting into the Poisson bracket equation $F_t + \{F,H\}$,
we formally obtain
\begin{align}
  \label{eq:ut}
  u_t + qDu^{\perp} + \nabla\left(\frac{1}{2}|u|^2 + g(D+b)\right) & = 0,\\
  D_t + \nabla\cdot(uD) & = 0,
\end{align}
which we recognise as the rotating shallow water equations in vector
invariant form. We define functionals $C_n$ (which will turn out to be
Casimirs) by
\begin{equation}
  C_n[u,D] = \int_\Omega Dq^n \diff x.
\end{equation}
To compute the variational derivatives, given $\epsilon>0$ and
$(v,\phi)\in H(\ddiv)\times L^2$, writing $D_\epsilon = D+\epsilon \phi$,
$u_\epsilon=u+\epsilon v$, we define $q_\epsilon \in H^1$ such that
\begin{equation}
  \label{eq:q eps}
  \left\langle \gamma, D_\epsilon q_\epsilon\right\rangle =
  -\left\langle \nabla^{\perp}\gamma, u_\epsilon \right\rangle + \left\langle \gamma, f\right\rangle,
  \quad \forall \gamma \in H^1,
\end{equation}
noting that $q_\epsilon|_{\epsilon=0}=q$ (after integrating by parts
and using the lack of surface term in a domain without
boundary). Then,
\begin{align}
  \left\langle \dede{C_n}{v}, v \right\rangle
  + \left\langle \dede{C_n}{D}, \phi\right\rangle
  & = \lim_{\epsilon\to 0}
  \frac{1}{\epsilon}
  \left(
  C_n[u+\epsilon v, D + \epsilon \phi] - C_n[u,D]
  \right), \\
  & = \int_\Omega \lim_{\epsilon \to 0}\frac{D_\epsilon q_\epsilon^n - Dq^n}{\epsilon}
  \diff x, \\
  &= \int_\Omega \dd{}{\epsilon}|_{\epsilon=0}D_\epsilon q_\epsilon^n
  \diff x, \\ \nonumber
  & = \left\langle q_\epsilon^{n-1}|_{\epsilon=0}, \dd{}{\epsilon}|_{\epsilon=0}D_\epsilon q_\epsilon \right\rangle \\
& \qquad  + \left\langle (n-1)q_\epsilon^{n-2}|_{\epsilon=0}\dd{}{\epsilon}|_{\epsilon=0} q_\epsilon,
  D_\epsilon q_\epsilon \right\rangle, \\
  & \nonumber = n\left\langle q_\epsilon^{n-1}|_{\epsilon=0}, \dd{}{\epsilon}|_{\epsilon=0}D_\epsilon q_\epsilon \right\rangle\\
& \qquad  - \left\langle (n-1)q_\epsilon^n|_{\epsilon=0}, \dd{}{\epsilon}|_{\epsilon=0} D_\epsilon) \right\rangle,
\end{align}
From \eqref{eq:q eps} we have
\begin{equation}
  \left\langle \gamma, \dd{}{\epsilon}|_{\epsilon=0} D_\epsilon q_\epsilon\right\rangle =
  -\left\langle \nabla^{\perp}\gamma, v \right\rangle,
  \quad \forall \gamma \in H^1,
\end{equation}
(correcting a typographic error in Equation 58 of \citet{mcrae2014energy}),
and hence we have
\begin{align}
  \label{eq:swe C_n derivs}
  \left\langle \dede{C_n}{v}, v \right\rangle
  + \left\langle \dede{C_n}{D}, \phi\right\rangle
 & = n\left\langle \nabla^{\perp}q^{n-1}, v \right\rangle
  - \left\langle (n-1)q^n, \phi \right\rangle,
\end{align}
i.e.
\begin{equation}
  \dede{C_n}{v} = n\nabla^{\perp}q^{n-1}, \quad
  \dede{C_n}{D} = -(n-1)q^n.
\end{equation}
Then, inserting into the Poisson bracket \eqref{eq:q bracket},
we obtain
\begin{align}
  \{C_n, G\} & = \left\langle
  \underbrace{nq\nabla^\perp q^{n-1}}_{=(n-1)\nabla^{\perp} q^n}, \dede{G}{u}^\perp\right\rangle \nonumber \\
&\qquad  - \left\langle \underbrace{\nabla\cdot n\nabla^{\perp}q^{n-1}}_{=0},
  \dede{G}{D} \right\rangle
  + \left\langle (n-1)q^n, \nabla\cdot \dede{G}{u} \right\rangle, \\
  & = \left\langle (n-1) \nabla q^n, \dede{G}{u} \right\rangle
  - \left\langle (n-1) q^n, \nabla\cdot\dede{G}{u} \right\rangle = 0,
  \label{eq:Cn in bracket}
\end{align}
where the last line is obtained by integrating by parts, and hence
$C_n$ is a Casimir of the bracket \eqref{eq:q bracket}, and so it is
a conserved quantity of the dynamics for any Hamiltonian. Another, more
simple, Casimir is the mass
\begin{equation}
  M[u,D] = \int_\Omega D \diff x,
\end{equation}
with variational derivatives
\begin{equation}
  \dede{M}{u}=0, \quad \dede{M}{D}=1,
\end{equation}
and hence
\begin{equation}
  \{M, G\} = \left\langle \nabla\cdot \dede{G}{u}, 1 \right\rangle
  = \int_{\Omega} \nabla \cdot \dede{G}{D} \diff x = 0,
\end{equation}
using the divergence theorem (assuming no boundary currently),
so mass is also a Casimir of the bracket and is conserved for
any Hamiltonian.

\subsection{Rotating shallow water equations: Poisson bracket discretisation}
To produce a discretisation of this structure, \citet{mcrae2014energy}
simply took the structure (\ref{eq:q H}-\ref{eq:q bracket}), restricted
$(u,D)$ to $\mathbb{V}_h^1\times \mathbb{V}_h^2$, and replaced $q$ with the discrete
approximation $q\in \mathbb{V}_h^0$, with
\begin{equation}
  \label{eq:discrete q}
  \left\langle \gamma, D q\right\rangle =
  -\left\langle \nabla^{\perp}\gamma, u \right\rangle + \left\langle \gamma, f\right\rangle,
  \quad \forall \gamma \in \mathbb{V}_h^0,
\end{equation}
following the exposition for the discretisation for incompressible
Euler equations above. Then, we obtain
\begin{align}
  \left\langle w, \dede{H}{u} - Du\right\rangle & = 0, \, \forall w\in \mathbb{V}_h^1, \\
  \left\langle \phi, \dede{H}{D} - \frac{1}{2}|u|^2 - g(D+b) \right\rangle
  &= 0,
  \quad \forall \phi\in \mathbb{V}_h^2,
\end{align}
i.e.,
\begin{equation}
  \dede{H}{u} = P_1(Du), \quad \dede{H}{D} = P_2
  \left(\frac{1}{2}|u|^2 + g(D+b)\right),
\end{equation}
where $P_1$ and $P_2$ are the $L^2$ projections into $\mathbb{V}_h^1$ and $\mathbb{V}_h^2$,
respectively. To derive the equations of motion, we take
$F[u,D] = \left\langle w, u \right\rangle + \left\langle \phi, D \right\rangle$ for $w,\phi\in \mathbb{V}_h^1\times \mathbb{V}_h^2$, so that $\dede{F}{u}=w$, $\dede{F}{D}=\phi$, and substitute
into the Poisson dynamics to obtain
\begin{align}
  \left\langle w, u_t \right\rangle + \left\langle \phi, D_t\right\rangle & = \pp{}{t}{F} \\
  & = -\{F, H\}, \\
  \nonumber  & = -\left\langle q, w\cdot P_1(Du)^\perp \right\rangle +
 \left\langle \nabla\cdot w, P_2\left(\frac{1}{2}|u|^2 
 + g(D+b)\right)\right\rangle \\
 & \qquad - \left\langle \nabla\cdot P_1(Du), \phi \right\rangle,\\
\nonumber  & = -\left\langle q, w\cdot P_1(Du)^\perp \right\rangle + \left\langle \nabla\cdot w, \frac{1}{2}|u|^2 + gD\right\rangle - \left\langle \nabla\cdot P_1(Du), \phi \right\rangle, \\
& \qquad\qquad  \quad \forall w,\phi \in \mathbb{V}_h^1\times \mathbb{V}_h^2,
\end{align}
where we were able to drop the $P_2$ since the result of the
projection was in an $L^2$ inner product with $\nabla\cdot w\in P_2$
(so the discrete de Rham complex is crucial here). Writing $m=P_1(Du)$
we put everything together as $(u,D,q,m)\in \mathbb{V}_h^1\times \mathbb{V}_h^2\times
\mathbb{V}_h^0\times \mathbb{V}_h^1$, such that
\begin{align}
  \label{eq:swe mcrae u}
  \left\langle w, u_t \right\rangle + \left\langle qw, m^\perp \right\rangle
  - \left\langle \nabla\cdot w, \frac{1}{2}|u|^2 + g(D+b) \right\rangle & = 0,\,
  \quad \forall w \in \mathbb{V}_h^1, \\
  \label{eq:swe mcrae D}
  \left\langle \phi, D_t + \nabla\cdot m \right\rangle & = 0,\,
  \forall \phi \in \mathbb{V}_h^2, \\
  \label{eq:swe mcrae q}
  \left\langle \gamma, qD \right\rangle + \left\langle \nabla^\perp \gamma, u\right\rangle
  - \left\langle \gamma, f \right\rangle & = 0, \,\forall \gamma \in \mathbb{V}_h^0, \\
  \left\langle v, m - uD \right\rangle & = 0, \,\forall v \in \mathbb{V}_h^1.
  \label{eq:swe mcrae m}
\end{align}
This is a set of coupled equations which must be solved together,
but since Equations (\ref{eq:swe mcrae q}-\ref{eq:swe mcrae m}) do
not contain time derivatives, $q$ and $m$ may be reconstructed at
any time from $u$ and $D$. Since the equations are derived from
a Poisson bracket formulation, we can immediately deduce that they
conserve the Hamiltonian. Concerning the Casimirs, we make
the same calculations for mass $M$  as for the undiscretised
case, obtaining
\begin{align}
  \dede{M}{u}=0, \quad \dede{M}{D}=P_1(D) = 1.
\end{align}
For $C_n$, we can only make use of \eqref{eq:q eps} when $n=1$ (so
that $\gamma=1$) or $n=2$ (so that $\gamma=q_\epsilon$), leading
to
\begin{equation}
\dede{C_n}{u}=n\left(\nabla^\perp q^{n-1}\right), \quad
\dede{C_n}{D}=-(n-1)P_2\left(q^n\right).
\end{equation}
\eqref{eq:Cn in bracket} then follows but only for $n=1,2$,
\begin{align}
  \nonumber
  \{C_n, G\} & = \left\langle
     nq \nabla^\perp q^{n-1}, \dede{G}{u}^\perp\right\rangle \\
 &    \qquad
  - \left\langle n\underbrace{\nabla\cdot\nabla^{\perp}q^{n-1}}_{=0},
  \dede{G}{u} \right\rangle
  + \left\langle P_2((n-1)q^n), \nabla\cdot \dede{G}{u} \right\rangle,  \\
  & = \left\langle
  nq \nabla^\perp q^{n-1}, \dede{G}{u}^\perp\right\rangle
  + \left\langle (n-1)q^n, \nabla\cdot \dede{G}{u} \right\rangle = 0,
\end{align}
where we may integrate by parts since $\dede{G}{u}\in
\mathbb{V}_h^1\subset H(\ddiv)$ and $q^n\in \mathbb{V}_h^0\subset
H^1$. Hence, $C_n$ is a Casimir for $n=1$ (total vorticity) and $n=2$
(enstrophy).

\citet{mcrae2014energy} verified these conservation properties in
numerical experiments, and showed second order convergence with $h$
for the scheme with the $BDFM_1$-$DG_1$ finite element spaces on triangles
($BDFM_1$ is a slightly more exotic variant which has an intermediate
number of degrees of freedom between $BDM_1$ and $RT_1$, which results
in a 2:1 ratio of velocity to pressure degrees of freedom).

The Poisson bracket approach has been extended to other compatible
spaces with various motivations. \citet{eldred2019quasi} used finite
element spaces built around splines to form higher order discrete de
Rham complexes. These spaces have the same degrees of freedom as the
lowest order $Q_1$-$RT_0$-$DG_0$ complex on quadrilaterals, and
achieve higher order by making use of degrees of freedom from a patch
of neighbouring cells. The advantage is that this removes the jump in
the dispersion relation for gravity waves in higher order spaces on
quadrilaterals, as discussed in Section \ref{ssec:spectral gap}. The
price to be paid is that there is increased interelement coupling, and
that there are some technicalities at the boundaries between patches
of structured quadrilaterals, e.g. at the edges and vertices of the
cube upon which a cubed sphere mesh is constructed. Since Hamiltonian
and the Poisson brackets are the same, the only thing that has changed
is the finite element spaces, which still satisfy the discrete de Rham
complex so energy and enstrophy conservation follows directly.

\citet{lee2018discrete} extended the method to mixed mimetic spectral
elements, which are a variant of mixed elements 
using spectral element histopolation functions to construct high
order spaces. The usual spectral element technique of using incomplete
quadrature then leads to diagonal mass matrices for the continuous
space $\mathbb{V}_h^0$ without losing the discrete de Rham complex
property. Again, since the Hamiltonian and Poisson brackets are the
same (excepting some details on quadrature rules, where care must be
taken), and the new finite element spaces still satisfy the discrete
de Rham complex, the energy and enstrophy conservation follows
directly. \citet{lee2018mixed} extended these spaces to the
surface of the sphere.

\subsection{Poisson integrators}
To extend these conservation properties to a fully discrete method
after time discretisation, we need to look beyond the implicit
midpoint rule into the more general case of Poisson integrators.
To make this generalisation, we write the Poisson bracket as
\begin{equation}
  \{F,G\} = A\left(\dede{F}{z}, \dede{G}{z}; z\right),
\end{equation}
where $z\in W$ comprises the dynamic fields (i.e., $z=(u,D)$ and $W=\mathbb{V}_h^1\times \mathbb{V}_h^2$ for the
case of the shallow water equations). We use this notation to
express that Poisson brackets are bilinear in $(\dede{F}{z},\dede{G}{z})$
but with possibly arbitrary additional dependence on $z$, which acts
as a coefficient. When there is no explicit dependence on $z$, we obtain
linear dynamics, and the $z$ dependency encodes nonlinear dynamics.
From the properties
of the Poisson bracket, $A$ is bilinear and antisymmetric in
$\dede{F}{z}$ and $\dede{G}{z}$. To derive one particular Poisson
integrator, we write $z(s)=z^n + s(z^n-z^{n+1})$, and seek
$z^{n+1}$ such that
\begin{equation}
  \int_0^{\Delta t}\left(\pp{}{t}{F}[z(s)] + A(\dede{F}{z}[z(s)], \dede{H}{z}[z(s)]; z^{n+1/2})\right)\beta(s)
  \diff s = 0, 
\end{equation}
for all linear functions $\beta(s)$.
In other words, we replace the bracket $\{F,G\}$ by
\begin{equation}
  \{F,G\} = A\left(\dede{F}{z}, \dede{G}{z}; z^{n+1/2}\right),
\end{equation}
where $z^{n+1/2}=(z^n + z^{n+1})/2$, and project the equation
onto linear dynamics in time. We observe energy conservation
since taking $F=H$ and $\beta=1$ leads to
\begin{align}
  \frac{H[z^{n+1}]-H[z^n]}{\Delta t} &= 
  \int_0^{\Delta t}\pp{}{s}{H}[z(s)]\diff s, \\
  &= 
  - \int_0^{\Delta t}\underbrace{
    A\left(\dede{H}{z}[z(s)], \dede{H}{z}[z(s)]; z^{n+1/2}\right)
  }_{=0}
  \diff s = 0, 
\end{align}
by antisymmetry.

This formulation leads to a practical method since taking
$F[z]=\left\langle w, z\right\rangle$ for $w\in W$ and $\beta=1$ gives
\begin{equation}
  \label{eq:poisson}
  \left\langle w, z^{n+1} - z^n \right\rangle + 
  A\left(w, \int_0^{\Delta t}\dede{H}{z}[z(s)]\diff s; z^{n+1/2}\right)  = 0,
\end{equation}
by linearity in the second argument. This scheme was introduced along
with higher order variants as a larger set of Poisson integrators in
\citet{hairer2010energy,cohen2011linear}. The easiest way to obtain
implementable formulae for the scheme is to choose a quadrature rule
for the time integral in \eqref{eq:poisson} such that the integral is
exact. This is possible whenever the Hamiltonian is polynomial.  For
example, when the Hamiltonian is quadratic, $\dede{H}{z}$ is
linear, and the midpoint rule
\begin{equation}
  \int_0^{\Delta t} \dede{H}{z}[z(s)]\diff s =
  \frac{1}{2}\left(
    \dede{H}{z}\left[z^n + \frac{1}{2}\left(z^{n+1}-z^n\right)\right]
  \right),
\end{equation}
(evaluated at $s=1/2$ with weight 1), is exact; the scheme is then
equivalent to the implicit midpoint rule.

For the rotating shallow water equation scheme described above,
the Hamiltonian is cubic, so a two point quadrature must be used
to compute the time integral involving the quadratic derivatives
of the Hamiltonian exactly.
We obtain the scheme
\begin{align}
  \nonumber
  \left\langle w, u^{n+1} \right\rangle
  + \left\langle \phi, D^{n+1}\right\rangle
  =&
  \left\langle w, u^{n} \right\rangle
  + \left\langle \phi, D^{n}\right\rangle \\
& \qquad  - \Delta t A\left((w,\phi),\left(\overline{\dede{H}{u}},\overline{\dede{H}{D}}\right);
  (u^{n+1/2}, D^{n+1/2})\right), \nonumber \\
  & \qquad
  \forall (w,\phi)\in \mathbb{V}_h^1\times \mathbb{V}_h^2,
\end{align}
where
\begin{align}
  \overline{\dede{H}{u}} & = P_1[m^{n+1/2}], \\
  &= P_1\frac{1}{3} \left(D^nu^n + \frac{1}{2}D^nu^{n+1}
  + \frac{1}{2}D^{n+1}u^n + D^{n+1}u^{n+1}\right), \\
  \overline{\dede{H}{D}} & = P_2(\pi^{n+1/2}), \\
&  = P_2 \left(\frac{1}{6}\left(|u^n|^2 +
  u^n\cdot u^{n+1} + |u^{n+1}|^2\right) + \frac{g}{2}(D^{n+1}+D^n)
  + b\right).
\end{align}
This can be implemented as $(u^{n+1},D^{n+1},q^{n+1/2},\overline{m}^{n+1/2})
\in \mathbb{V}_h^1\times \mathbb{V}_h^2\times \mathbb{V}_h^0\times \mathbb{V}_h^1$ such that
\begin{align}
  \nonumber
  \left\langle w, u^{n+1}-u^n\right\rangle
  + \Delta t \left\langle w, q^{n+1/2}(\overline{m}^{n+1/2})^{\perp} \right\rangle & \\
\qquad  - \Delta t \left\langle \nabla\cdot w, \pi^{n+1/2} \right\rangle
  & = 0, \quad \forall w \in \mathbb{V}_h^1, \\
  \left\langle \phi, D^{n+1} - D^n \right\rangle
  + \Delta t\left\langle \phi, \nabla\cdot \overline{m}^{n+1/2} \right\rangle & = 0, \quad
  \forall \phi \in \mathbb{V}_h^2, \\
  \left\langle \gamma, D^{n+1/2}q^{n+1/2}\right\rangle
  + \left\langle \nabla^{\perp}\gamma, u^{n+1/2} \right\rangle
  - \left\langle \gamma, f \right\rangle & = 0, \quad \forall \gamma \in \mathbb{V}_h^0, \\
  \left\langle v, \overline{m}^{n+1/2} - m^{n+1/2} \right\rangle & = 0,
  \quad \forall v \in \mathbb{V}_h^1,
\end{align}
with $m^{n+1/2}$, $\pi^{n+1/2}$ defined as above. This scheme is
similar, but not identical to, the implicit midpoint rule, but results
in energy conservation. \citet{cohen2011linear} showed that this class
of Poisson integrators also preserves Casimirs that are at most
quadratic functions of the state space variables. This covers mass and
total vorticity, but not enstrophy, which is a nonpolynomial function
of $u$ and $D$. Finding Poisson integrators that preserve enstrophy
for this discrete formulation of the rotating shallow water equations is
an open problem.

\subsection{Enstrophy conservation on domains with boundaries}
\label{sec:ens bdys}
Something that we have neglected from our discussion so far is the case
when the domain $\Omega$ has an exterior boundary. This is important
when extending these tools to ocean applications (where there are
coastlines). When boundaries are present, we consider the subcomplex
\begin{equation}
\begin{CD}
  \mathring{H}^1 @> \nabla^\perp >> \mathring{H}(\textrm{div}) @>
  \nabla\cdot
  >> L^2 \\
  @VV{\pi_0}V @VV{\pi_1}V @VV{\pi_2}V \\
  \mathring{\mathbb{V}}_h^0 @> \nabla^\perp >> \mathring{\mathbb{V}}_h^1 @> \nabla\cdot
  >> \mathbb{V}_h^2 \\
\end{CD}
\end{equation}
where
\begin{align}
  \mathring{H}_1 &= \{\phi \in H^1: \tr_{\partial\Omega}{\phi}=0\}, \\
  \mathring{\mathbb{V}}_h^0 &= \{\phi \in \mathbb{V}_h^0: \tr_{\partial\Omega}{\phi}=0\}, \\
  \mathring{H}(\ddiv) &= \{u \in H(\ddiv): \tr_{\partial\Omega}{u\cdot n}=0\}, \\
  \mathring{\mathbb{V}}_h^1 &= \{u \in \mathbb{V}_h^1: \tr_{\partial\Omega}{u\cdot n}=0\},
\end{align}
and where $\tr$ is the boundary trace operator. When a boundary
is present, we must modify \eqref{eq:discrete q} to incorporate
a boundary integral,
\begin{equation}
    \left\langle \gamma, D q\right\rangle =
    -\left\langle \nabla^{\perp}\gamma, u \right\rangle +
    \left\llangle \gamma, n^\perp\cdot u\right\rrangle +
    \left\langle \gamma, f\right\rangle,
  \quad \forall \gamma \in \mathbb{V}_h^0,
\end{equation}
where $\left\llangle\cdot,\cdot\right\rrangle$ is the usual $L^2$ inner product
on $\partial \Omega$. This is necessary because $q$ does not vanish on
the boundary in general. One can then try to proceed using this
modified definition of $q$ in the Poisson bracket, now defined over
$\mathring{\mathbb{V}}_1\times \mathbb{V}_h^2$. The proof of enstrophy conservation then
fails because although $\nabla^\perp q\in \mathbb{V}_h^1$, $\nabla^\perp q \notin
\mathring{\mathbb{V}}_1$ in general, because $q$ is not constant on the
boundary.  One possible solution is to restrict $q\in \mathring{\mathbb{V}}_h^0$,
and use Equation \eqref{eq:discrete q} with test functions in
$\mathring{\mathbb{V}}_h^0$, so that $\nabla^\perp q \in \mathring{\mathbb{V}}_h^1$. This
recovers enstrophy conservation, but at the expense of consistency as
it commits a first-order error in forcing $q$ to be zero on the
boundary.

An alternative solution presented in \citet{bauer2018energy}
is to split $\mathbb{V}_h^0=\mathring{\mathbb{V}}_h^0\oplus
(\mathring{\mathbb{V}}_h^0)^{\perp}$ where $(\mathring{\mathbb{V}}_h^0)^\perp$ is the
$L^2$-orthogonal complement of $\mathring{\mathbb{V}}_h^0$ in $\mathbb{V}_h^0$. We then
extend the solution space to $(u,D,Z')\in \mathring{\mathbb{V}}_h^1\times \mathbb{V}_h^2 \times
(\mathring{\mathbb{V}}_h^0)^{\perp}$, and define the bracket
\begin{align} \nonumber
  \{F,G\} & = \left\langle q,\dede{F}{u}\cdot\dede{G}{u}^\perp\right\rangle
  - \left\langle \nabla\cdot\dede{F}{u}, \dede{G}{D}\right\rangle
  + \left\langle \nabla\cdot\dede{G}{u}, \dede{F}{D}\right\rangle \\
  & \qquad +\left\langle \nabla\dede{F}{Z'}, q\dede{G}{u}\right\rangle
  -\left\langle \nabla\dede{G}{Z'}, q\dede{F}{u}\right\rangle,
\end{align}
with the same Hamiltonian, where $q\in \mathbb{V}_h^0$ such that
\begin{align}
  \label{eq:q with bc}
  \left\langle \gamma, qD \right\rangle 
  - \left\langle \gamma, f + \mathring{Z}+Z'\right\rangle & = 0, \forall \gamma \in
  \mathbb{V}_h^0, \\
  \left\langle \gamma, \mathring{Z} \right\rangle 
  - \left\langle \nabla^\perp \gamma, u \right\rangle & = 0, \quad \forall \gamma\in
  \mathring{\mathbb{V}}_h^0.
\end{align}
Here, $Z'$ represents the contribution to the vorticity from the
boundary, which would normally be given by
\begin{equation}
  \left\langle \gamma, Z'\right\rangle = -\left\langle \nabla^\perp \gamma, u\right\rangle
  + \left\llangle \gamma, n^\perp\cdot u\right\rrangle + \left\langle \gamma, f\right\rangle,
  \quad \forall \gamma \in (\mathring{\mathbb{V}}_h^0)^\perp.
\end{equation}
However, here we just ensure that $Z'$ satisfies this condition
initially, and that at future times $Z'$ has its own dynamics
consistent with the potential vorticity equation as we shall see later
(and see \ref{sec:consistent PV} for further discussion of this).
Since the Hamiltonian does not depend on $Z'$, we have
$\dede{H}{u}$ and $\dede{H}{D}$ as before, and $\dede{H}{Z'}=0$.
Thus the $u_t$ and $D_t$ equations are unchanged, and we have
\begin{equation}
  \left\langle \gamma, Z'_t \right\rangle = \left\langle \nabla\gamma, q\dede{H}{u}\right\rangle
  = 
  \left\langle \nabla \gamma,
  mq\right\rangle, \forall \gamma \in (\mathring{\mathbb{V}}_h^0)^\perp.
\end{equation}
The bracket is antisymmetric, so the Hamiltonian is conserved.

To check whether this has repaired the Casimirs, we recompute the
derivatives of $C_n$ for this extended phase space, writing
$D_\epsilon = D + \epsilon \phi$ for $\phi\in \mathbb{V}_h^2$, $u_\epsilon = u +
\epsilon v$ for $v\in \mathring{\mathbb{V}}_h^1$, $Z'_\epsilon = Z' + \epsilon
\gamma'$ for $\gamma'\in (\mathring{\mathbb{V}}_h^0)^\perp$, and defining
$q_\epsilon\in \mathbb{V}_h^0$ and $\mathring{Z}\in \mathring{\mathbb{V}}_h^0$ such that
\begin{align}
  \left\langle \gamma, q_\epsilon D_\epsilon \right\rangle 
  - \left\langle \gamma, f+Z'_\epsilon + \mathring{Z}_\epsilon\right\rangle & = 0, \forall \gamma \in
  \mathbb{V}_h^0, \\
    \left\langle \gamma, \mathring{Z}_\epsilon \right\rangle 
  - \left\langle \nabla^\perp \gamma, u_\epsilon \right\rangle & = 0, \quad \forall \gamma\in
  \mathring{\mathbb{V}}_h^0.
\end{align}
Then, for $n\in (1,2)$,
\begin{align}
  \delta C_n & = \left\langle \dede{C_n}{u}, v \right\rangle
  + \left\langle \dede{C_n}{D}, \phi\right\rangle
  + \left\langle \dede{C_n}{Z'}, \gamma'\right\rangle \\
  & = \lim_{\epsilon\to 0}
  \frac{1}{\epsilon}
  \left(
  C_n[u+\epsilon v, D + \epsilon \phi, Z'+\epsilon \gamma']- C_n[u,D]
  \right), \\
  &  = \int_\Omega \dd{}{\epsilon}|_{\epsilon=0}D_\epsilon q_\epsilon^n
  \diff x, \\ \nonumber
  & = \left\langle q_\epsilon^{n-1}|_{\epsilon=0}, \dd{}{\epsilon}|_{\epsilon=0}D_\epsilon q_\epsilon \right\rangle \\
& \qquad  + \left\langle (n-1)q_\epsilon^{n-2}|_{\epsilon=0}\dd{}{\epsilon}|_{\epsilon=0} q_\epsilon,
  D_\epsilon q_\epsilon \right\rangle, \\ \nonumber
  & = n\left\langle q_\epsilon^{n-1}|_{\epsilon=0}, \dd{}{\epsilon}|_{\epsilon=0}D_\epsilon q_\epsilon \right\rangle \\
& \qquad  - \left\langle (n-1)q_\epsilon^n|_{\epsilon=0}, \dd{}{\epsilon}|_{\epsilon=0} D_\epsilon \right\rangle, \\
  & = n\left\langle q^{n-1}, \gamma' + \dd{}{\epsilon}|_{\epsilon=0}\mathring{Z}_\epsilon \right\rangle
  - \left\langle (n-1)q^n, \phi \right\rangle, \\
  & = n\left\langle \mathring{P}_0^{\perp}q^{n-1}, \gamma'\right\rangle
  + n\left\langle \mathring{P}_0q^{n-1},
  \dd{}{\epsilon}|_{\epsilon=0}\mathring{Z}_\epsilon \right\rangle
  - \left\langle (n-1)q^n, \phi \right\rangle, \\
  & = n\left\langle \mathring{P}_0^{\perp}q^{n-1}, \gamma'\right\rangle
  + n\left\langle \nabla^\perp\mathring{P}_0q^{n-1}, v\right\rangle
  - \left\langle (n-1)q^n, \phi \right\rangle,
\end{align}
where $\mathring{P}_0^{\perp}$ is the $L^2$ projection onto the orthogonal
subspace to $\mathring{\mathbb{V}}_h^0$. Hence,
\begin{equation}
  \dede{C_n}{u} = n\nabla^\perp \mathring{P}_0q^{n-1}, \quad
  \dede{C_n}{Z'} = n\mathring{P}_0^{\perp}q^{n-1}, \quad
  \dede{C_n}{D} = -(n-1)P_2(q^n),
\end{equation}
for $n\in (1,2)$.  Then,
\begin{align}
  \nonumber \{C_n,G\} & = \left\langle q,
  n\nabla^\perp\left(\mathring{P}_0q^{n-1}\right)\cdot\dede{G}{u}^\perp\right\rangle
  \\
  \nonumber & \qquad
  - \left\langle \underbrace{\nabla\cdot n\nabla^\perp\left(\mathring{P}_0q^{n-1}\right)}_{=0},\dede{G}{D}\right\rangle
  + \left\langle \nabla\cdot\dede{G}{u},-(n-1)P_2(q^n)
  \right\rangle \\
  & \qquad +\left\langle \nabla n\mathring{P}_0^{\perp}q^{n-1}, q\dede{G}{u}\right\rangle
  -\left\langle \nabla\dede{G}{Z'}, n\nabla^\perp\left(\mathring{P}_0q^{n-1}\right)\right\rangle, \\
\nonumber  & = \left\langle q,
  n\nabla
  \underbrace{\left(\mathring{P}_0q^{n-1} + \mathring{P}_0^\perp q^{n-1}\right)}_{=q^{n-1}}
  \cdot\dede{G}{u}\right\rangle \\
  \nonumber & \qquad 
  - \left\langle \underbrace{\nabla\cdot n\nabla^\perp q^{n-1}}_{=0},\dede{G}{D}\right\rangle
  + \left\langle \nabla\cdot\dede{G}{u},-(n-1)P_2(q^n)
  \right\rangle \\
  & \qquad
  +\left\langle \underbrace{\nabla^\perp\cdot\nabla\dede{G}{Z'}}_{=0}, n\left(\mathring{P}_0q^{n-1}\right)\right\rangle, \\
  \nonumber  & = -\left\langle (n-1)\nabla q^n,
\dede{G}{u}\right\rangle
  + \left\langle \nabla\cdot\dede{G}{u},-(n-1)q^n
  \right\rangle \\
& = \left\langle (n-1)q^n,
\nabla\cdot\dede{G}{u}\right\rangle
  - \left\langle \nabla\cdot\dede{G}{u},(n-1)q^n
  \right\rangle
  = 0,
\end{align}
where we have repeatedly used $q^n\in \mathbb{V}_h^0$ for $n\in (1,2)$. The
surface integral vanished in the penultimate line since
$\mathring{P}_0q^{n-1}\in \mathring{\mathbb{V}}_h^0$ vanishes on the boundary,
and the surface integral vanished in the final line since the normal
component of $\dede{G}{u}\in \mathring{\mathbb{V}}_h^1$ vanishes on the boundary.
Hence the total vorticity $C_1$ and the enstrophy $C_2$ are both
Casimirs and hence are conserved by the dynamics.

\citet{bauer2018energy}
demonstrated through numerical experiments that this scheme produces
convergent solutions. This idea is closely related to the
approach of \citet{ketefian2009mass}, who introduced vorticity
variables at the boundary to make an energy enstrophy conserving
staggered finite difference method. 

\subsection{Thermal shallow water equations}

The thermal shallow water equations provide a useful stepping stone
between the rotating shallow water equations and three dimensional
models. This is because they incorporate an additional advected tracer,
the temperature, whilst remaining in the two dimensional setting.
Additionally, they provide an interesting reduced model for describing
some atmospheric processes, especially when further augmented with
a moisture variable, as discussed in the excellent book by Zeitlin
\citep{zeitlin2018geophysical}. The variational derivation of these
equations originates from \citet{ripa1993conservation}, and can
be placed in the framework of Euler-Poincar\'e equations by treating
the buoyancy $s$ as an additional advected quantity satisfying
\begin{equation}
  s_t + u\cdot\nabla s = 0.
\end{equation}
Then, the equations take the form
\begin{align}
  \pp{u}{t} + u\cdot \nabla u + fu^\perp + s\nabla(D+b) + \frac{D}{2}\nabla s & = 0, \\
  \pp{D}{t} + \nabla\cdot (Du) & = 0, \\
  \pp{s}{t} + u\cdot \nabla s & = 0,
\end{align}
where $u$ is the horizontal velocity and $D$ is the layer thickness as
before.  Instead of $s$, one can instead work with the bouyancy
density $S=sD$, which satisfies
\begin{equation}
  S_t + \nabla\cdot(uS) = 0.
\end{equation}
As a consequence of the variational derivation, these equations
have Poisson bracket formulations. When $S$ is used as the prognostic
variable, the equations can be obtained from the following Poisson bracket,
\begin{align}
  \{F,G\}_S &= \{F,G\}_0
  + \left\langle \dede{F}{S},\nabla\cdot\left(s\dede{G}{u}\right)\right\rangle
  - \left\langle \dede{G}{S},\nabla\cdot\left(s\dede{F}{u}\right)\right\rangle,
  \label{eq:S bracket}
\end{align}
where $\{\cdot,\cdot\}_0$ is the Poisson bracket presented in
\eqref{eq:q bracket}, combined with the Hamiltonian
\begin{equation}
  H = \int_\Omega \frac{D|u|^2}{2} + S\left(\frac{D}{2} + b\right)\diff x.
\end{equation}
This Poisson bracket formulation has Casimirs of the form
\begin{equation}
  C[D,u,S] = \int_\Omega DqA\left(\frac{S}{D}\right) + DB\left(
  \frac{S}{D}\right)\diff x,
\end{equation}
where $A$ and $B$ are arbitrary functions.

Alternatively, when $s$ is used as the prognostic variable, we have
the following Poisson bracket formulation,
\begin{align}
  \label{eq:s bracket}
  \{F,G\}_s &= \{F,G\}_0
  - \left\langle \nabla\dede{F}{s},s\dede{G}{u}\right\rangle
  + \left\langle \nabla\dede{G}{s},s\dede{F}{u}\right\rangle, \\
  H & = \int_\Omega \frac{D|u|^2}{2} + Ds\left(\frac{D}{2} + b\right)\diff x.
\end{align}
Similarly, this Poisson bracket formulation has Casimirs of the form
\begin{equation}
  C[D,u,s] = \int_\Omega DqA\left(s\right) + DB\left(
  s\right)\diff x.
\end{equation}

The advantage of the $S$ formulation is that conservation of total buoyancy
\begin{equation}
  B = \int_\Omega S \diff x,
\end{equation}
is naturally incorporated, and local conservation is possible when
choosing $S\in \mathbb{V}_h^2$. In both cases, for a conforming discretisation,
$s\in H^1$ is required, i.e. we should take $s\in \mathbb{V}_h^0$. In the case of
the $S$ formulation, following the approach to $q$ in
\citet{mcrae2014energy}, \citet{eldred2019quasi} proposed to
introduce $s\in \mathbb{V}_h^0$ as a diagnostic quantity defined by
\begin{equation}
  \label{eq:s}
  \int_\Omega \gamma D s \diff x = \int_\Omega \gamma S \diff x,
  \quad \forall s \in \mathbb{V}_h^0.
\end{equation}
Alternatively, a nonconforming discretisation can be obtained by
introducing additional facet terms into the Poisson bracket; we shall
discuss this further in section \ref{sec:upwind tracers}.

In the case of the conforming $S$ formulation with $S\in \mathbb{V}_h^2$ and prognostic 
$s\in \mathbb{V}_h^0$, the variational derivatives of $H$ with respect to
$u\in \mathbb{V}_h^1$, $D \in \mathbb{V}_h^2$ and $S\in \mathbb{V}_h^2$ become
\begin{equation}
  \dede{H}{u} = m := P_1(uD), \,
  \dede{H}{D} = P_2(|u|^2/2) + S/2, \,
  \dede{H}{S} = \frac{D}{2} + b.
\end{equation}
Picking $F=\left\langle u, w\right\rangle + \left\langle \phi, D \right\rangle
+ \left\langle S, \alpha \right\rangle$ and using the Poisson bracket \eqref{eq:S bracket}
then gives the system of equations
\begin{align}
  \left\langle w, u_t \right\rangle + \left\langle w, qm^{\perp}\right\rangle & \nonumber \\
\quad  - \left\langle \nabla \cdot w, \frac{|u|^2}{2} + \frac{S}{2} \right\rangle
  - \left\langle \frac{D}{2} + b, \nabla\cdot\left(sw\right) \right\rangle
  & = 0, \quad \forall w \in \mathbb{V}_h^1, \\
  \left\langle \phi, D_t \right\rangle + \left\langle \phi, \nabla\cdot m \right\rangle=0, \quad
  \forall \phi\in \mathbb{V}_h^2, 
  \\
  \left\langle \alpha, S_t \right\rangle + \left\langle \alpha, \nabla\cdot (sD) \right\rangle
  & = 0, \quad \forall \alpha \in \mathbb{V}_h^2, \\
  \left\langle \gamma, qD \right\rangle - \left\langle \nabla^\perp \gamma, u \right\rangle
  - \left\langle \gamma, f \right\rangle & = 0, \quad \forall \gamma \in \mathbb{V}_h^0, \\
  \left\langle \beta, sD \right\rangle - \left\langle \beta, S\right\rangle &=0, 
  \quad \forall \beta \in \mathbb{V}_h^0.
\end{align}
The introduction of the thermal variable breaks the symmetry that
makes $C_n$ a Casimir for $n>1$ for the shallow water equations.
However, straightforward calculations show that mass $M$, total
buoyancy $B$ and total vorticity $Z$ are all Casimirs for this
discrete bracket.

A similar conforming discretisation obtained from the Poisson bracket
\eqref{eq:s bracket} also preserves all three of these quantities,
with the variation that now $B$ is a quadratic functional
\begin{equation}
  B = \int_\Omega sD \diff x,
\end{equation}
so requires a time integrator that preserves quadratic Casimirs.

\citet{eldred2019quasi} introduced these formulations plus a number of
nonconforming versions with facet integrals. They introduced time
integration methods that conserve the relevant Casimirs and
demonstrated all of these schemes in convergence tests and other
benchmarks.

\subsection{Rotating compressible Euler equations}
Using the approaches described in this section, Poisson bracket
discretisations are possible for any variational fluid model with a
Hamiltonian being a function of velocity $u$, density $D$, and a
thermal field $\theta$ (as well as extensions to e.g.  magnetic flux
$B$ etc.). In this section we briefly discuss such discretisations for
the compressible Euler equations that are the basis for atmospheric
dynamical cores in weather and climate models. These equations are
given (in the ``$\theta$-$\Pi$'' formulation) by
\begin{align}
  u_t + (u\cdot\nabla)u + 2\Omega \times u + c_p\theta \nabla \Pi & =
  -g\hat{z}, \\
  \theta_t + u\cdot\nabla\theta & = 0, \\
  D_t + \nabla\cdot(uD) & = 0, \\
  \Pi^{\frac{1-\kappa}{\kappa}} & = \frac{R}{p_0}D\theta,
\end{align}
where $u$ is the velocity, $\Omega$ is the rotation vector for the
Earth, $\theta$ is the potential temperature (a scaling of temperature
that absorbs the changes in temperature due to changes in pressure),
$\Pi$ is the Exner function, $g$ is the acceleration due to gravity,
$\hat{z}$ is the unit vector pointing away from the centre of the
Earth, $D$ is the density, $\kappa=R/c_p$, $R$ is the ideal gas
constant, $c_p=R+c_v$ is the specific heat at constant pressure, $c_v$
is the specific heat at constant volume, and $p_0$ is a reference
pressure used to define $\theta$.

One Poisson bracket formulation for these equations based around the
three dimensional vorticity vector $\omega=\nabla\times u + 2\Omega$ is
\begin{align}
  \label{eq:3d poisson} \nonumber
  \{F,G\} & = \left\langle \dede{F}{u},
  \omega \times \dede{G}{u} \right\rangle
  + \left\langle \dede{F}{D}, \nabla\cdot\dede{G}{u} \right\rangle
  + \left\langle \frac{1}{D}\dede{F}{\theta}\nabla\theta, \dede{G}{u}\right\rangle \\
  & \qquad
  - \left\langle \dede{G}{D}, \nabla\cdot\dede{F}{u} \right\rangle
  - \left\langle \frac{1}{D}\dede{G}{\theta}\nabla\theta, \dede{F}{u}\right\rangle,
\end{align}
with Hamiltonian
\begin{equation}
  H = \int_\Omega \frac{D|u|^2}{2} + Dgz + \underbrace{c_vD\theta\Pi}_{\mbox{thermal energy}} \diff x,
\end{equation}
where $z$ is the height above some reference altitude. Other Poisson
bracket formulations are also possible, notably with $\Theta=D\theta$
instead of $\theta$, but we do not intend to be encyclopaedic here.

Similarly to the treatment of the thermal shallow water equations, a conforming
discretisation requires that $\theta\in H^1$, i.e. we should take
$\theta\in \mathbb{W}_h^0$.  We take $u\in \mathbb{W}_h^2$, $D\in \mathbb{W}_h^3$. The discrete
variational derivatives of $H$ are then
\begin{align}
  \dede{H}{u} &= m := P_2(Du), \\
  \dede{H}{D} &= P_3\left(\frac{1}{2}|u|^2 + gz + c_p\theta\Pi\right), \\
  \dede{H}{\theta} &= c_p P_0\left(D\Pi\right), 
\end{align}
where the derivatives of the internal energy, the third term in $H$,
require a little algebra. Using the Poisson bracket \eqref{eq:3d poisson}
as the basis for a discretisation requires an approximation
of $\nabla\times u$ since $\mathbb{W}_h^2$ is not a curl-conforming space. Analogously
to \eqref{eq:discrete q}, we approximate $\omega\in \mathbb{W}_h^1$
such that
\begin{equation}
  \left\langle v, \omega \right\rangle - \left\langle \nabla\times v, u \right\rangle = 0,
  \quad \forall v \in \mathbb{W}_h^1,
\end{equation}
for the case of a domain without boundaries. When boundaries are present,
as indeed they must be since the gravitational potential energy term $Dgz$
does not work if the domain is periodic in the vertical, then we must
again define $\mathring{\mathbb{W}}_h^1$ as
\begin{equation}
  \mathring{\mathbb{W}}_h^1 = \left\{
  \omega \in \mathbb{W}_h^1: \omega\times n = 0 \mbox{ on }\partial \Omega\right\},
\end{equation}
where $n$ is the outward pointing normal to $\Omega$. Then the state
space must be extended to include $Z'\in
(\mathring{\mathbb{W}}_h^1)^\perp$, the $L^2$-orthogonal complement to
$\mathring{\mathbb{W}}_h^1$ in $\mathbb{W}_h^1$, which represents
vorticity components on the boundary which have their own dynamics
consistent with the conservation of total vorticity. We ignore this
aspect for now, but return to it in Section \ref{sec:vorticity bdys}.

This construction leads to the discretisation
\begin{align}
\nonumber
  \left\langle w, u_t \right\rangle + \left\langle w,
  \omega \times m \right\rangle & \\
  \label{eq:euler 3D poisson}
\qquad  - \left\langle \nabla\cdot w, \frac{1}{2}|u|^2 + gz + c_p\theta\Pi
  \right\rangle
  - \left\langle w, \frac{1}{D}s\nabla\theta\right\rangle
  & = 0, \quad \forall w \in \mathbb{W}_h^2, \\
  \left\langle \phi, D_t + \nabla\cdot m \right\rangle
  & =0, \quad \forall \phi \in \mathbb{W}_h^3, \\
  \label{eq:poisson theta}
  \left\langle \gamma, \theta_t\right\rangle + \left\langle \gamma\frac{1}{D}\alpha\nabla\theta, m\right\rangle & = 0, \quad \forall \gamma \in \mathbb{W}_h^0, \\
  \left\langle r, m - uD \right\rangle & = 0, \quad \forall r \in \mathbb{W}_h^2, \\
  \left\langle v, \omega \right\rangle - \left\langle \nabla\times v, u \right\rangle
  & = 0, \quad \forall v \in \mathbb{W}_h^1, \\
  \left\langle \alpha, s - c_p D\Pi \right\rangle & = 0, \quad
  \forall \alpha \in \mathbb{W}_h^0,  \\
  \Pi^{\frac{1-\kappa}{\kappa}} & = \frac{R}{p_0}D\theta.
  \label{eq:euler 3D poisson Pi}
\end{align}
In this formulation, there is no approximation in the definition of
$\Pi$.  The equations have nonpolynomial terms due to the fractional
powers in the definition of $\Pi$, which cannot be integrated exactly.
This can be dealt with by replacing the exact integral in the
definition of the variational derivatives with a quadrature rule. Care
must be taken to use this quadrature rule consistently in all of the
terms to obtain an energy conserving formulation.  In contrast to
$\dede{H}{D}$, it is not possible to remove the projection $P_\theta$ in
the definition of $\dede{H}{\theta}$ from the equations, because
$\dede{H}{\theta}$ does not appear in an inner product with a function
from $\mathbb{V}_h^0$. Hence, we have to introduce a third auxiliary
variable $s$.  \citet{lee2021exact} used a related formulation to
build a discretisation using mimetic spectral elements.

Focussing on the approximation of the pressure gradient term
$c_p\theta\nabla \Pi$, the relevant terms are
\begin{align}
  \nonumber
  -\left\langle \nabla\cdot w, c_p\theta\Pi \right\rangle
  -\left\langle w, \frac{1}{D}s\nabla\theta\right\rangle
  &= -\left\langle \nabla\cdot (\theta w), c_p\Pi \right\rangle \\
& \qquad  + \left\langle w, \left(c_p\Pi - \frac{1}{D}s\right)
  \nabla\theta\right\rangle,
\end{align}
which is a consistent approximation to $c_p\theta\nabla \Pi$ since
$s/D$ only differs from $c_p\Pi$ by multiplication by $D$, projection
to $\mathbb{W}_h^0$ and division by $D$ again.

As discussed earlier in this article, it can be preferable to use the
temperature space $\mathbb{W}_\theta$ for $\theta$, which is more compatible
with hydrostatic balance. Since $\mathbb{W}_\theta$ allows discontinuities in
the horizontal direction, we need to modify the Poisson bracket
formulation to incorporate the nonconforming discretisation.
This is done by focussing on the following term in the Poisson bracket,
\begin{equation}
  \left\langle \frac{1}{D}\dede{F}{\theta}\nabla \theta, \dede{G}{u} \right\rangle
  = \left\langle \dede{F}{\theta}, \frac{1}{D}\dede{G}{u}\cdot\nabla \theta \right\rangle,
\end{equation}
which leads to the term in \eqref{eq:poisson theta}
approximating $u\cdot\nabla \theta$. To adapt this to the partially
discontinuous space $\mathbb{W}_\theta$, we replace with the discontinuous Galerkin
discretisation,
\begin{equation}
  L\left[\dede{G}{u},\dede{F}{\theta}; \theta\right]=
  -\left\langle \nabla_h\cdot\left(\frac{1}{D}\dede{G}{u}\dede{F}{\theta}\right), 
  \theta \right\rangle
  + \left\llangle \jump{\frac{1}{D}\dede{G}{u}\dede{F}{\theta}}, \{\theta\}
  \right\rrangle_\Gamma,
\end{equation}
using the discontinuous Galerkin notation as introduced in Section
\ref{sec:transport}. Here we have chosen a centred flux $\{\theta\}$
but will discuss upwind fluxes in \ref{sec:upwind tracers}.
As usual, this is a consistent approximation with
the facet integrals vanishing when $u$, $D$, $\theta$, etc are all
smooth functions. Then, we use $L$ in a modified Poisson bracket,
\begin{align}\nonumber
  \{F,G\} & = \left\langle \dede{F}{u},
  \omega \times \dede{G}{u} \right\rangle
  + \left\langle \dede{F}{D}, \nabla\cdot\dede{G}{u} \right\rangle
  + L\left[\dede{G}{u}, \dede{F}{\theta}; \theta\right] \\
  & \qquad 
  - \left\langle \dede{G}{D}, \nabla\cdot\dede{F}{u} \right\rangle
  - L\left[\dede{F}{u}, \dede{G}{\theta}; \theta\right].
\end{align}
This leads to the modified velocity equation,
\begin{align}
  \nonumber
    \left\langle w, u_t \right\rangle + \left\langle w,
  \omega \times m \right\rangle
  - \left\langle \nabla\cdot w, \frac{1}{2}|u|^2 + gz + c_p\theta\Pi
  \right\rangle & \\
 \quad +\left\langle \nabla_h\cdot\left(\frac{1}{D}ws\right), 
  \theta \right\rangle - \left\llangle \jump{\frac{1}{D}ws}, \{\theta\}
  \right\rrangle_\Gamma
  & = 0, \quad \forall w \in \mathbb{W}_h^2,
\end{align}
in which we recognise another consistent approximation of the pressure
gradient term. This concept for introducing any chosen $\theta$
advection scheme into a Poisson bracket formulation was introduced in
\citet{gassmann2008towards}, with application to a global finite
difference model in \citet{gassmann2013global}.

Here we should note that formulating a Poisson time integrator for this
system is challenging, because the integrands are nonpolynomial. As noted
to us by Chris Eldred, and implemented in \citet{wimmer2020energy},
a Poisson integrator must be approximated by using an incomplete quadrature
rule in the time averaged variational derivatives. This can be done to
high order at the expense of a more complicated assembly.

\subsection{Upwinding for incompressible Euler: SUPG}
\label{sec:upwinding incompressible SUPG}

One of the interesting and useful features of the Poisson bracket
formulation is that it can be modified to incorporate stabilisation of
the transport schemes whilst remaining antisymmetric and hence energy
conserving. In the context of two dimensional incompressible
turbulence (and geostrophic turbulence), this is useful because energy
cascades to large scales (where functions can be well approximated by
finite element functions), whilst enstrophy cascades to small scales
(where they can not). In energy enstrophy conserving schemes,
vorticity features pile up at the grid scale, leading to unphysical
noise, when really they should be cascading to scales below the grid
scale. Here a scheme that conserves energy whilst dissipating
enstrophy at the small scale through upwind stabilisation is
appropriate. In three dimensional isotropic turbulence, energy is also
cascading towards small scales, and so additional dissipative
gridscale closures or parameterisations are necessary. If these are
added to a Poisson bracket formulation with upwind stabilisation, we
know that there are no spurious energy transfers between scales and
between potential, kinetic and internal energy, and the only energy
changes are due to the additional dissipative closures and
parameterisations. If desired, the energy dissipated from those terms
can be collected and recycled into subgrid closures, as is done in
\citet{gassmann2013global}.

First we discuss energy conserving upwinding techniques for the
advection term in the velocity equation. For incompressible
quasigeostrophic models, \citet{sadourny1985parameterization} proposed
a subgrid closure within the Arakawa Jacobian finite difference
formulation by replacing $q\to q - \tau u\cdot\nabla u$ in the Poisson
bracket, where $\tau$ is a chosen timescale. This provides upwinding
by approximating the value of $q$ taken upstream along the streamline
passing through the gridpoint, hence the name Anticipated Potential
Vorticity Method (APVM). This was included in a rotating shallow water
formulation by \citet{arakawa1990energy}, and has been included in
more recent unstructured grid formulations in
\citet{ringler2010unified,chen2012scale}. 

When regarded as a numerical
scheme instead of a turbulence closure, this modification appears as a
$\mathcal{O}(\tau)$ consistency error. Instead, in the context of
incompressible Euler, one can replace $\omega \to \omega -
\tau(\omega_t + u\cdot \nabla \omega)$. Here, the idea is that this
term vanishes when the approximation of the solution is accurate and
smooth, since then the vorticity equation $\omega_t + u\cdot \nabla
\omega = 0$ is well approximated. Hence, the approximation is
consistent. Following this idea, we can modify the incompressible
Euler bracket \eqref{eq:omega incompressible bracket} to become
\begin{equation}
  \label{eq:omega modified}
  \{F,G\} = \int_\Omega
  \left(\omega - \tau\left(\omega_t + u\cdot\nabla\omega\right)\right)
  \dede{F}{u}\cdot \dede{G}{u}^\perp
  \diff x,
\end{equation}
where $\omega$ is obtained from \eqref{eq:omega_h} as usual.
The discretisation becomes
\begin{equation}
  \label{eq:vorticity SUPG}
  \left\langle w, u_t \right\rangle + \left\langle \left(\omega - \tau\left(\omega_t +
  u\cdot\nabla\omega\right)\right)w, u^\perp\right\rangle = 0, \quad \forall
  w\in \zeta_h.
\end{equation}
Writing $u=\nabla^\perp\psi$, $w=\nabla^\perp\phi$ for $\psi,\phi\in
\mathbb{V}_h^0$, we get
\begin{equation}
  \underbrace{\left\langle \nabla\phi, \nabla\psi_t \right\rangle}_{= -\left\langle
    \phi, \omega_t \right\rangle} + \left\langle \left(\omega -
  \tau\left(\omega_t + u\cdot\nabla\omega\right)\right)\nabla\phi,
  u\right\rangle = 0, \quad \forall \phi\in \mathbb{V}_h^0,
\end{equation}
which we rewrite as 
\begin{equation}
  \label{eq:omega SUPG}
  \left\langle \phi + \tau u\cdot\nabla\phi, \omega_t \right\rangle +
  \left\langle \phi + \tau u\cdot\nabla\phi, u\cdot\nabla\omega\right\rangle = 0, \quad \forall \phi\in \mathbb{V}_h^0,
\end{equation}
after integrating by parts in the term $-\left\langle \omega\nabla\phi,
u\right\rangle$, which is permissible since $\omega,\phi\in H^1(\Omega)$
and $u\in H(\ddiv)$.
This is the SUPG discretisation of the incompressible Euler equation,
which is obtained by replacing the test function $\phi$ with
$\phi+\tau u\cdot\nabla\phi$. The additional term leads to streamwise
diffusion of the vorticity $\omega$ without harming the consistency of
the scheme.

Since the scheme as written here is derived from a Poisson bracket
formulation, it conserves energy by construction. Regarding Casimirs
$C_n$, we now get
\begin{align}
  \left\{F,C_n\right\}
  &= -n\int_\Omega(\omega+ \tau u\cdot\nabla\omega)
  \dede{F}{u}\cdot \left(\nabla^\perp\omega^{n-1}\right)^{\perp}\diff x, \\
  &= n\int_\Omega(\omega+ \tau u\cdot\nabla\omega)
  \dede{F}{u}\cdot \nabla\omega^{n-1}\diff x, \\
  & = (n-1)\int_\Omega \dede{F}{u}\cdot \nabla\omega^n \diff x
  - n\tau\int_\Omega u\cdot\nabla\omega
  \dede{F}{u}\cdot \left(\nabla^\perp\omega^{n-1}\right)^{\perp}\diff x,
\end{align}
and this latter term only vanishes when $n=1$, hence we have
conservation of total vorticity but not enstrophy. By substituting
$\phi=\omega$ into \eqref{eq:omega SUPG}, we
can obtain the enstrophy dynamics
\begin{align}
  \dd{}{t}\int_\Omega \frac{1}{2}\omega^2 \diff x &= 
  \left\langle \omega, \omega_t \right\rangle, \\
  & = \left\langle \left(\omega -
  \tau\left(\omega_t + u\cdot\nabla\omega\right)\right)\nabla\omega,
  u\right\rangle = 0, \\
  & = \underbrace{\left\langle u, \frac{1}{2}\nabla \omega^2 \right\rangle}_{=0}
  - \tau \left\langle \omega_t, u\cdot\nabla \omega \right\rangle
  - \tau \left\langle u\cdot\nabla \omega, u\cdot \nabla \omega \right\rangle,
\end{align}
where the first term in the last line vanishes after integration by
parts and noting that $u$ is divergence free. The last term is
negative semidefinite, and corresponds to diffusion of enstrophy along
streamlines, as occurs in APVM. The middle term, which is what we get
if we change from APVM to SUPG to ensure consistency of the scheme, is
indefinite.  However, it only contains one derivative, so it is a
lower order term compared to the streamwise diffusion, and hence the
streamwise component of $\omega$ is kept smooth.

\subsection{Upwinding for incompressible Euler: vorticity free formulation}
\label{sec:upwinding incompressible natale}
An alternative Poisson bracket with upwinding for the two dimensional
incompressible Euler equations stems from the variational formulation
of \citep{natale2018variational}, discussed in the previous section.
Since that discretisation conserves energy, it should not come as a surprise
that it has a Poisson bracket formulation, given by
\begin{align}
  \{F,G\} & = -\left\langle u, \nabla^\perp\left(
  \dede{F}{u}^\perp\cdot \dede{G}{u}\right) \right\rangle
  + \left\llangle \{u\}, \jump{n^\perp\cdot\left(
  \dede{F}{u}^\perp\cdot \dede{G}{u}\right)} \right\rrangle_\Gamma, \\
  H & = \frac{1}{2}\int_\Omega |u|^2 \diff x.
\end{align}
This formulation can be thought of as an alternative way to obtain an
approximation to $\omega$ when $u\in \mathbb{V}_h^1$. Instead of using an
auxiliary equation to define $\omega$, here we integrate the curl by
parts in each cell and choose an approximation to $u$ on the facets.
This approximation is necessary because although $u$ has continuous
normal components, it does not have continuous tangential components;
the tangential component of $u$ is multivalued on the boundary. The
variational derivation leads to a centred approximation $\{u\}$, but we
can equally take an upwind approximation $\tilde{u}$ (where
$\tilde{u}$ is the value of $u$ on the upwind side of the facet).
This leads to
\begin{equation}
  \left\langle w, u_t \right\rangle - \left\langle u, \nabla^\perp
  \left(w^{\perp}\cdot u\right) \right\rangle
  + \left\llangle \{u\}, \jump{n^\perp\cdot\left(
    w^\perp\cdot u\right)} \right\rrangle_\Gamma = 0,
  \quad \forall w \in \zeta^\perp.
\end{equation}

It is more difficult to diagnose the enstrophy budget for this scheme
than it was for the SUPG scheme. However, numerical experiments in
\citet{natale2018variational} showed that this scheme does indeed tend
to reduce enstrophy whilst exactly conserving energy. They also proved
convergence of the upwinded scheme, albeit at a suboptimal rate; numerical
experiments showed convergence at the rate expected given the degree
of the polynomials (\emph{i.e.}, second order $L^2$ convergence for
$BDM_1$).

\subsection{Scale selective dissipation}
\citet{natale2017scale} investigated the ability of the schemes in
Sections \ref{sec:upwinding incompressible SUPG} and
\ref{sec:upwinding incompressible natale} to produce energy
backscatter from small to large scales consistently with two
dimensional turbulence in the forced dissipative setting.

For the scheme of \ref{sec:upwinding incompressible natale} they
provided a multiscale interpretation of the discretisation, in the
case of the $BDM_1$ space (the approach is general for $BDM$ spaces, but
discussion is simplified by just describing the lowest order
case). Since $\mathbb{V}_h^1=$$BDM_1$ contains
$\mathbb{V}_h^{1,l}=$P1$^2$ (the space of vector valued continuous
piecewise linear functions) as a subspace, we can write an $L^2$
orthogonal decomposition,
\begin{equation}
  \mathbb{V}_h^1 = \mathbb{V}_h^{1,l}\oplus \mathbb{V}_h^{1,s},
\end{equation}
where $\mathbb{V}_h^{1,s}$ is the $L^2$ orthogonal complement of $\mathbb{V}_h^{1,l}$ in $\mathbb{V}_h^1$.
Since $\mathbb{V}_h^{1,l}$ is a continuous finite element space and therefore
contains functions that are smoother than $\mathbb{V}_h^{1,s}$, we can consider
$\mathbb{V}_h^{1,l}$ to be a subspace of larger scale fields whilst $\mathbb{V}_h^{1,s}$
contains the small scales. This decomposition is not compatible with
the decomposition $\mathbb{V}_h^1=\zeta\oplus \zeta^\perp$, so we have to use the
mixed formulation where the divergence free condition is enforced
explicitly \emph{via} the pressure gradient term. The formulation may
be written as
\begin{align}
  \label{eq:hdiv decomp}
  \left\langle v, u_t \right\rangle + a(u; u, v) + s(u; u,v)  - \left\langle P, \nabla
  \cdot v \right\rangle & = 0, \quad \forall v \in \mathbb{V}_h^1, \\
  \left\langle \nabla \cdot u, \phi \right\rangle & = 0, \quad
  \forall \phi\in \mathbb{V}_h^0,
\end{align}
where
\begin{align}
  a(\hat{u}; u,v ) &= \left\langle \hat{u}^\perp, \nabla(u^\perp\cdot v) \right\rangle -
  \left\llangle\{\hat{u}^\perp\}\cdot n_+\jump{u^\perp\cdot{v}}\right\rrangle_\Gamma, \\
  s(\hat{u};u,v) & = -
  \left\llangle c_+ \jump{\hat{u}^\perp}\cdot n_+\jump{u^\perp\cdot{v}}\right\rrangle_\Gamma,
\end{align}
noting that we use $\hat{u}=u$ in \eqref{eq:hdiv decomp},
and  where $c_+$ being equal to 1 if $u\cdot n_+\geq 0$ and 0 otherwise.

If we now write $u=u^l+u^s$, with $u^l\in \mathbb{V}_h^{1,l}$ and $u^s\in \mathbb{V}_h^{1,s}$,
we observe that
\begin{align}
  s(u;u,u^l) & = -
  \left\llangle c_+ \jump{u^\perp}\cdot n_+,\jump{u^\perp\cdot{u^l}}\right\rrangle_\Gamma, \\
  & = -
  \left\llangle c_+ \jump{u^\perp}\cdot n_+,\jump{u^\perp}\cdot u^l\right\rrangle_\Gamma, \\
  & = -
  \left\llangle c_+ \jump{u^\perp\cdot n_+},\jump{u^\perp\cdot n_+}n_+\cdot u^l\right\rrangle_\Gamma, \\
  & = -
  \left\llangle c_+ \jump{u^\perp\cdot n_+}n_+,n_+\jump{u^\perp\cdot n_+}n_+\cdot u^l\right\rrangle_\Gamma, \\
  & = -\left\llangle c_+ \jump{u^\perp},\jump{u^\perp}n_+\cdot u^l\right\rrangle_\Gamma, \\
  & = -\left\llangle c_+ \jump{u},\jump{u}n_+\cdot u^l\right\rrangle_\Gamma, 
\end{align}
where we used $\jump{u^\perp}=\jump{u^\perp\cdot n_+}n_+$ since
$u^\perp$ has continuous tangential components. Consequently,
\begin{align}
  s(u;u,u^s) & = -
  \left\llangle c_+ \jump{u^\perp}\cdot n_+,\jump{u^\perp\cdot{u^s}}\right\rrangle_\Gamma, \\
  & = -
  \left\llangle c_+ \jump{u^\perp}\cdot n_+,\jump{u^\perp}\cdot (u-u^l)\right\rrangle_\Gamma, \\
  & = 
  \left\llangle c_+ \jump{u^\perp}\cdot n_+,\jump{u^\perp}\cdot u^l\right\rrangle_\Gamma
  -
  \left\llangle c_+ \jump{u^\perp}\cdot n_+,\underbrace{\jump{u^\perp\cdot u}}_{=0}
  \right\rrangle_\Gamma, \\
  & = \left\llangle c_+ \jump{u},\jump{u}n_+\cdot u^l\right\rrangle_\Gamma.
\end{align}
Since $\left\langle u^l,u^s\right\rangle = 0$, we can obtain equations for the
evolution of $E^l = \|u^l\|^2/2$ and $E^s=\|u^s\|/2$ by setting
$v=u^l$ and $v=u^s$, respectively,
\begin{align}
  \dd{E^l}{t} + a(u;u, u^l) - \left\langle p, \nabla\cdot u^l \right\rangle
  & = \left\llangle c_+ u^l\cdot n_+\jump{u}, \jump{u}\right\rrangle_\Gamma, \\
  \dd{E^s}{t} + a(u;u, u^s) - \left\langle p, \nabla\cdot u^s \right\rangle
  & = -\left\llangle c_+ u^l\cdot n_+\jump{u}, \jump{u}\right\rrangle_\Gamma.
\end{align}
Since $c_+ u^l\cdot n_+\geq 0$, the upwinding creates an energy transfer
from small to large scales.

\citet{natale2017scale} demonstrated that this leads to energy
backscatter in practice, by using the analysis technique of
\citet{thuburn2014cascades}. This technique involves simulating two
dimensional incompressible turbulence with Newtonian damping and
wavenumber 16 forcing. The instantaneous rate of change of local
energy is computed, and then Fourier transformed to obtain
$\pp{}{t}{E}(k)$, the rate of change of energy at wavenumber $E$.  Then
the same solution is filtered by removing all wavenumbers above a
cutoff $k_T$, and the rate of change of energy is recomputed using
this filtered solution, to obtain $\pp{}{t}{E}_T(k)$. Then, the rate of
change of subgrid energy $E_{SG}(k)= \pp{}{t}{E}(k) - \pp{}{t}{E}_T(k)$ is
computed. This shows the rate of change of energy at wavenumber $k$
due to wavenumbers $>k_T$. This can then be compared with a high
resolution reference solution which has a much larger range of scales
to support backscatter from. This computation is then repeated using
enstrophy instead of energy, obtaining
$Z_{SG}(k)$. \citet{natale2017scale} examined the upwind scheme of
\citet{natale2018variational}, together with the SUPG scheme
\eqref{eq:vorticity SUPG}, using this technique, as well as comparing
the upwind momentum flux formulation of \citet{guzman2017h}, which does
not conserve energy. The experiments showed that all three schemes
exhibit the trough in $Z_{SG}(k)$ near $k=k_T$, which demonstrates
that enstrophy is being transported to smaller scales, consistent with
the enstrophy cascade. However, the upwind momentum flux formulation
showed no peak at low $k$ in $E_{SG}(k)$ that is indicative of the
energy inverse cascade in the reference solution. Both the
\citet{natale2018variational} and SUPG schemes showed such a peak,
although it is stronger and closer to the reference solution for the
SUPG scheme. When the upwinding is replaced by centred approximation
in the \citet{natale2018variational} scheme, and when the $\tau$
parameter is set to zero in the SUPG scheme, no statistical
equilibrium is reached and the implicit solvers for the systems
eventually fail. Hence, we conclude that energy conservation and
some form of stabilisation by upwinding or SUPG are both critical
to obtaining these important features in two dimensional forced
dissipative turbulence.

\subsection{Upwinding for rotating shallow water equations using potential vorticity}
\label{ssec:mcrae}
Following \citet{arakawa1990energy}, \citet{mcrae2014energy}
demonstrated that an energy conserving enstrophy dissipating scheme
for the rotating shallow water equations is possible using the APVM
technique where $q$ is replaced by $q - \tau u\cdot \nabla q$ in
\eqref{eq:q bracket}, extending the APVM idea discussed above to the
rotating shallow water equations.  This was demonstrated to have a
beneficial effect on the smoothness of the solution whilst still
preserving energy. \citet{natale2018variational} proposed to replace
$q$ in \eqref{eq:q bracket} by $q - \tau((Dq)_t + \nabla\cdot(mq))/D$
(or equivalently, by $q - \tau(q_t + m\cdot \nabla q)$, since
$D_t+\nabla\cdot m=0$ in $L^2$) in order to obtain streamwise
stabilisation within a consistent scheme, since smooth solutions of
the rotating shallow water equations satisfy $(Dq)_t +
\nabla\cdot(Duq)=0$. This leads to the system
\begin{align}
  \label{eq:SUPG q u}
  \left\langle w, u_t \right\rangle + \left\langle
  \left(q - {\tau}(q_t + m\cdot\nabla q\right)
  w, m^\perp \right\rangle & \nonumber \\
  \qquad - \left\langle \nabla\cdot w, \frac{1}{2}|u|^2 + g(D+b) \right\rangle & = 0,\,
  \quad \forall w \in \mathbb{V}_h^1, \\
  \left\langle \phi, D_t \right\rangle + \left\langle \nabla\cdot m \right\rangle & = 0,\,
  \forall \phi \in \mathbb{V}_h^2, \\
  \label{eq:SUPG q q}
  \left\langle \gamma, qD \right\rangle + \left\langle \nabla^\perp \gamma, u\right\rangle
  - \left\langle \gamma, f \right\rangle & = 0, \,\forall \gamma \in \mathbb{V}_h^0, \\
  \left\langle v, m - uD \right\rangle & = 0, \,\forall v \in \mathbb{V}_h^1.
  \label{eq:SUPG q m}
\end{align}
To obtain the enstrophy dynamics, we use $w=\nabla^\perp q$ in
\eqref{eq:SUPG q u} and substitute with $\gamma=q$ in
\eqref{eq:SUPG q q} to obtain
\begin{align}
  \dd{}{t}\int_{\Omega} \frac{1}{2}q^2D\diff x
&= \left\langle
  \left(q - \tau(q_t + m\cdot \nabla q)\right),
  m\cdot \nabla q \right\rangle, \\
&= -\left\langle
  {\tau}(qD)_t
m\cdot \nabla q \right\rangle 
- \left\langle \nabla\cdot(qm)),
  m\cdot \nabla q \right\rangle, 
\end{align}
where we again have an indefinite consistency term and a streamwise
diffusion term that is always $\leq 0$. The energy conservation and
stabilisation of enstrophy, with decay of enstrophy when gridscale
features arise through vortex stretching, was demonstrated for this
upwind stabilised rotating shallow water scheme in \citet{bauer2018energy} using
numerical experiments.

The enrichment of this scheme with upwind stabilisation presents an
opportunity for a stability analysis which seems attainable at the
time of writing but is presently open. In particular, it would be
interesting to consider the stability of a backward Euler step for
these equations.

\subsection{Upwind discontinuous Galerkin methods for active tracers}
\label{sec:upwind tracers}
Inspired by the observation of \citet{gassmann2008towards} that
upwinding for advected quantities (such as layer depth, density,
temperature etc.) can be incorporated into an energy conserving scheme
by simply ensuring that the antisymmetry is maintained in the bracket,
\citet{wimmer2020benergy} and \citet{wimmer2020energy} examined energy
conserving tracer upwinding using upwind discontinuous Galerkin schemes
and SUPG schemes, respectively. \citet{wimmer2020benergy} considered
upwind discontinuous Galerkin schemes for the rotating shallow water
equations, applied to both the layer depth (which is in $\mathbb{V}_h^2$
which allows arbitrary discontinuities between cells) and the velocity (which is
in $\mathbb{V}_h^1$, and hence allows discontinuity in the tangential
components, which is sufficient to allow for some dissipation of small
scale enstrophy).

The goal is to obtain an energy conserving formulation for which the
layer depth equation takes the upwind discontinuous Galerkin
form,
\begin{equation}
  \label{eq:D dg for poisson}
  \left\langle \phi, D_t \right\rangle - \left\langle \nabla\phi, Du \right\rangle
  + \left\llangle \jump{\phi u}, \tilde{D} \right\rrangle_\Gamma = 0,
  \quad \forall \phi \in \mathbb{V}_h^2.
\end{equation}
To do this, we have to realise it as a modification of the component
of the Poisson bracket given by
\begin{equation}
  \label{eq:poisson D bit}
  \left\langle \dede{F}{D}, \nabla\cdot \dede{H}{u} \right\rangle,
\end{equation}
where $\dede{F}{u}=\phi$ and $\dede{H}{u}=P_1(Du)$. The problem is
that $u$ explicitly appears in \eqref{eq:D dg for poisson}, whilst
$\dede{H}{u}$ involves a projection of $uD$. \citet{wimmer2020benergy}
solved this problem by introducing a recovery operator
$\mathcal{U}:(D,m)\in \mathbb{V}_h^2\times \mathbb{V}_h^1 \to u\in\mathbb{V}_h^1$
defined by
\begin{equation}
  \left\langle Dv, u \right\rangle = \left\langle v, m\right\rangle, \quad \forall v \in
  \mathbb{V}_h^1,
\end{equation}
which is well defined provided that $D>0$ (breaking this condition will
cause the scheme to fail anyway). Note in particular that if $m=P_1(Dw)$
for any $w\in \mathbb{V}_h^1$, then $\mathcal{U}(D, m) = w$.

Thus, we can rewrite \eqref{eq:D dg for poisson} as
\begin{align}
  \nonumber
  \left\langle \dede{F}{D}, D_t \right\rangle \underbrace{- \left\langle \nabla\dede{F}{D}, D
  \mathcal{U}\left(D,\dede{H}{u}\right) \right\rangle
  + \left\llangle \jump{\dede{F}{D} \mathcal{U}\left(D,\dede{H}{u}\right)},
  \tilde{D} \right\rrangle_\Gamma} = 0, \\
  \quad \forall \dede{F}{D} \in \mathbb{V}_h^2,
\end{align}
with the term marked with the underbrace replacing \eqref{eq:poisson D bit}.
In order to keep the bracket antisymmetric, we also make the same substitutions
in the corresponding term with $F$ and $G$ exchanged, and the Poisson bracket
becomes
\begin{align}
  \{F,G\} = & \left\langle \dede{F}{u}, q \dede{G}{u}^\perp \right\rangle
  {- \left\langle \nabla\dede{F}{D}, D
  \mathcal{U}\left(D,\dede{G}{u}\right) \right\rangle
  + \left\llangle \jump{\dede{F}{D} \mathcal{U}\left(D,\dede{G}{u}\right)},
  \tilde{D} \right\rrangle_\Gamma} \nonumber \\
  & \quad  {+ \left\langle \nabla\dede{G}{D}, D
  \mathcal{U}\left(D,\dede{G}{u}\right) \right\rangle
  - \left\llangle \jump{\dede{G}{D} \mathcal{U}(D,\dede{F}{u})},
  \tilde{D} \right\rrangle_\Gamma}.
\end{align}

The velocity upwinding used is a development of the
\citet{natale2018variational} scheme, which does not require an auxiliary
potential vorticity variable. To achieve this, we apply the modification
in the $q$ term of the Poisson bracket,
\begin{equation}
  \left\langle \dede{F}{u}, q\dede{G}{u} \right\rangle
  \mapsto \left\langle \dede{F}{u}, \frac{1}{D}\nabla^\perp\cdot u \dede{G}{u}^\perp
  \right\rangle
  + \left\langle \dede{F}{u}, \frac{f}{D}\dede{G}{u}^\perp \right\rangle.
\end{equation}
Following a discontinuous Galerkin methodology, we then integrate by
parts separately in each cell, and select $\tilde{u}$, the value of $u$
from the upwind side, in the corresponding facet integral, to obtain
\begin{equation}
  -\left\langle u, \nabla^\perp\left(\dede{F}{u}\cdot
  \frac{1}{D}\dede{G}{u}^\perp\right) \right\rangle
  + \left\llangle \tilde{u}, \jump{n^\perp\cdot\left(
    \dede{F}{u}\cdot \frac{1}{D}\dede{G}{u}^\perp
    \right)}\right\rrangle_\Gamma.
\end{equation}
Since we are already using $\mathcal{U}$ in the layer depth terms,
we might as well avoid additional projections and use it to replace
$\dede{G}{u}/D$. Putting all of this together gives the
\begin{align}
  \{F,G\} \nonumber= &
    -\left\langle u, \nabla^\perp\left(D\mathcal{U}\left(D,\dede{F}{u}\right)\cdot
  \mathcal{U}\left(D,\dede{G}{u}\right)^\perp\right) \right\rangle \\
& \quad   + \left\llangle \tilde{u}, \jump{n^\perp\cdot\left(
    D\mathcal{U}\left(D,\dede{F}{u}\right)\cdot \mathcal{U}\left(D,\dede{G}{u}\right)^\perp
    \right)}\right\rrangle_\Gamma \nonumber
  \\
  & \quad {- \left\langle \nabla\dede{F}{D}, D
  \mathcal{U}\left(D,\dede{G}{u}\right) \right\rangle
  + \left\llangle \jump{\dede{F}{D} \mathcal{U}\left(D,\dede{G}{u}\right)},
  \tilde{D} \right\rrangle_\Gamma} \nonumber \\
  & \quad  {+ \left\langle \nabla\dede{G}{D}, D
  \mathcal{U}\left(D,\dede{F}{u}\right) \right\rangle
  - \left\llangle \jump{\dede{G}{D} \mathcal{U}\left(D,\dede{F}{u}\right)},
  \tilde{D} \right\rrangle_\Gamma}.
\end{align}
Now, if we use this Poisson bracket to generate the dynamical
equations for $D$ and $u$, we get
\begin{align}
  \left\langle w, u_t \right\rangle
  -\left\langle u, \nabla^\perp\left(D\mathcal{U}(D,w)\cdot u^\perp\right) \right\rangle
  \nonumber & \\
 \qquad + \left\llangle \tilde{u}, \jump{Dn^\perp\cdot\left(
    \mathcal{U}(D,w)\cdot u^\perp
    \right)}\right\rrangle_\Gamma \nonumber \\
  \quad  {+ \left\langle \nabla\dede{H}{D}, D
  \mathcal{U}(D,w) \right\rangle
  - \left\llangle \jump{\dede{H}{D} \mathcal{U}(D,w)},
  \tilde{D} \right\rrangle_\Gamma} & = 0, \quad \forall w \in \mathbb{V}_h^1, \\
  \left\langle \phi, D_t \right\rangle
  - \left\langle \nabla\phi, Du \right\rangle
  + \left\llangle \jump{\phi u}, \tilde{D} \right\rrangle_\Gamma & = 0,
  \quad \forall \phi \in \mathbb{V}_h^2,
\end{align}
which conserves energy by construction, despite the presence
of the upwind terms. As written, it is not clear how to
implement this scheme, due to the presence of the application
of $\mathcal{U}$ to the test function $w$. It is not scalable
to compute and store $\mathcal{U}(D,w)$ for each test function
$w$, since $\mathcal{U}(D,w)$ is globally supported in general.
\citet{wimmer2020benergy} solved this problem by introducing
an auxiliary variable $r\in \mathbb{V}_h^1$ such that
\begin{align}
  \nonumber
  \left\langle r, Dv \right\rangle &=
  -\left\langle u, \nabla^\perp\left(Dv\cdot u^\perp\right) \right\rangle \\
&\quad  + \left\llangle \tilde{u}, \jump{Dn^\perp\cdot\left(
    v\cdot u^\perp
    \right)}\right\rrangle_\Gamma, \quad \forall v \in \mathbb{V}_h^1.
\end{align}
Then,
\begin{align}
  \left\langle r, w \right\rangle & = \left\langle r, \mathcal{U}(D, w)D \right\rangle \\
  & =   -\left\langle u, \nabla^\perp\left(D\mathcal{U}(D,w)\cdot u^\perp\right) \right\rangle\nonumber \\
& \quad  + \left\llangle \tilde{u}, \jump{Dn^\perp\cdot\left(
    \mathcal{U}(D,w)\cdot u^\perp
    \right)}\right\rrangle_\Gamma, \quad \forall w \in \mathbb{V}_h^1,
\end{align}
as required. Hence, we obtain the three coupled equations for
$(u,r,D,\dede{H}{D})\in \mathbb{V}_h^1\times\mathbb{V}_h^1\times
\mathbb{V}_h^2 \times \mathbb{V}_h^2$, where
\begin{align}
  \left\langle w, u_t \right\rangle - \left\langle Dw, r \right\rangle & = 0, \forall
  w \in \mathbb{V}_h^1, \\
  \left\langle r, v\right\rangle 
  -\left\langle u, \nabla^\perp\left(Dv\cdot u^\perp\right) \right\rangle
  + \left\llangle \tilde{u}, \jump{Dn^\perp\cdot\left(
    v\cdot u^\perp
    \right)}\right\rrangle_\Gamma \nonumber \\
  \quad  {+ \left\langle \nabla\dede{H}{D}, D
  v \right\rangle
  - \left\llangle \jump{\dede{H}{D} v},
  \tilde{D} \right\rrangle_\Gamma} & = 0, \quad \forall v \in \mathbb{V}_h^1, \\
  \left\langle \phi, D_t \right\rangle
  - \left\langle \nabla\phi, Du \right\rangle
  + \left\llangle \jump{\phi u}, \tilde{D} \right\rrangle_\Gamma & = 0,
  \quad \forall \phi \in \mathbb{V}_h^2, \\
  \left\langle s,\dede{H}{D} - \frac{1}{2}|u|^2 - gD \right\rangle & = 0,
  \quad \forall s\in \mathbb{V}_h^2.
\end{align}
Note that the use of upwinding in the $D$ equation has altered the
$\dede{H}{D}$ term in the $u$ equation, and it now appears in terms
other than an inner product with $\nabla\cdot w$. This means we have
to solve for the projection to $\mathbb{V}_h^1$ in an additional
equation, which complicates the solution of this system after
discretisation with an implicit Poisson integrator (which is necessary
for exact energy conservation for the fully discrete scheme).
\citet{wimmer2020benergy} addressed the solution by using a Picard
iteration on $(u,D)$, keeping the standard linearisation about a state
of rest as an update equation and forming the nonlinear residual by
first computing $\dede{H}{D}$ and $r$ and substituting them into the
residuals for $u$ and $D$.  The linear system can then be solved for
using hybridisation techniques as described in Section
\ref{sec:examples}. Another possible approach is to apply Newton's method
to the full four component system. Then, when solving the Jacobian
linear system for the update, $\delta r$ and $\delta \dede{H}{D}$ can
be eliminated as part of a Schur complement preconditioner; the Schur
complement in $u$ and $D$ can then be approximated by the corresponding
Schur complement arising from a more standard discretisation such
as those discussed in earlier sections.

\citet{wimmer2020benergy} demonstrated that this scheme combined with
a Poisson integrator produces relative energy conservation error of
size $10^{-9}$ after 4 Picard iterations when applied to a standard
test case used in the development of numerical schemes for numerical
weather prediction (the flow over a mountain testcase 5 from
\citet{williamson1992standard}). In fact, similar sized energy
errors were even obtained when using the implicit midpoint rule (also
approximated using 4 Picard iterations), showing that the most important
aspect of energy conservation is in the spatial discretisation, at least
for this test problem. Additionally, if a midpoint rule is used then there
is no need for the additional auxiliary variable $r$ above. This is an
example where the reduction in energy error through space and time
discretisation tricks is lost due to the energy error from doing a
finite number of nonlinear iterations\footnote{Thanks to Golo Wimmer
for pointing out this observation.}.

\subsection{SUPG methods for active tracers}
\citet{wimmer2020energy} examined similar approaches when SUPG schemes
are applied to the advected quantities. The motivation for this is
that an SUPG scheme is required to stabilise the vertical transport of
temperature when it is approximated in the $\mathbb{W}_h^\theta$ space
proposed in Section \ref{sec:spaces}. SUPG schemes have
additional complications as they modify the test function
throughout the equation, including in the time derivative term.
\citet{wimmer2020energy} addressed this problem as follows. First,
we adopt the notation that SUPG makes the modification $\gamma\mapsto
\gamma + \tau S(u; \gamma)$ to the test function $\gamma\in\mathbb{W}_h^\theta$.
This general notation is to cover different possibilities; for standard
SUPG for continuous finite element spaces the modification takes the
form
\begin{equation}
  S(u;\gamma) = u\cdot \nabla\gamma.
\end{equation}
We note that $S(u;\gamma)$ will always be linear in the test function
$\gamma$, but we do not require linearity in the other
variables. Then, if the plain version of the discrete transport
equation takes the form
\begin{equation}
  \left\langle \gamma, \theta_t \right\rangle + L\left(u, \theta; \gamma\right) = 0,\quad
  \forall \gamma \in \mathbb{W}_h^\theta,
\end{equation}
then the SUPG version takes the form
\begin{equation}
  \left\langle s  +\tau S(u; s), \theta_t \right\rangle + L\left(u, \theta;
  s + \tau S(u; s)\right) = 0,\quad
  \forall s \in \mathbb{W}_h^\theta.
\end{equation}
For example, with a standard continuous Galerkin approximation, we would
have
\begin{equation}
  L(u, \theta; \gamma) = -\left\langle \nabla \gamma, u\theta \right\rangle,
  \quad \forall \gamma \in \mathbb{W}_h^\theta.
\end{equation}
We then define $s(u;\gamma) \in \mathbb{W}_h^\theta$ according to 
\begin{equation}
  \left\langle s(u; \gamma)  + \tau S(u; s(u,\gamma)), \sigma \right\rangle 
  = \left\langle \gamma, \sigma\right\rangle,\quad
  \forall \sigma \in \mathbb{W}_h^\theta.
\end{equation}
\citet{wimmer2020energy} proved the well posedness of this definition.
Our SUPG transport equation then becomes
\begin{equation}
  \left\langle \gamma, \theta_t \right\rangle + L\left(u, \theta;
  s(u; \gamma) + \tau S(u; s(u; \gamma))\right) = 0,\quad
  \forall \gamma \in \mathbb{W}_h^\theta.
\end{equation}
Following the principle of \citet{gassmann2008towards} again,
we substitute this form into the relevant terms in the Poisson bracket,
which are
\begin{align}
  \{F,G\} =& \ldots -\left\langle \frac{1}{D}\dede{G}{\theta}
  \nabla \theta, \dede{F}{u} \right\rangle 
   +\left\langle \frac{1}{D}\dede{F}{\theta}
  \nabla \theta, \dede{G}{u} \right\rangle.
\end{align}
Since we know that $\dede{H}{u}/D=u$ in the unapproximated case, 
the substitution gives
\begin{align}
  \nonumber
  \{F,G\} =& \ldots - L\left(\frac{1}{D}\dede{G}{u}, \theta;
s\left(u; \dede{F}{\theta}\right) + \tau S\left(u; s\left(u; \dede{F}{\theta}\right)\right)\right) \\
  & \quad + L\left(\frac{1}{D}\dede{F}{u}, \theta;
  s\left(u; \dede{G}{\theta}\right) + \tau S\left(u; s\left(u; \dede{G}{\theta}\right)\right)\right).
\end{align}
This produces dynamical equations of the form
\begin{align}
  \left\langle \gamma + \tau S\left(\frac{1}{D}\dede{H}{u}; \gamma\right), \theta_t \right\rangle
  + L\left(\frac{1}{D}\dede{H}{u}, \theta; \gamma + \tau S\left(\frac{1}{D}\dede{H}{u}; \gamma\right)\right) & = 0, \quad \forall \gamma \in \mathbb{W}_\theta, \\
\nonumber  \left\langle w, u_t \right\rangle + \ldots & \\
\qquad  - L\left(\frac{1}{D}w, \theta;
  s\left(u; \dede{H}{\theta}\right) + \tau S\left(u; s\left(u; \dede{H}{\theta}\right)\right)\right) &= 0,
  \quad \forall w \in \mathbb{W}_h^1,
\end{align}
where ``$\ldots$'' represents the terms coming from the other parts of the
Poisson bracket. The latter term in the $u_t$ equation is an approximation
of the term
\begin{equation}
  -c_p\theta\nabla \Pi,
\end{equation}
in the case of the compressible Euler equations. To replace
$\dede{H}{u}/D$ in the $\theta_t$ equation back with $u$ again,
\citet{wimmer2020energy} also used the $\mathcal{U}$ operator as
described above for the shallow water equations, but we will not
incorporate that additional complexity here.

\citet{wimmer2020energy} demonstrated robust upwind stabilisation
combined with energy conservation in various test problems using this
method.  For the thermal shallow water equations, the upwinded version
demonstrated much smoother temperature fields than in the standard
energy conserving version. This is significant because the thermal
shallow water equations exhibit very fine structures in the
temperature field which lead to the accumulation of numerical noise at
the gridscale, if upwinding is not used.  For the compressible Euler
equations in a vertical slice formulation in a falling bubble
configuration, they showed that the energy conserving form of the
upwinded scheme leads to the appearance of secondary Kelvin Helmholtz
vortices that appear in much higher resolution simulations of the same
problem (but do not appear with upwinding schemes on the same
resolution that do not conserve energy). This latter result seems to
suggest that the energy conserving formulation is transferring
potential energy dissipated from the upwind transport scheme and
injecting it into the kinetic energy in a manner that is consistent
with subscale processes in the higher resolution
simulation. \citet{wimmer2020energy} also developed formulae that
showed that the energy conservation is indeed maintained by the
transfer of dissipated potential energy into kinetic energy. Since the
potential energy dissipation occurs at the gridscale, this raises the
concern that the energy conservation leads to the production of noise
in the velocity field. By careful measurement of the gridscale
component of the velocity field in numerical experiments, it was shown
that the energy injection is at a larger scale (similar to what was
observed by \citet{natale2017scale}). The na\"ive intuition that the
energy conserving formulation balances upwinding diffusion with
antidiffusion in the velocity term is incorrect, as the additional
terms are not second order derivatives in velocity.

\section{Consistent vorticity and potential vorticity transport}

\label{sec:consistent PV}

Another subtopic in structure preserving schemes is schemes that have
consistent vorticity transport. This means that although potential
vorticity or vorticity is not one of the prognostic variables, the
discretised dynamics imply a discretisation for the (potential)
vorticity transport equation. This is useful because it can imply
additional control on the smoothness of the velocity field, especially
when upwinding is incorporated into the implied vorticity dynamics.
To quote \citet{ringler2010unified}: ``Given the fundamental
importance of PV in geophysical flows, numerical models are sometimes
constructed to faithfully represent some aspects of the PV dynamics
within the discrete system''.  Although this property is closely
linked with the Poisson bracket formulations described in Section
\ref{sec:poisson}, we have chosen to discuss it in a separate
section, because some papers have emphasised the importance of this
aspect whilst not strictly building in energy conservation.

The history again comes through the ``C grid'' finite difference
school of methods for numerical weather prediction, with Sadourny
providing key ideas
\citep{sadourny1972conservative,sadourny1985parameterization},
and more modern application to unstructured grids taking place in
\citep{ringler2010unified}.

For the incompressible Euler equations, we start from \eqref{eq:weak
  2d euler}. Applying the 2D curl $\nabla^\perp\cdot$ leads to
the law of conservation of vorticity,
\begin{equation}
  \label{eq:vorticity eqn}
  \omega_t + \nabla\cdot(\omega u) = 0,
\end{equation}
from which conservation of the Casimirs $C_n$ can be directly derived.
We have already seen in the previous section that the energy
conserving discretisation \eqref{eq:u zeta_h} implies a consistent
discretisation of \eqref{eq:vorticity eqn}, and that the modification
of the bracket \eqref{eq:omega modified} leads to a consistent SUPG
discretisation. It is this idea that we seek to translate to the
equations of geophysical fluid dynamics.

\localtableofcontents

\subsection{Consistent potential vorticity transport in the energy-enstrophy conserving framework}
For the rotating shallow water equations, the starting point is the
``vector invariant'' form, which is
\begin{equation}
  u_t + qu^\perp + g\nabla (D+b) = 0,
\end{equation}
where $q=(\nabla^\perp\cdot u + f)/D$ as before. Applying
$\nabla^\perp\cdot$ leads to
\begin{equation}
  \label{eq:qD}
  (qD)_t + \nabla\cdot(qDu) = 0,
\end{equation}
which is the law of conservation of potential vorticity (from which
the conservation of the Casimirs $C_n$ can also be derived directly).
Here, the identity $\nabla^\perp\cdot \nabla$ is used to eliminate the
pressure gradient term $g\nabla(D+b)$ from this conservation law, and
so compatible discretisations that preserve this identity are
important to maintain this structure.

This occurs straightforwardly in the energy-enstrophy conserving
formulation of \citep{mcrae2014energy}, presented in Section
\ref{ssec:mcrae}. Taking $w=-\nabla^\perp \gamma$ for $\gamma \in
\mathbb{V}_h^0$ in \eqref{eq:swe mcrae u} gives
\begin{equation}
  \left\langle -\nabla^\perp\gamma, u_t \right\rangle - \left\langle \nabla \gamma,
  q m \right\rangle = 0,
  \quad \forall \gamma \in \mathbb{V}_h^0.
  \label{eq:curl u}
\end{equation}
Then, taking the time derivative of \eqref{eq:swe mcrae q}
and substituting into  \eqref{eq:curl u}
gives
\begin{equation}
  \left\langle \gamma, (qD)_t \right\rangle - \left\langle \nabla \gamma, q m \right\rangle = 0,
  \quad \forall \gamma \in \mathbb{V}_h^0,
\end{equation}
which is a standard finite element discretisation of \eqref{eq:qD}. This
becomes even clearer after integrating by parts in the second term
to obtain
\begin{equation}
  \left\langle \gamma, (qD)_t + \nabla\cdot(q m) \right\rangle = 0,
  \quad \forall \gamma \in \mathbb{V}_h^0,
\end{equation}
which is simply the projection of \eqref{eq:qD} into $\mathbb{V}_h^0$.
For the compatible finite element discretisation, this integration by
parts is an identity (provided that $u\cdot n=0$ on domain boundaries)
since $\gamma, q \in \mathbb{V}_h^0\subset H^1$ and $m\in
\mathbb{V}_h^1\subset H(\ddiv)$. If we apply similar manipulations to
these in the case of the upwind stabilised system (\ref{eq:SUPG q
  u}-\ref{eq:SUPG q m}), we obtain the following stabilised
discretisation of the potential vorticity conservation law,
\begin{equation}
  \left\langle \gamma, (qD)_t \right\rangle
  - \left\langle \nabla \gamma, mq - \tau m(q_t + m\cdot\nabla q)\right\rangle = 0,
  \quad \forall \gamma \in \mathbb{V}_h^0.
\end{equation}
Some rearrangement and integration by parts (using $m\cdot n=0$ on the
boundary) then gives
\begin{equation}
  \left\langle \gamma, (qD)_t + \nabla\cdot (mq) \right\rangle
  - \left\langle \tau m\cdot \nabla \gamma, m\cdot \nabla q\right\rangle = 0,
  \quad \forall \gamma \in \mathbb{V}_h^0.
\end{equation}
As we have previously noted, $D_t+\nabla\cdot m =0$ in $L^2$, so we can
write
\begin{equation}
  \left\langle \gamma + \frac{1}{D}\tau m\cdot \nabla \gamma, (qD)_t + \nabla\cdot (mq) \right\rangle = 0,
  \quad \forall \gamma \in \mathbb{V}_h^0,
\end{equation}
which we observe is an SUPG discretisation of the law of conservation of
potential vorticity. Even though the potential vorticity is not a prognostic
variable, we obtain consistent dynamics for this diagnosed quantity.

Also in a similar direction, \citet{lee2021petrov} designed a scheme
that applies Petrov-Galerkin style upwinding by evaluating mass flux
test functions at downstream locations along advective characteristics
(similar to a semi-Lagrangian scheme). This was demonstrated to have a
beneficial effect on the implied potential vorticity dynamics
(although the scheme does dissipate energy, unlike those described
above).

In the case of the thermal shallow water equations, there is a
source term in the potential vorticity
equation. \citet{eldred2019quasi} constructed a scheme based on the
ideas above so that the diagnostic potential vorticity satisfies a
discretised version of this conservation law with sources, which
preserves constant potential vorticity in the case when entropy $s$ is
constant.

\subsection{Consistent potential vorticity transport using primal dual grids}
The previous section suggests an approach to designing numerical schemes where one
selects an advection scheme for potential vorticity (which could be
higher order accurate, with limiters etc.). Then the corresponding
scheme for velocity that is consistent with the chosen advection
scheme is deduced. This makes use of the compatible spaces, since an
equation for vorticity can be immediately obtained by choosing a test
function $w=\nabla^\perp\gamma$ for $\gamma\in \mathbb{V}_h^0$ in the
velocity equation.
\citet{thuburn2015primal} took this approach to
designing a compatible scheme based on combinations of spaces on
overlaid primal and dual grids, necessitating lowest order spaces.  To
raise the order of accuracy of the velocity advection scheme, they
chose a third-order upwind finite volume scheme for the potential
vorticity on the dual grid (a swept area scheme in this case) to
obtain higher order accurate potential vorticity dynamics in space and
time. The advantage of third order (and, generally, odd order) transport
schemes is that a backward error analysis shows that the leading order
error is diffusive rather than dispersive: this reduces grid scale
oscillations in the numerical solution. This was achieved by
considering both the potential vorticity and velocity equations after
time discretisation. This was possible because the finite volume
scheme can be reinterpreted as an equation of the form (or a discrete
time formulation that is analogous to)
\begin{equation}
  \label{eq:q upwind}
  \left\langle \gamma, (qD)_t \right\rangle - \left\langle \nabla \gamma, mq^* \right\rangle
  = 0, \quad \forall \gamma \in \mathbb{V}_h^0,
\end{equation}
for some chosen $q^*$, where $m$ is the mass flux such that $D_t +
\nabla\cdot m=0$ in $L^2$ for the corresponding discrete scheme for
$D$. Then, this equation can be obtained from the following velocity
equation
\begin{equation}
  \left\langle w, u_t \right\rangle + \left\langle w, q^*m^{\perp} \right\rangle
  + \left\langle \nabla\cdot w, P \right\rangle = 0, \quad \forall w\in\mathbb{V}_h^1,
\end{equation}
for some $P$. This can be checked upon taking $w=\nabla^\perp\gamma$.

\subsection{Consistent potential vorticity transport using Taylor-Galerkin schemes}
In a similar direction, \citet{shipton2018higher} used a third order
Taylor-Galerkin scheme for the diagnostic potential vorticity equation
and constructed the prognostic velocity equation accordingly. A
Taylor-Galerkin scheme is an extension of the Lax-Wendroff technique
in which one expands a Taylor series in time, transforming higher
order time derivatives into space derivatives using the advection
equation, before discretising in space to obtain stable schemes.

\subsection{Preservation of constant potential vorticity}
Another aspect of these schemes is that although we wish to solve
the equation in conservative form \eqref{eq:q upwind} to conserve
total vorticity, the potential vorticity also solves
\begin{equation}
  \label{eq:q adv}
  q_t + u\cdot \nabla q = 0,
\end{equation}
obtained by combining \eqref{eq:qD} with the continuity equation for
$D$. This means that the value of $q$ is preserved along
characteristics moving at speed $u$. In particular, it can be
desirable that if $q$ is constant then it remains constant, which is a
property of \eqref{eq:q adv}. To obtain this at the discrete level, we
need a mass flux $m$ such that $D_t + \nabla\cdot m=0$. This is
straightforward when the scheme \eqref{eq:swe mcrae D} is used, since
one can take $\phi=D_t+\nabla\cdot m$ implying that $D_t+\nabla\cdot
m=0$ in $L^2$ as we have previously discussed. This becomes more
complicated when discontinuous Galerkin methods are used (or finite
volume methods for lowest order spaces as used in
\citet{thuburn2015primal}) because the equation is not immediately in
that form. However, such a form can be deduced using compatible
properties of the spaces.  \citet{thuburn2015primal} used such an
approach coming from finite volume schemes \citep{ringler2010unified},
which was translated to discontinuous Galerkin methods by
\citet{shipton2018higher}.  To describe this, we consider an upwind
discontinuous Galerkin method for $D$,
\begin{equation}
  \left\langle \phi, D_t \right\rangle - \left\langle \nabla \phi, Du \right\rangle
  + \left\llangle \jump{\phi u}, \tilde{D} \right\rrangle_\Gamma = 0, \quad
  \forall \phi \in \mathbb{V}_h^2.
\end{equation}
For each cell $K$, we then define $m\in \mathbb{V}_h^1(K)$ from
\begin{align}
  \int_f \gamma (m-u\tilde{D})\cdot n \diff S & = 0,  \quad
  \forall \gamma \in T(f), \ \forall f \in K, \\
  \int_K w \cdot (m-uD)\diff x & = 0, \quad \forall w \in \mathbb{V}_h^{1,-}(K),
\end{align}
where $f$ are all the facets of $K$, $T(f)$ is the appropriate trace
space on $f$ spanned by $u\cdot n|_f$ with $u \in \mathbb{V}_h^1(K)$,
and $\mathbb{V}_h^{1,-}(K)$ is the appropriately sized curl conforming
space to close the system, as used in the definition of the commuting
projection into $\mathbb{V}_h^1$.  For example, if $\mathbb{V}_h^1$ is
$BDM_k$ then $\mathbb{V}_h^{1,-}$ is the (rotated) Raviart Thomas space of
degree $k-1$. This is a local projection that can be evaluated
independently in each cell $K$, with $m\in\mathbb{V}_h^1$ (since the
normal components agree on facets).  Then we have
\begin{align}
  \left\langle \phi, \nabla\cdot m \right\rangle
  &= -\left\langle \nabla \phi, m \right\rangle
  + \left\llangle \jump{\phi}, m \right\rrangle_\Gamma, \\
  &= -\left\langle \nabla \phi, uD \right\rangle
  + \left\llangle \jump{\phi}, u\tilde{D} \right\rrangle_\Gamma,
  \quad \forall \phi \in \mathbb{V}_h^2,
\end{align}
where we used that $\phi \in \mathbb{V}_h^2\implies \nabla\phi \in
\mathbb{V}_h^{1,-}$ and $\phi|_f \in T(f)$, following standard
calculations defining the commuting projection. Hence, we obtain
\eqref{eq:swe mcrae D} with $m$ defined as above, and so
$D_t+\nabla\cdot m = 0$ in $L^2$. We then aim to construct our scheme
to have the form \eqref{eq:q upwind}, with the property that
$q=1\implies q^*=1$.  Then, if $q=1$, we have
\begin{align}
  \left\langle \gamma, Dq_t \right\rangle & = - \left\langle \gamma, D_t \right\rangle +
  \left\langle \nabla \phi, m \right\rangle, \\ & = -\left\langle \gamma, D_t +
  \nabla\cdot m \right\rangle = 0, \quad \forall \gamma \in \mathbb{V}_h^0.
\end{align}
In other words, we have that $q$ constant implies that $q$ stays
constant. Both \citet{thuburn2015primal,shipton2018higher} implemented
this approach in a semi implicit framework by writing the timestepping
scheme as a fixed number of Picard iterations, taking care to show
that the properties above are enforced at each Picard iteration.

\subsection{Consistent potential vorticity transport with boundaries}
\label{sec:pv bdys}
Earlier we discussed the extension in \citet{bauer2018energy} of the
scheme of \citet{mcrae2014energy} to the case of the presence of
boundaries. This extension also solves the problem of how to obtain a
scheme with a consistent discretisation for the implied potential
vorticity equation in that case, including the use of SUPG
stabilisation there.  Taking the time derivative of \eqref{eq:q with
  bc} gives
\begin{align}
  \left\langle \gamma, (qD)_t \right\rangle
  & = \left\langle \gamma,
  \mathring{Z}_t + Z'_t \right\rangle, \\
  & = \left\langle \mathring{P}_h^1\gamma, \mathring{Z}_t \right\rangle + \left\langle
  \mathring{P}^{\perp}_0\gamma,
  Z'_t \right\rangle, \\
  & = \left\langle \nabla^\perp\mathring{P}_h^1\gamma, u_t \right\rangle + \left\langle
 \mathring{P}^{\perp}_0\gamma,
 Z'_t \right\rangle, \\  
 & = -\left\langle \nabla^\perp\mathring{P}_h^1\gamma, qm^{\perp} \right\rangle - \left\langle
 \nabla \mathring{P}^{\perp}_0\gamma, qm \right\rangle, \\
 & = -\left\langle \nabla\gamma, qm\right\rangle, \quad \forall \gamma \in \mathbb{V}_h^0,
\end{align}
as required.
Similarly, the equations can be modified so that $q$ is replaced by
$q^*$, the SUPG modified potential vorticity, to obtain a consistent
SUPG stabilised diagnostic potential vorticity equation (but we do not
discuss it here). The existence of this implied potential vorticity
equation actually also provides a useful equivalent formulation that
avoids explicit computation of $Z'$ (which is in the rather cumbersome
space $\mathring{P}^\perp_0$ which is not efficient to compute with).
Rather than separately incrementing $Z'$, we can redundantly increment
$q$ on the whole domain, solving
\begin{align}
  \left\langle w, u_t \right\rangle + \left\langle qw, m^\perp \right\rangle
  - \left\langle \nabla\cdot w, \frac{1}{2}|u|^2 + g(D+b) \right\rangle & = 0,\,
  \quad \forall w \in \mathring{\mathbb{V}}_h^1, \\
  \left\langle \phi, D_t \right\rangle + \left\langle \nabla\cdot m \right\rangle & = 0,\,
  \forall \phi \in \mathbb{V}_h^2, \\
  \left\langle \gamma, (qD)_t \right\rangle + \left\langle \nabla \gamma, qm\right\rangle
& = 0, \,\forall \gamma \in \mathbb{V}_h^0,
\end{align}
having initialised $q$ from \eqref{eq:discrete q}. This provides a
computationally feasible technique for a scheme with consistent
potential vorticity dynamics (and indeed a scheme that conserves
energy and even enstrophy if the SUPG form is not used). In a
practical implementation, if errors from roundoff or truncated
numerical solvers cause $q$ and $u$ to diverge from
\begin{equation}
  \left\langle \gamma, qD \right\rangle + \left\langle \nabla^\perp \gamma, u \right\rangle
  - \left\langle \gamma, f \right\rangle = 0, \, \forall \gamma \in \mathring{\mathbb{V}}_h^0
\end{equation}
at any point, then we may replace $q\gets q'+q$, where
$q'\in \mathring{\mathbb{V}}_h^0$ satisfies
\begin{equation}
  \left\langle \gamma (q'+q)D \right\rangle + \left\langle \nabla^\perp \gamma, u \right\rangle
  - \left\langle \gamma, f \right\rangle = 0, \, \forall \gamma \in \mathring{\mathbb{V}}_h^0.
\end{equation}
This is equivalent to solving for the new $q$ from \eqref{eq:discrete
  q} with Dirichlet boundary condition obtained from the old $q$.

\subsection{Consistent vorticity transport in three dimensions}
\label{sec:vorticity bdys}
Following a similar route, we can show that equations (\ref{eq:euler
  3D poisson}-\ref{eq:euler 3D poisson Pi}) have a consistent
discretisation of vorticity transport. Initially, to see this, assume
that there are no boundaries. Then we choose $w=\nabla\times \Sigma$
in \ref{eq:euler 3D poisson Pi} with $\Sigma\in \mathbb{W}_h^1$.  This
gives
\begin{align}
  \label{eq:vorticity sigma}
  \left\langle \Sigma, \omega_t \right\rangle + \left\langle \nabla\times \Sigma,
  \omega \times m \right\rangle
  - \left\langle \nabla\times \Sigma, \frac{s}{D}\nabla\theta\right\rangle
  & = 0, \quad \forall \Sigma \in \mathbb{W}_h^2,
\end{align}
where either $s=c_pD\Pi$ or $s=\dede{H}{\theta}$, depending on whether
the Poisson bracket formulation is used or not.
\eqref{eq:vorticity sigma} is an integral form of the vorticity
equation
\begin{equation}
  \label{eq:vorticity}
  \omega_t + \nabla\times (\omega\times m) + \nabla\left(
  \frac{s}{D}\right)\times \nabla\theta = 0.
\end{equation}
The last term on the left hand side of \eqref{eq:vorticity} is known
as the baroclinic torque. If we dot \ref{eq:vorticity} with $\nabla\theta$
and use \eqref{eq:theta transport}, we obtain the potential
vorticity conservation law
\begin{equation}
  (D q)_t + \nabla\cdot (qm) = 0,
\end{equation}
where $q$ is the Ertel potential vorticity,
\begin{equation}
  q = \frac{\omega \cdot \nabla \theta}{D}.
\end{equation}
It would be wonderful to have a discretisation that has a consistent
conservation of Ertel's potential vorticity, in the manner of
this section. This could be done if one could choose test functions
$\Sigma = \psi\nabla\theta$ in \eqref{eq:vorticity sigma} with
$\psi\in \mathbb{W}_h^0$, with $\theta$ also satisfying an exact
advection equation $\theta_t + u\cdot \nabla\theta=0$ in $L^2$.
Neither of these properties appear to be possible in the present
framework, and this remains a challenging unsolved problem. What is
possible, if $\theta\in \mathbb{W}_h^0$, is to obtain conservation
of total potential vorticity,
\begin{align}
  \dd{}{t}\int_\Omega D q \diff x & =  \dd{}{t}\langle \nabla\theta,
  \omega \rangle, \\
  & = \langle \nabla \theta_t, \omega \rangle
  + \langle \nabla\theta,  \omega_t \rangle, \\
  & = -\langle \theta_t, \underbrace{\delta_1\omega}_{=\delta_1\delta_2u=0} \rangle
  + \langle \nabla\theta,  \omega_t \rangle, \\
  & = -\langle \underbrace{\nabla\times \nabla \theta}_{=0}, \omega \times m
  - \frac{s}{D}\nabla\theta \rangle,
\end{align}
upon choosing $\Sigma=\nabla\theta$. Note that this also works for any
quantity $\nabla\psi$ for $\psi\in \mathbb{W}_h^0$, dynamical or not.

Returning to the conservation of vorticity in \eqref{eq:vorticity
  sigma}, we can develop
upwind stabilisations of the vorticity equation using a residual based
approach (similar to that of \citet{bendall2023improving} by 
replacing
\begin{equation}
  \left\langle \nabla\times \Sigma,
  \omega \times m \right\rangle
  \mapsto
  \left\langle \nabla\times \Sigma,
  \omega^* \times m \right\rangle,
\end{equation}
where
\begin{equation}
  \omega^* = \omega -
  \tau\left(\omega_t + \nabla\times (\omega\times m) + \nabla\left(
    \frac{s}{D}\right)\times \nabla\theta\right),
\end{equation}
for a stabilisation parameter $\tau$.

This approach can also be extended to the case of domains with
boundaries, which of course is important for atmosphere models that
have a top and bottom surface, where we assume that the boundary
condition is $u\cdot n=0$. This requires that
$u\in\mathring{\mathbb{W}_h^2}$, and we have a similar difficulty
as with the shallow water equations with defining vorticity,
since a consistent approximation requires
\begin{equation}
  \label{eq:omega bdys}
  \langle \Sigma, \omega \rangle = \langle \Sigma, 2\Omega\rangle
  -\langle \nabla\times \Sigma,
  u \rangle + \llangle n\times \Sigma, u\rrangle,
  \quad \forall \Sigma\in\mathbb{W}_h^1,
\end{equation}
and $\omega\in \mathbb{W}_h^1$, not $\mathring{\mathbb{W}}_h^1$.  We
have the same solution to the difficulty, which is to introduce the
space $(\mathring{\mathbb{W}}_h^1)^\perp$, the $L^2$ complement of
$\mathring{\mathbb{W}}_h^1$ in $\mathbb{W}_h^1$, and writing
$\omega=\omega' + \mathring{\omega}$, with
$\omega'\in\mathring{\mathbb{W}}_h^1$ and $\mathring{\omega}\in
\mathring{\mathbb{W}}_h^1$. Then, $\omega'$ has its own dynamics
defined by
\begin{equation}
  \langle \Sigma, \omega'_t \rangle + \langle \nabla\times \Sigma,
  \omega\times m \rangle - \langle \nabla\times \Sigma,
  \frac{1}{D}s\nabla\theta \rangle = 0, \quad \forall \Sigma'\in
  (\mathbb{W}_h^1)^\perp.
\end{equation}
In the Poisson bracket setting, the bracket \eqref{eq:3d poisson}
is then extended as
\begin{align}
  \{F,G\} = \ldots & + \left\langle \nabla\times \dede{F}{\omega'},
  \omega \times \dede{G}{u} -   \frac{1}{D}\dede{G}{\theta}
  \nabla\theta \right\rangle \nonumber \\
  & - \left\langle \nabla\times \dede{G}{\omega'},
  \omega \times \dede{F}{u} -   \frac{1}{D}\dede{F}{\theta}
  \nabla\theta \right\rangle.
\end{align}
Similar calculations to those in Sections \ref{sec:ens bdys} and
\ref{sec:pv bdys} then lead from this Poisson bracket to
\eqref{eq:vorticity sigma}. Further, just as in Sections \ref{sec:ens
  bdys} and \ref{sec:pv bdys}, solving the resulting equation set is
equivalent to solving \eqref{eq:vorticity sigma} in place of
\eqref{eq:omega bdys}. This type of scheme is as yet not explored in
numerical computations.

\section{Structure preserving schemes on non-affine meshes}
\label{sec:nonaffine}
The interpretation of the compatible finite element discrete de Rham
complexes and the finite element exterior calculus (FEEC) is well
known. \citet{arnold2006finite,arnold2010finite,arnold2018finite}
provided a comprehensive unifying treatment of the stability and error
analysis of numerical approximations of the Hodge Laplacian, which is
now finding applications in the design of stable discretisations for
the Stokes equation and in elasticity. In applications to geophysical
fluid dynamics, one of the main applications for finite element
exterior calculus has been in establishing a clean separation between
topological and geometric aspects of the formulation. In particular,
when nonaffine meshes are used, at first sight it seems that it is not
possible to perform exact integration when assembling the discrete
operators on a computer, because the integrands are nonpolynomial.
For example, when a function $u\in \mathbb{V}_h^1$ is transformed back to a
reference cell, a contravariant Piola transform must be used, so it
will take the form $J^{-T}\hat{u}/\det(J)$, where $J$ is the Jacobian of
the transformation from the reference cell to the mesh cell.  $J$ is
nonconstant for nonaffine meshes, and so $1/\det(J)$ is nonpolynomial.
In general, this creates potential problems because we rely on clean
separation in the Helmholtz decomposition of $\mathbb{V}_h^1$ between
divergence-free and rotational components to obtain good long time
behaviour for all methods (not just structure preserving ones). These
type of errors have the potential to break structure preserving
properties of all the schemes discussed in this section. However, all
is not lost because many of the terms in our fluid dynamics equations
result in cancellation of geometric factors (factors involving $J$
when the equations are transformed back to the reference cell). This
cancellation of geometric factors can be derived using standard vector
calculus but they are most transparently established under the
invariance of various operations involving differential forms under
pullback (wedge product, exterior derivative, etc.). In fact, the
author of this review only became aware of the possibility of some of
them after computing with the differential form formulations.
In particular, we have the following formulae,
\begin{align}
  \label{eq:magic 1}
  \int_{\hat{K}} \hat{\phi} \hat{u}\cdot \hat{w}^\perp \diff x
  & = \int_K \phi u \cdot w^\perp \diff x, \\
  \int_{\hat{K}} \hat{\phi} \nabla \cdot \hat{w} \diff x
  & = \int_K \phi \nabla\cdot w \diff x, \\
  \int_{\hat{f}} \hat{\phi}\hat{u}\cdot \hat{n} \diff S
  & =
  \int_f \phi u \cdot n \diff S,
  \label{eq:magic 3}
\end{align}
where $g_K:\hat{K}\to K$ is the mapping from reference cell $\hat{K}$
to a mesh cell $K$, $\hat{f}$ is a facet of cell $\hat{K}$ with
normal $\hat{n}$, $f$ is the image of $\hat{f}$ under $g$ with normal
$n$, and 
\begin{equation}
  \hat{\phi} = \phi\circ g, \, \frac{J\hat{w}}{\det(J)} = w \circ g,
  \, \frac{J\hat{u}}{\det(J)} = u \circ g.
\end{equation}
Following computations in \citet{thuburn2012framework} for the finite
difference case (interpreted as discrete exterior calculus (DEC)),
\citet{cotter2014finite} presented a finite element exterior calculus
formulation in the setting of the family of methods related to
\citet{mcrae2014energy}. Relevant to the discussion of this section,
they noted the existence of these Jacobian free pullback formulae as
derived from the properties of pullbacks of differential forms.
\citet{eldred2022interpretation} provided further insight into this
family of DEC and FEEC schemes for rotating shallow water equations.
They noted that a special case of the Leibniz rule for the wedge
product, in which one of the two terms is a constant, is underpinning
the presence of an implied conserved potential vorticity as described in
Section \ref{sec:consistent PV}, and hence the conservation of total
vorticity and potential vorticity for associated Poisson bracket schemes.

In the context of the Poisson bracket formulation for rotating shallow
water equations, the pullback formulae (\ref{eq:magic 1}-\ref{eq:magic
  3}) can be used to maintain a structure preserving formulation.
This is achieved by replacing the $L^2$ inner product by a quadrature
rule,
\begin{equation}
  \left\langle \phi, p \right\rangle_q = \sum_i \phi(x_i)q(x_i)w_i,
\end{equation}
and similar for vector valued functions. In practice, this quadrature
is defined cellwise as usual, as the image of quadrature points on the
reference cell $\hat{K}$ under the map $g$. This incomplete quadrature
causes a ``variational crime'' and then analysis is required to
demonstrate whether the convergence rate is affected; a minimal
condition is that it still satifies the definition of an inner product
on the relevant finite element spaces. This modification produces
a modified Helmholtz decomposition,
\begin{equation}
  \mathbb{V}_h^1 = B^1_h\oplus\tilde{\mathfrak{h}}^1_h\oplus
  (\tilde{B}^*)^1_h,
\end{equation}
where $\tilde{\mathfrak{h}}^k_h$ and $(\tilde{B}^*)^k_h$ are modified
spaces constructed using $\tilde{\delta}_k$, the dual operator defined
using the modified inner product, \emph{i.e.},
\begin{equation}
  \left\langle \phi, \tilde{\delta}^k_hu\right\rangle_q
  = -\left\langle d^{k+1}\phi, u \right\rangle_q.
\end{equation}
It also produces a modified definition of the variational derivative,
\begin{equation}
  \left\langle \dede{F}{u}, v\right\rangle_q = \lim_{\epsilon\to 0}
  \frac{1}{\epsilon}\left(F[u+\epsilon v] - F[u]\right).
\end{equation}
However, we maintain the usual $L^2$ inner product in \eqref{eq:q
  bracket}, since all of the terms in the bracket have polynomial
integrands when transformed back to the reference cell, and hence can
be evaluated exactly. Further, the integrands in \eqref{eq:discrete q}
can also be integrated exactly for similar reasons. If we assume that
the modified inner product is exact for these terms as well (it just
needs to be a sufficiently high order Gaussian quadrature) then the
formulation (\ref{eq:swe mcrae u}-\ref{eq:swe mcrae u}) is still
energy and enstrophy preserving by the above arguments. This
modification of the inner product can be extended to more complicated
formulations involving temperature and in three dimensions.

\section{Summary and outlook}
\label{sec:summary}
In this survey, we have introduced the application of compatible
finite element methods to the world of geophysical fluid dynamics,
with applications to oceans, weather and climate. We have introduced
the main properties of the spaces and their application to
understanding the discrete wave propagation properties when they are
used for linearised models. We have discussed how to build compatible
finite element methods for nonlinear models, focussing on the
transport and pressure gradient terms; we have also discussed the
approach to solving the linear and nonlinear systems that arise from
certain timestepping schemes. Then, we have surveyed the use of
compatible finite element methods in structure preserving methods:
variational integrators, Poisson integrators and schemes with
consistent potential vorticity transport. There is much more work to
be done in the analysis of all of these schemes, considering
stability, convergence of solutions, and mesh independence of
preconditioners, \emph{etc.} There are also plenty of research
directions in finding practical approaches that incorporate as much of
the structure preserving properties as possible. The finite element
exterior calculus continues to be an important guiding principle for
designing compatible finite element methods for geophysical fluid
dynamics. It should prove a useful tool for the rigorous analysis of
stability of these methods for fully nonlinear systems, where only
limited progress has been made so far.  The author looks forward to
many fruitful collaborations on compatible finite element methods
for geophysical fluid dynamics in the future.

\paragraph{Acknowledgements} The author would like to thank
Werner Bauer, Thomas Bendall, Darryl Holm, Ruiao Hu, Andrea Natale,
Oliver Street and Golo Wimmer for their very useful feedback and
discussions on draft versions of this article. The author would also
like to thank David Ham and the Firedrake Team, whose work facilitated
many of the numerical results surveyed in this article. The author is
grateful for funding of research described in this article from the
Engineering and Physical Sciences Research Council, the Natural
Environment Research Council, UK Research and Innovation, The Grantham
Institute for Climate Change, the Met Office, and Imperial College
London.

\clearpage 
\addcontentsline{toc}{section}{References}
\bibliography{cjc-acta}
\label{lastpage}
\end{document}